\documentclass{amsart}

\RequirePackage{amsmath}
\RequirePackage{amssymb}
\RequirePackage{amsxtra}
\RequirePackage{amsfonts}
\RequirePackage{latexsym}
\RequirePackage{euscript}
\RequirePackage{amscd}
\RequirePackage{amsthm}
\RequirePackage{xypic}
\xyoption{all}
\CompileMatrices

\theoremstyle{plain}
\newtheorem{thm}{Theorem}[subsection]
\newtheorem{prop}[thm]{Proposition}
\newtheorem{lemma}[thm]{Lemma}
\newtheorem{cor}[thm]{Corollary}
\newtheorem{conj}[thm]{Conjecture}
\newtheorem{mthm}{Theorem}
\newtheorem{mcor}{Corollary}
\newtheorem*{mconj}{Mod $p$ Main Conjecture}

\theoremstyle{definition}
\newtheorem{definition}[thm]{Definition}
\newtheorem{example}[thm]{Example}
\newtheorem{remark}[thm]{Remark}

\numberwithin{equation}{section}

\renewcommand{\SS}{\mathcal{S}}
\renewcommand{\O}{\mathcal{O}}
\renewcommand{\a}{\mathfrak{a}}
\renewcommand{\b}{\mathfrak{b}}
\newcommand{\p}{\wp}
\newcommand{\m}{\mathfrak{m}}

\newcommand{\Q}{\mathbf{Q}}
\newcommand{\T}{\mathbf{T}}
\newcommand{\Z}{\mathbf{Z}}

\newcommand{\inj}{\hookrightarrow}

\newcommand{\too}{\longrightarrow}

\newcommand{\F}{\mathbf{F}}
\newcommand{\Fp}{\F_p}

\newcommand{\Fpbar}{\bar{\F}_p}

\newcommand{\K}{\mathcal K}
\renewcommand{\L}{\mathcal L}

\newcommand{\Lbar}{\bar{L}}
\renewcommand{\P}{\mathbf{P}}
\newcommand{\Qbar}{\bar{\Q}}

\newcommand{\Qp}{\Q_{p}}

\newcommand{\Zp}{\Z_{p}}

\newcommand{\w}{\omega}

\newcommand{\Iw}{\Lambda}
\newcommand{\BbbIw}{\Z_p[[\Z_p^{\times}]]}
\newcommand{\prodNp}{\Z_{p,N}^{\times}}
\def\BigHecke{{\mathcal H}}

\newcommand{\Afi}{A_{f,i}}
\newcommand{\AfOi}{A_{{f_0},i}}

\newcommand{\HfS}{H_{f,\SS}}
\newcommand{\HfSf}{H_{f,\SS(f,i)}}
\newcommand{\HfSm}{H_{f,\SS_{\min}}}
\newcommand{\Hs}{H_{s}}
\newcommand{\HsS}{H_{s,\SS}}
\newcommand{\HsSf}{H_{s,\SS(f,i)}}

\newcommand{\Qi}{\Q_{\infty}}
\newcommand{\Qinf}{\Q_{\infty}}
\newcommand{\Qiv}{\Q_{\infty,v}}
\newcommand{\Qivp}{\Q_{\infty,v_{p}}}
\newcommand{\QS}{\Q_{\Sigma}}

\renewcommand{\a}{\frak a}

\newcommand{\Gl}{G_{\ell}}
\newcommand{\Gp}{G_{p}}
\newcommand{\GQ}{G_{\Q}}
\newcommand{\GQi}{G_{\Q_{\infty}}}
\newcommand{\Gv}{G_{v}}
\newcommand{\Il}{I_{\ell}}

\newcommand{\Iv}{I_{v}}
\newcommand{\Ivp}{I_{v_{p}}}

\newcommand{\Ab}{\bar{A}}

\newcommand{\rhob}{\bar{\rho}}
\newcommand{\rhobar}{\bar{\rho}}

\newcommand{\vep}{\varepsilon}

\newcommand{\Hida}{\mathcal H}

\DeclareMathOperator{\alg}{alg}

\DeclareMathOperator{\an}{an}

\DeclareMathOperator{\chr}{char}

\DeclareMathOperator{\cond}{cond}

\DeclareMathOperator{\End}{End}

\DeclareMathOperator{\Frob}{Frob}
\DeclareMathOperator{\Gal}{Gal}
\DeclareMathOperator{\GL}{GL}
\DeclareMathOperator{\Hom}{Hom}
\DeclareMathOperator{\Id}{Id}
\DeclareMathOperator{\im}{im}

\DeclareMathOperator{\new}{new}
\DeclareMathOperator{\old}{old}
\DeclareMathOperator{\ord}{ord}
\DeclareMathOperator{\nord}{n.ord}

\DeclareMathOperator{\ps}{ps}

\DeclareMathOperator{\Sel}{Sel}

\DeclareMathOperator{\Spec}{Spec}
\DeclareMathOperator{\Spf}{Spf}
\DeclareMathOperator{\spe}{sp}
\DeclareMathOperator{\sss}{ss}
\DeclareMathOperator{\semisimp}{ss}

\DeclareMathOperator{\Trace}{Trace}

\newcommand{\dalg}{\delta}
\newcommand{\lalg}{\lambda^{\alg}}
\newcommand{\lan}{\lambda^{\an}}

\def\frak{\mathfrak}
\def\text{\textrm}

\def\heightone{\wp}
\def\free#1{\tilde{#1}}
\def\normalclosure#1{\hat{#1}}
\def\full#1{{#1}^{\circ}}
\def\et{\text{\'et}}
\def\iso{\cong}
\def\cotimes{\hat{\otimes}}
\def\formalilim#1
{\text{``$\underset{\substack{\longrightarrow\\ #1}}{\lim}$''}}
\def\ilim#1{\underset{\substack{\longrightarrow\\#1}}{\lim}}
\def\formalplim#1
{\text{``$\underset{\substack{\longleftarrow\\\ #1}}{\lim}$''}}
\def\plim#1{\underset{\substack{\longleftarrow\\ #1}}{\lim}}
\def\ilim#1{\underset{\substack{\longrightarrow\\ #1}}{\lim}}
\def\content#1{I(#1)}
\def\npL#1#2{L_{#1}(#2,\omega^i)}
\def\pL#1#2{L(#1,#2,\omega^i)}
\def\modpL#1{\overline{L}_{#1}(\rhobar,\omega^i)}
\def\hom#1#2{H_1(#1;#2)}
\def\homc#1#2{H_1(#1,\{\text{cusps}\};#2)}
\def\cohom#1#2{H^1(#1;#2)}

\def\pages#1{pp.~#1}

\title{Variation of Iwasawa invariants in Hida families}
\author{Matthew Emerton, Robert Pollack and Tom Weston}
\address[Matthew Emerton]{Mathematics Department, Northwestern
University, Evanston, IL}
\address[Robert Pollack]{Department of Mathematics, University of Chicago,
Chicago, IL}
\address[Tom Weston]{Department of Mathematics and Computer Science,
Amherst College, Amherst, MA}
\email[Matthew Emerton]{emerton@math.northwestern.edu}
\email[Robert Pollack]{pollack@math.uchicago.edu}
\email[Tom Weston]{taweston@amherst.edu}

\begin{document}

\maketitle 

\section{Introduction}

Let $\rhob : \GQ \to \GL_{2}(k)$ be an absolutely irreducible
modular Galois representation over a finite field $k$ of characteristic $p$.
Assume further that $\rhob$ is $p$-ordinary and $p$-distinguished in the sense
that the restriction of $\rhob$ to a decomposition group at $p$ is
reducible and non-scalar.
The {\it Hida family}
$\Hida(\rhob)$ of $\rhob$ 
is the set of all $p$-ordinary $p$-stabilized newforms $f$
with mod $p$ Galois representation isomorphic to $\rhob$.  (If $\rhob$ is
unramified at $p$, then one must also fix an unramified line in $\rhob$
and require that the ordinary line of $f$ reduces to this fixed line.)  
These newforms are a dense set of points in a certain $p$-adic
analytic space of overconvergent eigenforms, consisting of an
intersecting system of branches (i.e.\ irreducible components)
$\T(\a)$ indexed by the minimal primes $\a$ of a certain Hecke algebra.

To each modular form $f \in \Hida(\rhob)$ one may associate the
{\it Iwasawa invariants} $\mu^{\an}(f)$, 
$\lambda^{\an}(f)$, $\mu^{\alg}(f)$, and 
$\lambda^{\alg}(f)$.
The analytic (resp.\ algebraic) $\lambda$-invariants
are the number of zeroes of the
$p$-adic $L$-function (resp.\ of the characteristic power series of the
dual of the Selmer group) of $f$, while the 
$\mu$-invariants are the exponents of the powers
of $p$ dividing the same objects.
In this paper we prove the
following results on the behavior of these Iwasawa invariants as
$f$ varies over $\Hida(\rhob)$.  

\begin{mthm} \label{mthm:1}
Fix $* \in \{\alg,\an\}$.
If $\mu^{*}(f_{0}) = 0$ for some $f_{0} \in \Hida(\rhob)$,
then $\mu^{*}(f) = 0$ for all $f \in \Hida(\rhob)$.
(We then write simply $\mu^{*}(\rhob)=0$.)
\end{mthm}

It is conjectured by Greenberg that if a $p$-ordinary modular form $f$ of weight
two has a residually irreducible Galois representation, then
$\mu^{\an}(f)=\mu^{\alg}(f)=0$;
Theorem~\ref{mthm:1} thus shows that this conjecture is equivalent to
the corresponding conjecture for modular forms of arbitrary weight.

\begin{mthm} \label{mthm:2}
Fix $* \in \{\alg,\an\}$ and assume that $\mu^{*}(\rhob)=0$.  Let
$f_{1},f_{2} \in \Hida(\rhob)$ lie on the branches
$\T(\a_{1}),\T(\a_{2})$ respectively.  Then
\begin{equation} \label{eq:main}
\lambda^{*}(f_{1}) - \lambda^{*}(f_{2}) =
\sum_{\ell \mid N_{1}N_{2}} e_{\ell}(\a_{2}) - e_{\ell}(\a_{1});
\end{equation}
here the sum is over all primes dividing the tame level of $f_{1}$ or
$f_{2}$ and $e_{\ell}(\a_{j})$ is a certain explicit non-negative
invariant of the branch $\T(\a_{j})$ and the prime $\ell$.
\end{mthm}

The first (resp.\ second) theorem corresponds to Theorems
\ref{thm:prealgebraic} and
\ref{thm:mu vanishing} (resp.\ 
Theorems \ref{thm:algebraic} and \ref{thm:lambda}) with trivial
twist by the mod $p$ cyclotomic character.

Note that the right-hand side of (\ref{eq:main}) is identical both
algebraically and analytically.  In particular, using
work of Kato on the main conjecture of Iwasawa theory for
modular forms we obtain the following result.

\begin{mcor} \label{mcor:1}
Assume that $\mu^{\alg}(\rhob)=\mu^{\an}(\rhob)=0$.  
If the main conjecture holds for some
$f_{0} \in \Hida(\rhob)$, then it holds  for
all $f \in \Hida(\rhob)$.
\end{mcor}

The corollary in particular reduces the main conjecture
to the case of modular forms of weight two,
together with the conjecture on the vanishing of the $\mu$-invariant.
We also obtain the following result on the variation of
$\lambda$-invariants in a Hida family. 

\begin{mcor} \label{mcor:2}
Fix $* \in \{\alg,\an\}$ and assume that $\mu^{*}(\rhob)=0$.
\begin{enumerate}
\item $\lambda^{*}(\cdot)$ is constant on branches
of $\Hida(\rhob)$.
\item $\lambda^{*}(\cdot)$ is minimized on the branches of
$\Hida(\rhob)$ of minimal tame level.
\end{enumerate}
\end{mcor}

This work was motivated by the paper \cite{GV} of Greenberg and Vatsal
where Theorems~\ref{mthm:1} and~\ref{mthm:2} were obtained for $p$-ordinary
modular forms of weight two.
Our results here
help to illuminate the results of \cite{GV}; indeed,
if two congruent elliptic curves have different $\lambda$-invariants,
it follows from these results that they must lie on two branches
of the associated Hida family with different ramification behavior. 
One may thus think of the change of $\lambda$-invariants
in terms of ``jumps'' as one moves from one branch to another at
crossing points (which necessarily occur at non-classical\footnote{To
avoid circumlocutions, we adopt the convention that weight one is
a non-classical weight.} eigenforms).

\subsection{Iwasawa theory}
\label{subsec:iwasawa intro}

The papers \cite{greenberg} of Greenberg
and \cite{mazur2} of Mazur introduce the following point of view on
Iwasawa theory:  If $X$ is a $p$-adic analytic family
of $p$-adic representations of $G_{\Q}$ interpolating a collection
of motivic representations, then there should
exist an analytic function (or perhaps a section of an invertible sheaf)
$L$, which we call a $p$-adic $L$-function, 
defined on $X$, and
a coherent sheaf $H$ (the sheaf of Selmer groups) on $X$,
such that the zero locus of $L$ 
coincides with the codimension one part of the support of $H$ (thought of as a Cartier divisor on $X$).
The function $L$ should $p$-adically interpolate the (suitably modified)
values at $s = 1$ of the classical $L$-functions attached to 
the Galois representations at motivic points of $X$,
while at such a motivic point $x \in X$, the Selmer group
$H_x$ (the fibre of $H$ at $x$)
should coincide with the Bloch--Kato Selmer group (computed with respect to
the deRham condition at $p$).

The $L$-function considered above is what is usually known as
the {\em analytic} $p$-adic $L$-function in classical Iwasawa theory;
the equation of the codimension one part of the support of $H$ is
usually known as the {\em algebraic} $p$-adic $L$-function.  
The equality of the zero locus of $L$ with the divisorial
part of the support of $H$ is thus a statement of the 
{\em main conjecture of Iwasawa theory} for the family of 
Galois representations $X$.

In this paper, we restrict to the case of two dimensional nearly
ordinary modular representations, in which case the space
$X$ can be described via Hida theory.
Consider a nearly ordinary residual representation;
such a representation may be written in the form
$\rhobar\otimes \omega^i: G_{\Q} \rightarrow \GL_2(k),$
where $\rhobar$ is as above, $\omega$ denotes the mod $p$ cyclotomic
character, and $0 \leq i \leq p-1$.
As $\rhobar$ is ordinary, we may (and do) fix a 
{\em $p$-stabilisation}
of $\rhobar\otimes \omega^i$: a choice of one dimensional $G_{\Q_p}$-invariant
quotient of $\rhobar\otimes \omega^i$.

Consider the universal deformation
space $Y^{\nord}_{\Sigma}$
parameterising all nearly ordinary deformations of the
$p$-stabilized representation $\rhobar\otimes \omega^i$ which are
unramified at primes in $\Sigma \cup \{p\}$ for some fixed finite set
of primes $\Sigma$.   Any such deformation is of the form
$\rho\otimes\omega^i\otimes\chi$, where $\rho$ is an ordinary
deformation of $\rho$, and $\chi$ is a character of $\Gamma:=1 + p \Z_p$.
(We regard $\chi$ as a Galois character by composing it with the
projection onto $\Gamma$ of the  $p$-adic cyclotomic character.)
If we let $Y^{\ord}_{\Sigma}$ denote the universal deformation space
parameterising all ordinary deformations of the
$p$-stabilized representation $\rhobar$ which are
unramified at primes in $\Sigma \cup \{p\}$, then
$Y^{\nord}_{\Sigma}$ is 
equal to the product (as formal schemes) of $Y^{\ord}_{\Sigma}$
and $\Spf \left(\Z_p[[\Gamma]]\right).$
By results of Wiles and Taylor-Wiles \cite{wiles, taylor-wiles},
as strengthened by Diamond \cite{diamond2}, 
the deformation spaces $Y^{\nord}_{\Sigma}$ can be identified with 
various local
pieces of the universal ordinary Hecke algebra constructed by Hida
(at least under mild hypotheses on $\rhobar$).

If we allow $\Sigma$ to vary, the spaces $Y^{\nord}_{\Sigma}$
form a formal Ind-scheme over $W(k)$ parameterizing
nearly ordinary deformations of $\rhobar$ of arbitrary conductor.
Let $Y^{\nord}$
be an irreducible component of this formal Ind-scheme;
we may write $Y^{\nord}$ as the product of $Y^{\ord}$ (an irreducible
component of the the formal Ind-scheme parameterizing ordinary
deformations of $\rhobar$) and $\Spf \left(\Z_p[[\Gamma]]\right)$.
The motivic points
are dense on $Y^{\nord}$
and all have the same tame conductor,
which we will denote by $N$.  
If $\T_N^{\new}$ denotes the new-quotient of the Hida Hecke 
algebra of level $N$,
then $Y^{\ord}$
corresponds to a minimal prime $\a$ of $\T_N^{\new}$, and we set
$X^{\ord} = \Spf\left(\T_N^{\new}/\a\right)$ and
$X^{\nord} = \Spf \left(\T_N^{\new}/\a \, [[\Gamma]]\right)$. 
(We point out that $X^{\ord}$
is nearly the same as $Y^{\ord}$
-- and hence $X^{\nord}$ is nearly the same as $Y^{\nord}$.
Indeed, $X^{\ord}$ is a finite cover of $Y^{\ord}$ 
and maps isomorphically to it after inverting $p$;
see section \ref{sec:branches} for details.)

The generalities of Iwasawa theory
discussed above apply to this space $X^{\nord}$.
In section \ref{sec:analytic L-functions} we give
a construction of the analytic $p$-adic $L$-function
on $X^{\nord}$; as is usual, we regard it as a ``function of
two-variables'' -- one variable moving along $X^{\ord}$,
and the other along $\Spf \left(\Z_p[[\Gamma]]\right).$  In our discussion
of Selmer groups we do not go so far as to construct a
coherent sheaf of Selmer groups over the space $X^{\nord}$;
for this, the reader should refer to recent and forthcoming
work of Ochiai \cite{Ochiai}.

An important feature of the situation we consider is that
the nearly ordinary deformation space contains many different irreducible
components.  The main goal of this paper is to describe how to pass
Iwasawa theoretic information (such as knowledge of the
main conjecture) between components.

\subsection{Two-variable $p$-adic $L$-functions}
\label{subsec:two-variable intro}

To establish the analytic parts of Theorems \ref{mthm:1} and \ref{mthm:2}, we make use of
a two-variable $p$-adic $L$-function (playing the role of $L$ above).  
As already remarked upon in~\ref{subsec:iwasawa intro},
one variable is a parameter
of wild character space
(i.e.\ the standard cyclotomic variable) and the second variable
runs through the Hida family (which is a finite cover of weight space)
of some fixed residual representation $\rhobar$.
Of the many constructions of two-variable $p$-adic $L$-functions (see \cite{GS}, \cite{Stevens},
\cite{Pan}), our construction 
follows most closely that of Mazur \cite{mazur3} and 
Kitigawa \cite{Kitigawa}, with two main differences.

The first difference is that by using the fact that the Hecke algebras are
known to be Gorenstein (as in \cite{wiles})
we may construct these $L$-function with fewer assumptions on $\rhobar$.  
The second
difference is that we do not limit ourselves to working solely 
on a part of the Hida family parameterising newforms of some fixed tame level. 

For each branch $Y$ of the Hida family we construct a two-variable
$L$-function along $Y$, which at each classical point specialises to the
$p$-adic $L$-function of the corresponding newform computed with respect
to its canonical period.  (Actually, the two-variable $L$-function
is defined not on $Y$, but
on the partial normalisation $X$ of $Y$ discussed above
in~\ref{subsec:iwasawa intro}.)
However,
for any finite set of primes $\Sigma$ not containing $p$,
we also construct a two-variable $L$-function 
along $Y_\Sigma$ (the Hida family
of $\rhobar$ with ramification only at the primes of $\Sigma$).
If we specialize this two-variable $p$-adic $L$-function 
at a classical point,
we obtain a non-primitive $p$-adic $L$-function
attached to the corresponding newform
(i.e.\ the usual $p$-adic $L$-function stripped of
its Euler factors at primes of $\Sigma$).
In fact, we prove that this non-primitivity occurs in families:
the $p$-adic $L$-function on the Hida family $Y_{\Sigma}$,
restricted to some branch $Y$,
is equal to the $p$-adic $L$-function of that branch
stripped of its two variable 
Euler factors at the primes of $\Sigma$.  

With this construction in hand,
the analytic part of Theorem \ref{mthm:1} follows 
immediately.  Indeed, the vanishing of the analytic $\mu$-invariant for any particular
form $f$ in the Hida family is equivalent to the vanishing of the $\mu$-invariant
of the two-variable $p$-adic $L$-function along the branch $Y$ containing $f$,
which in turn is equivalent to the vanishing of the
$\mu$-invariant of the non-primitive
two-variable $L$-function on $Y_{\Sigma}$ for any sufficiently large
choice of $\Sigma$ (i.e.\ such that $Y_{\Sigma}$ contains $Y$). 
This last condition is independent of $f$, and so is equivalent to the
vanishing of the $\mu$-invariant for every modular form in the family.  

The analytic part of Theorem \ref{mthm:2} also follows from this construction.
The analytic $\lambda$-invariant of a modular form $f$ is equal to the
$\lambda$-invariant of the two-variable $L$-function attached to the branch containing $f$.
To compare the $\lambda$-invariants of two branches, we choose $\Sigma$ large enough so 
that $Y_\Sigma$ contains both of these branches.  We then relate each $\lambda$-invariant
to the $\lambda$-invariant of the two-variable $p$-adic $L$-function attached to $Y_\Sigma$.
As explained above, the difference of the two $\lambda$-invariants is then realized
in terms of the $\lambda$-invariants of certain Euler factors.
The quantity $e_{\ell}(\a_{i})$ appearing in the formula of
Theorem~\ref{mthm:2} is precisely the $\lambda$-invariant of the
Euler factor at $\ell$ along the branch corresponding to $\a_{i}$.

We point out
that we found it necessary to use two-variable $L$-functions and Hida theory
in order to compare the $\lambda$- and $\mu$-invariants of two congruent modular forms
of different weights.  It would be interesting to know if these results could be obtained
with different methods that do not depend upon using families of modular forms.

\subsection{Residual Selmer groups}

Let $\rhob : \GQ \to \GL_{2}(k)$ be as above.  For any $f \in \Hida(\rhob)$,
say with Fourier coefficients in the finite extension $K$ of $\Qp$,
Greenberg has defined the Selmer group $\Sel(\Qi,\rho_{f})$ 
as the subgroup of $H^{1}(\Qi,A_{f})$ cut out by local conditions, all of
which are as strong as possible except for that at $p$; here
$\Qi$ is the cyclotomic $\Zp$-extension of $\Q$ and $A_{f}$ is
$(K/\O_{K})^{2}$ with Galois action via $\rho_{f}$.  If one defines
$\Sel(\Qi,\rhob)$ in the analogous way, one is
confronted with the fundamental problem that the $\pi$-torsion on
$\Sel(\Qi,\rho_{f})$ may be larger than 
$\Sel(\Qi,\rhob)$; this is, of course, precisely the reason why
congruent modular forms need not have isomorphic Selmer groups.

In \cite{GV} this issue is overcome by introducing non-primitive
Selmer groups of $\rhob$; essentially, if $f_{1}, f_{2} \in \Hida(\rhob)$ 
have tame levels $N_{1}$ and $N_{2}$ respectively, then the $\pi$-torsion on
$\Sel(\Qi,\rho_{f_{1}})$ and $\Sel(\Qi,\rho_{f_{2}})$ can both be
compared to the non-primitive Selmer group of $\rhob$ obtained by
ignoring the local conditions
at primes dividing $N_{1}N_{2}$.

Although this approach can also be made to work
in the higher weight case, here we proceed
somewhat differently.  Following Mazur \cite{MazurES},
we allow non-strict, but not
necessarily vacuous, local conditions $\SS$ 
on the cohomology of $\rhob$, resulting
in a family of residual Selmer groups $\Sel_{\SS}(\Qi,\rhob)$.
We show that for any
$f \in \Hida(\rhob)$ there is a local condition $\SS(f)$ such that
$$\Sel_{\SS(f)}(\Qi,\rhob) \cong \Sel(\Qi,\rho_{f})[\pi].$$
It follows that the Iwasawa invariants of $f$ can be recovered from the
residual Selmer group $\Sel_{\SS(f)}(\Qi,\rhob)$.

In fact, if $f$ has tame level equal to the conductor of $\rhob$,
we show that $\Sel_{\SS(f)}(\Qi,\rhob)$ agrees with
the naive residual Selmer group $\Sel(\Qi,\rhob)$.
Since such $f$ always exist by level lowering, we are then able to use
duality results to show that the difference
$$\dim_{k} \Sel_{\SS(f)}(\Qi,\rhob) - \dim_{k} \Sel(\Qi,\rhob)$$
is precisely given by the dimensions of the local conditions $\SS(f)$.
The main algebraic results follow from this.

We note that Hida theory plays little overt role in these results.
In particular, these methods may well apply more generally for
algebraic groups other than $\GL_{2}$.  The key inputs are
automorphic descriptions of the Galois representation at ramified primes,
level lowering results, and the fact that the Selmer groups of interest
are $\Lambda$-cotorsion; the remainder of the argument is essentially
formal.

\subsection{$L$-functions modulo $p$}

In the paper \cite{mazur1}, Mazur raises the question of
whether one can define a mod $p$ $L$-function attached
to the residual representation $\rhobar$
and a choice of tame
conductor $N$ (divisible by the tame conductor of $\rhobar$). 
The construction discussed in~\ref{subsec:two-variable intro}
gives a positive answer to
this question, with the caveat that the appropriate extra data is not
simply a choice of tame conductor $N$,  but rather the
more precise data of a component $Y$ of the universal deformation space
of $\rhobar$.  We may then specialise the $p$-adic $L$-function
at the closed point of the partial normalisation $X$ of $Y$
on which it is defined, so as to obtain a mod $p$ $L$-function
attached to $\rhobar$ and $Y$.

We also show that the local condition $\SS(f)$ discussed in
the preceding section depends only on the component $Y$ of the Hida
family that contains $f$; we may thus write $\SS(Y)$ for $\SS(f)$.  
Assuming that both the analytic 
and algebraic $\mu$-invariant of $\rhobar$ vanish, Theorem~\ref{mthm:2}
then shows that the main conjecture for any member of $\Hida(\rhobar)$
is equivalent to the following mod $p$ main conjecture (for one,
or equivalently every, choice of $Y$):

\begin{mconj} Let $Y$ be a branch of $\Hida(\rhobar)$.  Then
the mod $p$ $L$-function $L_{p}(Y)$ of $Y$ is non-zero, the Selmer group
$\Sel_{\SS(Y)}(\Q_{\infty})$ is finite, and
$$\lambda\bigl(L_{p}(Y)\bigr) = \dim_k \Sel_{\SS(Y)}(\Q_{\infty},\rhobar).$$
\end{mconj}

\subsection{Overview of the paper}

In the following section, we recall the Hida theory used
in this paper.  In particular, we discuss the Hida family attached
to a residual representation $\rhobar$, its decomposition into
irreducible components and the various Galois
representations attached to this family, with an emphasis on the
integral behavior.

In the third section, we prove the algebraic parts of
Theorems \ref{mthm:1} and \ref{mthm:2} via
the theory of residual Selmer groups and the use of level lowering
to reduce to the minimal case.

In the fourth section, we construct two-variable $p$-adic $L$-functions
on the Hida family of $\rhobar$ and on each of its irreducible components.
By relating these two constructions to each other and to
the construction of classical $p$-adic $L$-functions, we obtain 
proofs of the analytic parts of
Theorems \ref{mthm:1} and \ref{mthm:2}.

In the final section we give applications 
to the main conjecture.   We discuss some
explicit examples to illustrate the general theory, including 
the congruence modulo $11$ between 
the elliptic curve $X_0(11)$ and the modular form $\Delta$.

\subsection{Acknowledgments}\label{subsec:acknowledgments}
The authors wish to thank Ralph Greenberg, Masato Kurihara, Barry Mazur
and Chris Skinner for helpful conversations.


\subsection*{Notation}

We write $\GQ$ for the absolute Galois group of $\Q$; if $\Sigma$ is a finite
set of places of $\Q$, we write $\QS$ for the maximal extension of $\Q$
unramified away from $\Sigma$.  We fix an odd prime $p$ and let
$\Qi$ denote the cyclotomic $\Zp$-extension of $\Q$.
For each prime
$\ell$ of $\Q$ (resp.\ place $v$ of $\Qi$) fix a decomposition group
$\Gl \inj \GQ$ (resp.\ $\Gv \inj \GQi$) with inertia group
$\Il$ (resp.\ $\Iv$) such that $\Gv \subseteq \Gl$ whenever
$v$ divides $\ell$; note that $\Il = \Iv$ for a place $v$ dividing
a prime $\ell \neq p$.
Let $v_{p}$ denote the unique place of $\Qi$ above
$p$.

We write $\vep : \GQ \to \Zp^{\times}$ for the $p$-adic cyclotomic character.
We let $\Gamma$ denote the group of 1-units in $\Z_p^{\times}$; the
cyclotomic character induces a canonical isomorphism
$\Gal(\Qi/\Q) \overset{\cong}{\too} \Gamma$.
If $\O$ is a $\Zp$-algebra, we let $\Lambda_\O$ denote the
completed group ring $\O[[\Gamma]]$; we simply write $\Iw$ for
$\Lambda_{\Zp}$.
We denote the natural map $\Gamma \rightarrow \Iw^{\times}$
by $\gamma \mapsto \langle \gamma \rangle$. 
Recall that $\Iw$ is a complete regular local
ring of dimension two, non-canonically isomorphic to the power series
ring $\Z_p[[T]]$.  (Such an isomorphism is obtained by choosing a topological
generator $\gamma$ of $\Gamma$ and mapping $\langle \gamma \rangle $
to $T +1$.)
We let ${\mathcal L}$ denote the field of fractions of $\Iw$.

We let $\Delta$ denote the group of cyclotomic units of $\Z_p$;
there is then an isomorphism $\Gamma \times \Delta \iso \Z_p^{\times}.$ 
We let $\omega$ denote the inclusion of $\Delta$ in $\Z_p^{\times}$,
so that 
$$\Hom(\Delta,\Z_p^{\times}) = \{ \omega^i \, | \, 0 \leq i \leq p-2\}.$$
If $\Z_{p,(i)}$ denotes $\Z_p$ regarded as a $\Delta$-module
via the character $\omega^i$, then there is a natural isomorphism
of $\Z_p[\Delta]$-algebras
$\Z_p[\Delta] \iso \prod_{i = 0}^{p-2} \Z_{p,(i)}.$
Combining this with the natural isomorphism
$\BbbIw \iso \Z_p[\Delta] \otimes_{\Z_p} \Iw,$
we obtain an isomorphism of $\BbbIw$-algebras
$$\BbbIw \iso \prod_{i=0}^{p-2} \Iw_{(i)},$$
where $\Iw_{(i)}$ denotes a copy of $\Iw$, enhanced
to a $\BbbIw$-algebra by having $\Delta$ act
through~$\omega^i$.

If $\heightone$ is a height one prime of $\Iw$,
then we write $\O(\heightone) := \Iw/\heightone$.
Note that $\O(\heightone)$ is a complete local domain of dimension one,
either isomorphic to $\F_p[[\Gamma]]$ (if $\heightone = p\Iw$)
or else free of finite rank over
$\Z_p$ (in which case we will say that $\heightone$ is of residue
characteristic zero).
The natural embedding $\Gamma \rightarrow \Iw^{\times}$
induces a character
$\kappa_{\heightone}:\Gamma \rightarrow \O(\heightone)^{\times}$.
We say that $\heightone$ is {\em classical} if it is of residue
characteristic zero, and if there is a finite
index subgroup $\Gamma'$ of $\Gamma$ such that $\kappa_{\heightone},$
when restricted to $\Gamma'$, coincides with the character
$\gamma \mapsto \gamma^k \in \Z_p^{\times} \subset \O(\heightone)^{\times},$
for some integer $k\geq 2$.  If we wish to specify $k$,
then we will say that $\heightone$ is {\em classical of weight~$k$}.

More generally, if $\heightone$ is a height one prime ideal 
of a finite $\Iw$-algebra $\T$
then we will again write $\O(\heightone) := \T/\heightone,$
and if $\heightone' := \heightone \bigcap \Iw$,
we will write
$\kappa_{\heightone} := \kappa_{\heightone'}
: \Gamma \rightarrow \O(\heightone')^{\times} \subset \O(\heightone)^{\times}.$
We say that $\heightone$
is classical (of weight $k$)
if $\heightone'$ 
is classical (of weight $k$).   

If $N$ is a natural number prime to $p$, then we will
write $\prodNp := (\Z/N)^{\times} \times \Z_p^{\times}$.
We will have occasion to
consider the completed group ring $\Z_p[[\prodNp]].$
In this situation,
we will extend the diamond bracket notation introduced above, and write
$x \mapsto \langle x \rangle$ to denote the natural map
$\prodNp \rightarrow
\Z_p[[\prodNp]]^{\times}.$

\section{Hida theory}\label{sec:hida}

\subsection{The universal ordinary Hecke algebra}\label{subsec:hida}

We begin by briefly recalling Hida's construction of the
universal ordinary Hecke algebra of tame level $N$.

Given an integer $k$,
we let $S_k(N p^{\infty},\Z_p)$ denote the
space of weight $k$ cusp forms that are on $\Gamma_1(N p^r)$ for some
$r \geq 0,$ and whose
$q$-expansion coefficients (at the cusp $\infty$) lie in $\Z_p$.
We define an action of $\Z_p^{\times}$ on $S_k(N p^{\infty},\Z_p)$
by taking the product of the nebentypus action
with the character $\gamma \mapsto \gamma^k$;
in this way $S_k(N p^{\infty}, \Z_p)$ becomes a $\BbbIw$-module.
Also, for each prime $\ell \neq p$, the Hecke operator $T_{\ell}$ acts on
$S_k(N p^{\infty}, \Z_p),$ as does the Hecke operator $U_p$.
Finally, the group $(\Z/N)^{\times}$ acts on $S_k(N p^{\infty}, \Z_p)$,
via the nebentypus action.
Thus if we write
$$\BigHecke := \Z_p[[\prodNp]] [\{T_\ell\}_{\ell \neq p},U_p],$$
then
$S_k(N p^{\infty}, \Z_p)$ is naturally an
$\BigHecke$-module.
More generally, for any $\Z_p$-algebra $R$,
we write $S_k(N p^{\infty},R) :=
S_k(N p^{\infty},\Z_p) \otimes_{\Z_p} R.$ 
The $\BigHecke$-action on $S_k(N p^{\infty},\Z_p)$
extends uniquely to an $R$-linear action on $S_k(N p^{\infty},R)$.

Taking $q$-expansions yields an injection of $\Q_p$-algebras
\begin{equation}\label{eqn:q-exp}
\bigoplus_{k\geq 0} S_k(N p^{\infty},\Q_p) \rightarrow \Q_p[[q]];
\end{equation}
we let $D$ denote the $\Z_p$-subalgebra of the source of~(\ref{eqn:q-exp})
obtained as the preimage of the subalgebra $\Z_p[[q]]$ of its target.
(This is a variant of Katz's ring of {\em divided congruences} \cite{katz}.)
If $\hat{D}$ denotes the $p$-adic completion of $D$, then $\hat{D}$
is naturally isomorphic to Katz's ring of {\em generalised $p$-adic
modular functions} \cite{katz}.
One shows that the action of
$\BigHecke$
on the source of~(\ref{eqn:q-exp})
(induced by the action on each of the direct summands) 
restricts to an action of
$\BigHecke$
on $D$.
This in turn induces an action of
$\BigHecke$
on $\hat{D}$.

Hida has defined the ordinary projector $e^{\ord}$ acting on
$\hat{D}$; it cuts out the submodule $\hat{D}^{\ord}$
of $\hat{D}$ on which $U_p$ 
acts as an isomorphism.  The projector $e^{\ord}$ commutes with
the operators in $\BigHecke,$ and so
$\hat{D}^{\ord}$ is closed under the action of this algebra.
The action of
$\BigHecke$
on $\hat{D}$ induces a map
\begin{equation}\label{eqn:hecke map}
\BigHecke \rightarrow \End_{\Z_p}(\hat{D}^{\ord}).
\end{equation}

\begin{definition} We let $\T_N$ denote the image of the map~(\ref{eqn:hecke map}),
and refer to $\T_N$ as the universal ordinary Hecke algebra of tame level~$N$.
\end{definition}

Before stating Hida's results concerning $\T_N$,
we introduce some additional notation.
For any $\Z_p$-algebra $R$, and any integer $k$,
we let $S_k(N p^{\infty}, R)^{\ord}$ denote the $R$-submodule
of $S_k(N p^{\infty}, R)$ obtained as the image of the idempotent $e^{\ord}$.
This is an $R$-submodule of $S_k(N p^{\infty}, R)$,
closed under the action of
$\BigHecke$, and by construction, the action of $\BigHecke$ on
$S_k(N p^{\infty}, R)^{\ord}$ factors through its quotient $\T_N$.

If $\kappa$ is any element of $\Hom(\Gamma,R^{\times})$, we let
$S_k(N p^{\infty},R)^{\ord}[\kappa]$ denote the $R$-submodule
of $S_k(N p^{\infty}, R)^{\ord}$ on which $\Gamma$ acts via the character
$\kappa$.  Note that this is in fact a finitely generated $R$-module,
because the $p$-part of the conductor of a $p$-ordinary 
eigenform is bounded by the $p$-part of the conductor of its nebentypus
(unless the eigenform has trivial $p$-nebentypus, in which case the $p$-part of
its conductor is bounded by $p$).

\begin{thm}\label{thm:hida}
~
\begin{enumerate}
\label{part:hida}
\item \label{part:hida1}
 The algebra $\T_N$ is free of finite rank over $\Iw$.
(It is regarded as a $\Iw$-algebra via the inclusion $\Iw\subset
\Z_p[[\prodNp]]$.)

\item \label{part:hida2}
If $\heightone$ is a classical height one prime in $\Iw$ of weight $k$,
then the surjection
$$\BigHecke \rightarrow \T_N \rightarrow \T_N/\heightone \T_N$$
identifies 
$\T_N/\heightone \T_N$ with
the quotient of the Hecke algebra $\BigHecke$ that acts faithfully on
the space 
$S_k(N p^{\infty}, \O(\heightone))^{\ord}[\kappa_{\heightone}].$
\end{enumerate}
\end{thm}
\begin{proof}
See~\cite[Thm.~3.1, Cor.~3.2]{hida1},
and~\cite[Thm.~1.1, Thm.~1.2]{hida2}.
\end{proof}

Part~\ref{part:hida1} of this theorem implies that $\T_N$ is a complete Cohen-Macaulay
semi-local ring of dimension two. 

It follows from part~\ref{part:hida2} of the theorem,
together with the usual duality between spaces of modular forms
and Hecke algebras, that a classical height one prime ideal
$\heightone$ in $\T_N$ gives rise to a normalised Hecke eigenform
in $S_k(N p^{\infty}, \O(\heightone))^{\ord}[\kappa_{\heightone}].$
We denote this eigenform by $f_{\heightone}$.

Suppose conversely that $f$ is a normalised Hecke eigenform
in $S_k(N p^{\infty}, K)^{\ord}[\kappa]$, for some finite extension
$K$ of $\Q_p$, and some character $\kappa: \Gamma \rightarrow K^{\times}$.
If $\O(f)$ denotes the finite extension of $\Z_p$ generated by the
$q$-expansion coefficients of $f$, then $\kappa$ in fact takes values
in $\O(f)^{\times}$, and $f$ lies in $S_k(N p^{\infty}, \O(f))^{\ord}[\kappa]$.
The Hecke action on $f$ induces a homomorphism
$\BigHecke \rightarrow \O(f)$,
which factors as 
$$\BigHecke \rightarrow \T_N \rightarrow \T_N/\heightone_f
\buildrel \iota \over \iso \O(f)$$
for some weight $k$ classical height one prime ideal $\heightone_f$ of 
$\T_N$, and some isomorphism $\iota:~\O(\heightone_f) \iso \O(f)$.
By construction, the isomorphism $\iota$ identifies the normalised
eigenform $f_{\heightone_f}$ corresponding to $\heightone_f$ with $f$,
and the character $\kappa_{\heightone_f}$ with the character $\kappa$.

To summarise, the classical height one primes in $\T_N$ correspond to
Galois conjugacy classes of $p$-ordinary normalised Hecke eigenforms of tame level $N$
and weight $k \geq 2$ defined over an algebraic closure of $\Qbar_p$.

We close this section by remarking that
it follows from Theorem~\ref{thm:hida}
that the algebra $\T_N$ acts faithfully on the space
$S_k(N p^{\infty})^{\ord}$, for any (fixed) weight~$k \geq 2$.
(Part~(2) of the theorem
shows that any element that annihilates this space lies in
the intersection of the ideals $\heightone \T_N$,
where $\heightone$ ranges over all classical height one prime
ideals of weight $k$.  Since $\T_N$ is free over $\Iw$, and
since the intersection of all such ideals $\heightone$ is the
zero ideal of $\Iw$, we see that any such element must vanish.)
Thus we could have defined $\T_N$ by considering 
any fixed weight $k\geq 2$, rather than taking the direct sum over all
weights.  (This is the approach to defining $\T_N$ taken in \cite{hida2}.)
Similarly, we could have defined $\T_N$ by forming the
direct sum over all weights $k$ of the spaces of cusp forms
with trivial $p$-power nebentypus and $\Q_p$ coefficients,
and then forming the analogue of $D$ and $\hat{D}$ above.
(This is the approach to defining $\T_N$ taken in \cite{hida1}.)  
We chose to incorporate all weights and $p$-power levels in
our definition of $\T_N$ so as to have specialisations
to all weights and nebentypus characters be made apparent from the outset.

\subsection{The new quotient}\label{subsec:new}


Given an integer $k$,
we denote by $S_k(N p^{\infty}, \Z_p)^{\ord}_{\old}$ 
the subspace of $S_k(N p^{\infty}, \Z_p)^{\ord}$ that consists
of forms that are spanned (over $\Q_p$) by forms that are
old at level $N$ (i.e.\ that arise from $p$-ordinary
forms of level $p^{\infty} M$,
for some proper divisor $M$ of $N$, by applying one of the various
degeneracy maps $q\mapsto q^d,$ corresponding to the divisors $d$ of $N/M$).
Note that $S_k(N p^{\infty}, \Z_p)^{\ord}_{\old}$ is closed under the action
of $\T_N$ on $S_k(N p^{\infty}, \Z_p)^{\ord}$.
We let $S_k(N p^{\infty}, \Z_p)^{\ord}_{\new}$ denote the
$\T_N$-module obtained as the quotient of $S_k(N p^{\infty}, \Z_p)^{\ord}$ by
its submodule $S_k(N p^{\infty}, \Z_p)^{\ord}_{\old}$.
Note that by construction
$S_k(N p^{\infty}, \Z_p)^{\ord}_{\new}$ is $p$-torsion free.

If $R$ is any $\Z_p$-algebra, then we may apply an extension of scalars
from $\Z_p$ to $R$ to each of the modules of oldforms and newforms
constructed in the preceding paragraph.
We denote the resulting modules
in the obvious way; they sit in the short exact sequence
$$0 \rightarrow S_k(N p^{\infty}, R)^{\ord}_{\old} \rightarrow 
S_k(N p^{\infty}, R)^{\ord} \rightarrow S_k(N p^{\infty},R)^{\ord}_{\new}
\rightarrow 0.$$
If $\kappa \in \Hom(\Gamma,R^{\times})$, let
$S_k(N p^{\infty},R)^{\ord}_{\new}[\kappa]$ denote the submodule of 
$S_k(N p^{\infty},R)^{\ord}_{\new}$ on which $\Gamma \subset \Iw^{\times}
\subset \T_N^{\times}$ acts through $\kappa$.

\begin{definition} We let $\T_N^{\new}$ denote the maximal quotient of $\T_N$
that acts faithfully on any one (or equivalently all)  of the modules
$S_k(N p^{\infty}, \Z_p)^{\ord}_{\new}.$  (The claimed equivalence follows
from the discussion at the end of the preceding section.)
\end{definition}

The following result is the analogue for $\T_N^{\new}$ of
Theorem~\ref{thm:hida}.
We include a proof, since it is phrased slightly differently to the
corresponding results in \cite{hida1} and \cite{hida2}.

\begin{thm}\label{thm:new} 
~
\begin{enumerate}
\item $\T_N^{\new}$ is a finite and reduced torsion free $\Iw$-algebra.
 
\item
The classical height one primes of $\T_N^{\new}$ correspond 
(under pull-back) to those classical height one primes $\heightone$ of $\T_N$
for which the corresponding normalised eigenform
$f_\heightone$
(as described in the discussion at the end of section~\ref{subsec:hida})
is of tame conductor $N$.

\item If $\heightone$ is a classical height one prime ideal of $\Iw$,
then $\T_N^{\new} \otimes_{\Iw} \Iw_{\heightone}$ is a finite \'etale
extension of the discrete valuation ring $\Iw_{\heightone}$.
\end{enumerate}
\end{thm}
\begin{proof}
Let us define $\T_N^{\old}$ in analogy to $\T_N^{\new}$,
and let $I^{\old}$ denote the kernel of the surjection
$\T_N \rightarrow \T_N^{\old}.$  
Recalling that $\mathcal L$ denotes the field of fractions of $\Iw$,
consider the short exact sequence
$$0 \rightarrow I^{\old} \otimes_{\Iw} {\mathcal L} \rightarrow
\T_N\otimes_{\Iw} {\mathcal L} \rightarrow
\T_N^{\old} \otimes_{\Iw} {\mathcal L} \rightarrow 0.$$

The tensor product
$\T_N \otimes_{\Iw} \mathcal L$ is a product of Artin local $\mathcal L$-algebras;
more precisely,
we have an isomorphism
$$\T_N \otimes_{\Iw} {\mathcal L} \iso
\prod_{j \in J} (\T_N)_{\frak a_j},$$
where $\{\frak a_j\}_{j \in J}$ is the set of minimal primes
of $\T_N$.
From \cite[Cor.~3.3]{hida1} we see that $I^{\old}$ has trivial
intersection with the nilradical of $\T_N \otimes_{\Iw} {\mathcal L}$,
and thus that $I^{\old}\otimes_{\Iw} {\mathcal L}$
is isomorphic to a product $\prod (\T_N)_{\frak a_j},$
where $\frak a_j$ ranges over a subset of the set of minimal primes,
indexed by some subset $J'$ of $J$.  The components $(\T_N)_{\frak a_j}$
for $j \in J'$ are precisely the {\em primitive} components of
$\T_N\otimes_{\Iw} {\mathcal L}$, in the terminology of \cite{hida1}.
These components are finite field extensions of $\mathcal L$ (since they
have trivial nilradical).

Write $\T_N\otimes_{\Iw} \mathcal L$ as a product of two factors, as follows:
$$\T_N\otimes_{\Iw} {\mathcal L} \iso
\left( \prod_{j \in J'} (\T_N)_{\frak a_j} \right ) \times
\left( \prod_{j \in J \setminus J'} (\T_N)_{\frak a_j} \right  ).$$
We let $\free{\T}$ denote the {\em free closure} (in the terminology
of \cite{hida1}) of the image
of $\T_N$ under the projection onto the first factor.  
If $\heightone$ is a classical height one prime of $\Iw$,
then part~(2) of Theorem~\ref{thm:hida} identifies
$(\T_N/\heightone\T_N)[1/p]$ 
with the $\Q_p$-Hecke algebra acting on the space
$S_k(N p^{\infty}, \O(\heightone) [1/p])^{\ord}[\kappa_\p]$.
We see from 
\cite[Cor.~3.7]{hida1} that under this isomorphism,
the surjective map 
$ (\T_N/\heightone\T_N)[1/p] \rightarrow
(\free{\T}/\heightone\free{\T})[1/p]$
is identified with the surjection from the
$\Q_p$-Hecke algebra acting on
$S_k(N p^{\infty}, \O(\heightone) [1/p])^{\ord}[\kappa_\p]$
onto the $\Q_p$-Hecke algebra acting on
$S_k(N p^{\infty}, \O(\heightone) [1/p])^{\ord}_{\new}[\kappa_\p]$.

Thus we see that the image of $\T_N$ under the map
$\T_N \rightarrow \free{\T}$, that is,
the image of $\T_N$ in $\prod_{j \in J'} (\T_N)_{\frak a_j},$
is precisely equal to $\T_N^{\new}$.  
Since this product is a product of fields, we see that $\T_N^{\new}$ is
$\Iw$-torsion free and reduced.   Being a quotient of the
finite $\Iw$-module $\T_N$, it is also finite over $\Iw$,
and so~(1) is proved. 
We also see that $\T_N^{\new}/ (\T_N^{\new} \bigcap \heightone
\free{\T})$
maps isomorphically onto the Hecke algebra acting on
$S_k(N p^{\infty}, \O(\heightone))^{\ord}_{\new}[\kappa_{\heightone}],$
proving~(2).  (Note that since $\free{\T}$ contains $\T_N^{\new}$
with finite index, the height one primes of $\free{\T}$ correspond
bijectively to those of $\T_N^{\new}$ under restriction.)

Since $\free{\T}$ contains $\T_N^{\new}$ with finite index,
the natural map $\T_N^{\new} \otimes_{\Iw} \Iw_{\heightone}
\rightarrow \free{\T} \otimes_{\Iw} \Iw_{\heightone}$ is an isomorphism.
Hecke algebras acting on spaces of newforms are reduced,
and so we see that the quotient
$
(\T_N^{\new} \otimes_{\Iw} \Iw_{\heightone})/\heightone
(\T_N^{\new} \otimes_{\Iw} \Iw_{\heightone}) \iso
(\free{\T} \otimes_{\Iw} \Iw_{\heightone})/\heightone
(\free{\T} \otimes_{\Iw} \Iw_{\heightone})$
is reduced.  This proves part~(3) (which is simply a rewording
of \cite[Cor.~1.4]{hida2}).
\end{proof}

\subsection{Galois representations}\label{subsec:galois}

In this section we recall the basic facts regarding
Galois representations attached to Hida families.
As above, we fix a tame level $N$ prime to $p$.
Recall that $\mathcal L$ denotes the field of fractions of $\Iw$.

\begin{thm}\label{thm:galois}  
There is a continuous
Galois representation 
$$
\rho: G_{\Q} \rightarrow \GL_2(\T_N^{\new}\otimes_{\Iw} \mathcal L),
$$
characterised by the following properties:
\begin{enumerate}

\item $\rho$ is unramified away from $N p$.

\item If $\ell$ is a prime not dividing $N p$ then $\rho(\Frob_{\ell})$ has
characteristic polynomial equal to
$X^2 - T_{\ell} X + \langle \ell \rangle \ell^{-1}
\in \T_N^{\new}[X].$
(Recall that $\langle \ell \rangle$ denotes the unit
in $\Z_p[[\prodNp]]$ corresponding to the element $\ell \in \prodNp$.)
\end{enumerate}
The representation $\rho$ satisfies the following additional properties:

\begin{enumerate}
\setcounter{enumi}{2}

\item $\rho$ is absolutely irreducible.

\item The determinant of $\rho$ is equal to the following character:
$$G_{\Q} \rightarrow \hat{\Z}^{\times} \rightarrow \prodNp 
\buildrel x \mapsto \langle x \rangle x^{-1} \over \longrightarrow
\Z_p[[\prodNp]]^{\times} \subset (\T_N^{\new})^{\times}.$$
(Here the first arrow is the full cyclotomic character.)

\item The space of $I_p$-coinvariants of $\rho$ is free of rank one
over $\T_N^{\new}\otimes_{\Iw} \mathcal L,$
and $\Frob_p$ acts on this space through the eigenvalue 
$U_p \in \T_N^{\new}$.
%
%
\end{enumerate}
\end{thm}
\begin{proof}
See \cite[Thm.~2.1]{hida2} and \cite[Thm.~2.6]{GS}.
\end{proof}

From the Galois representation $\rho$ we may construct
various related Galois representations.  We recall some of these
constructions here, and introduce the corresponding notation.

First, we may localise $\rho$ in various ways.
If $\frak a$ is a minimal prime ideal of $\T_N^{\new}$,
then we let $\rho_{\frak a}$ denote the representation
$$\rho_{\frak a}: G_{\Q} \rightarrow
\GL_2((\T_N^{\new})_{\frak a})$$ 
obtained by composing $\rho$ with the projection
$\T_N^{\new} \otimes_{\Iw} {\mathcal L} \rightarrow (\T_N^{\new})_{\frak a}.$

If $\frak m$ is a maximal ideal of $\T_N^{\new}$, then
the localisation $(\T_N^{\new})_{\frak m}$ is a direct factor
of $\T_N^{\new}$, and so $(\T_N^{\new})_{\frak m} \otimes_{\Iw} \mathcal L$
is a direct factor of $\T_N^{\new}\otimes_{\Iw} \mathcal L$.
We let $\rho_{\frak m}$ denote the representation
$$\rho_{\frak m}: G_{\Q} \rightarrow
\GL_2((\T_N^{\new})_{\frak m}\otimes_{\Iw} \mathcal L)$$
obtained by composing $\rho$ with the projection
$\T_N^{\new} \otimes_{\Iw} {\mathcal L} \rightarrow (\T_N^{\new})_{\frak m}
\otimes_{\Iw} \mathcal L.$

If $\normalclosure{\T}_N^{\new}$ denotes the normalisation of $\T_N^{\new}$,
then we may descend $\rho$ to a two dimensional representation
defined over $\normalclosure{\T}_N^{\new}$.  If $\heightone$ is a classical
height one prime ideal in $\T_N^{\new}$,
then we have an isomorphism
$(\T_N^{\new})_{\heightone} \iso (\normalclosure{\T}_N^{\new})_{\heightone}$,
by Theorem~\ref{thm:new}, part~(3),
and consequently may also descend $\rho$ to a two dimensional
representation over $(\T_N^{\new})_{\heightone}$.
Reducing this representation modulo the maximal ideal of this local
ring, we obtain a representation
$\rhobar_{\heightone}: G_{\Q} \rightarrow \GL_2(\O(\heightone)[1/p]).$
Part~(2) of Theorem~\ref{thm:galois} shows that this is the usual
absolutely irreducible
two dimensional Galois representation attached to the newform $f_{\heightone}$.
We may also use an integral model of $\rho$ to construct a residual
representation $\rhobar_{\frak m}$ attached to any maximal ideal
of $\T_N^{\new}$.  We recall this construction in more detail
in Theorem~\ref{thm:residual} below.

If $\frak a$ is a minimal prime of $\T_N^{\new}$, then the representation
$\rho_{\frak a}$ takes values in a field, and so we may define its
tame conductor by the usual formula.  That is,  if $\rho_{\frak a}$ acts
on the two dimensional $(\T_N^{\new})_{\frak a}$-vector space $\mathcal V$, then
exponent of a prime $\ell \neq p$ in the tame conductor is equal to
$$ \dim {\mathcal V} - \dim {\mathcal V}_{I_{\ell}} + \int_0^{\infty}
\left(\dim {\mathcal V} - \dim {\mathcal V}_{I_{\ell}^u} \right)\, du,$$
where $\{I_{\ell}^u\}_{u\geq 0}$ denotes the usual filtration
of $I_{\ell}$ by higher ramification groups, indexed by the upper numbering.
(Note that this formula is often expressed in terms of invariants
under the ramification groups, rather than coinvariants.  
One obtains the same value with either formulation, since
$\dim {\mathcal V}_{I_{\ell}^u} = \dim {\mathcal V}^{I_{\ell}^u}$ for
any value of $u$.)

\begin{prop}\label{prop:conductor}
If $\frak a$ is any minimal prime of $\T_N^{\new}$,
then the tame conductor of $\rho_{\frak a}$ is equal to $N$.
\end{prop}
\begin{proof}
Write $\normalclosure{\T}$ to denote the normalisation
of $\T_N^{\new}/\frak a$; this is the component of the normalisation
of $\T_N^{\new}$ cut out by the minimal prime $\frak a$.
Since it is Cohen-Macaulay (being normal and of dimension two)
and finite over $\Iw$, it is in fact finite flat over $\Iw$ by
\cite[Thm.~23.1]{matsumura}.

Let $\mathcal V$ denote a two dimensional vector space over
the fraction field of $\T^{\new}_N/\frak a$ on which $\rho_{\frak a}$
acts,
and let $M$ be a free rank two $\normalclosure{\T}$-lattice
in $\mathcal V$ invariant under $G_{\Q}$.
If $\heightone$ is a classical height one prime of $\T_N^{\new}$
containing $\frak a$,
then as was observed above, the injection
$$(\T_N^{\new}/\frak a)_{\heightone} \rightarrow
\normalclosure{\T}_{\heightone}$$
is an isomorphism.  Moreover, 
$G_{\Q}$ acts on the $\O(\heightone)[1/p]$-vector space 
$V(\heightone):= (M/\heightone M) [1/p] $
via the usual Galois representations $\rhobar_{\heightone}$
attached to newform $f_{\heightone}$ of tame conductor $N$.

The ratio of the tame conductor of $\rho_{\frak a}$ and the tame conductor
of $\rhobar_{\heightone}$ (that is, $N$) is equal to
\begin{equation}\label{eqn:conductors}
\prod_{\ell \neq p}
\ell^{\,\dim V(\heightone)_{I_{\ell}} -\dim {\mathcal V}_{I_{\ell}}};
\end{equation}
furthermore, each of the exponents appearing in this expression is non-negative
(i.e.\ $\dim V(\heightone)_{I_{\ell}}$ bounds $\dim {\mathcal V}_{I_{\ell}}$
from above).
(See \cite[\S 1]{livne}).
Thus to prove the proposition, we must show that
$\dim V(\heightone)_{I_{\ell}} = \dim {\mathcal V}_{I_{\ell}}$
for each prime $\ell \neq p$.

If $\dim V(\heightone)_{I_{\ell}} = 2$, then $\ell$ does not
divide $N$.  In this case, $\rho_{\frak a}$ is unramified at $\ell$,
and so $\dim {\mathcal V}_{I_{\ell}} = 2$ as well.  
Conversely, if $\dim {\mathcal V}_{I_{\ell}} = 2$, then the same is true
of $\dim V(\heightone)_{I_{\ell}}$ (since the latter dimension bounds the former
dimension from above).
If $\dim V(\heightone)_{I_{\ell}} = 0$, then also $\dim {\mathcal V}_{I_{\ell}}=0$ 
(again, since the latter dimension is bounded above by the former).  
Thus it remains to show that if $\dim V(\heightone)_{I_{\ell}} > 0$,
then $\dim {\mathcal V}_{I_{\ell}} > 0$.  This we now do.

Note that the formula~(\ref{eqn:conductors}) (and the fact that
the tame conductor of $V(\heightone)$ is equal to $N$,
and thus is independent of the particular choice of $\heightone$ containing
$\a$) shows that
$\dim V(\heightone)_{I_{\ell}}$ is independent of the choice of the
classical height one prime $\heightone$ containing $\frak a$.
Thus we may assume that $\dim V(\heightone)_{I_{\ell}}$ is positive
for every classical height one prime $\heightone$ containing $\frak a$. 
There is a natural map $I_{\ell} \rightarrow \Z_p(1)$, given by projection
onto the $p$-Sylow subgroup of the tame quotient of $I_{\ell}$.
Let $J_{\ell}$ denote the kernel of this projection.
Since $J_{\ell}$ has order prime to $p$, it acts on $M$
through a finite
quotient of order prime to $p$, and so certainly
$\dim {\mathcal V}_{J_{\ell}} = \dim V(\heightone)_{J_{\ell}},$
for any classical height one prime $\heightone$.
Let $\sigma$ denote a topological generator of $\Z_p(1),$
and consider the matrix $\rho_{\frak a}(\sigma) - I$ acting on
${\mathcal V}_{J_{\ell}}$.  (Here $I$ denotes the identity matrix.)
By assumption, the determinant of this matrix (which lies in $\hat{\T}$)
lies in $\heightone\hat{\T}$ for every classical height one prime $\heightone$
of $\hat{\T}$.  Thus it lies in $\heightone' \hat{\T}$ for every
classical height one prime $\heightone'$ of $\Iw$.  (Recall that these
primes are unramified in $\hat{\T}$, by part~(3) of Theorem~\ref{thm:new}.)
Since $\hat{\T}$ is finite flat over $\Iw$, we see that this determinant
vanishes, and thus that $\rho_{\frak a}(\sigma)$ admits a non-zero 
coinvariant quotient of ${\mathcal V}_{J_{\ell}}$.   This proves that
$\dim {\mathcal V}_{I_{\ell}} > 0,$ as required.
\end{proof}

As a byproduct of the proof of the preceding proposition,
we also obtain the following useful result.

\begin{prop}\label{prop:coinvariants}  Let $\heightone$ be a classical
height one prime ideal in $\T_N^{\new}$,
and let $L$ denote a choice of two dimensional
$(\T_N^{\new})_{\heightone}$-lattice on which $\rho$ acts,
so that $L/\heightone L$ is the two dimensional $\O(\heightone)[1/p]$-vector
space on which $\rhobar_{\heightone}$ acts.
If $\ell$ is any prime distinct from $p$, then 
$L_{I_{\ell}}$ is a free $(\T_N^{\new})_{\heightone}$-module,
and there are natural isomorphisms 
$L_{I_{\ell}} \otimes_{\Iw} {\mathcal L}
\rightarrow (L\otimes_{\Iw} {\mathcal L})_{I_{\ell}}$ and
$L_{I_{\ell}}/\heightone L_{I_{\ell}} \rightarrow (L/\heightone L)_{I_{\ell}}.$
\end{prop}
\begin{proof}
We let $\frak a$ denote the minimal ideal contained in $\heightone$
(unique, since $(\T_N^{\new})_{\heightone}$ is a discrete valuation ring),
and employ the notation introduced in the proof of the preceding
proposition.  The lattice $L$ is a free rank two
$(\T_N^{\new})_{\heightone}$-module, which we may regard as being embedded
in a $G_{\Q}$-equivariant fashion in $\mathcal V$.

For general reasons,
the composite surjection $L \rightarrow L/\heightone L \rightarrow
(L/\heightone L)_{I_{\ell}}$ induces an isomorphism
$L_{I_{\ell}}/\heightone L_{I_{\ell}}
\iso (L/\heightone L)_{I_{\ell}}.$
Thus a minimal set of generators for $L_{I_{\ell}}$ as a
$(\T_N^{\new})_{\heightone}$-module contains at most
$\dim (L/\heightone L)_{I_{\ell}}$ elements.
On the other hand the surjection ${\mathcal V} \rightarrow {\mathcal V}_{I_{\ell}}$
induces a surjection
$L_{I_{\ell}} \otimes_{\Iw} \mathcal L \rightarrow {\mathcal V}_{I_{\ell}}$,
and hence $L_{I_{\ell}} \otimes_{\Iw} \mathcal L$ is of dimension
at least $\dim {\mathcal V}_{I_{\ell}}$.  
The proof of Proposition~\ref{prop:conductor}
shows that $\dim (L/\heightone L)_{I_{\ell}} = \dim {\mathcal V}_{I_{\ell}}$,
and hence 
(taking into account the fact that $(\T_N^{\new})_{\heightone}$
is a discrete valuation ring) we may conclude that $L_{I_{\ell}}$ is free
over $(\T_N^{\new})_{\heightone}$, of rank equal to
$\dim {\mathcal V}_{I_{\ell}}.$
Thus the surjection
$L_{I_{\ell}} \otimes_{\Iw} \mathcal L \rightarrow {\mathcal V}_{I_{\ell}}$
is an isomorphism, as claimed.
\end{proof}

The next result describes the residual representation
$\rhobar_{\frak m}$ attached to a maximal ideal of~$\T_N^{\new}$.

\begin{thm}\label{thm:residual}
If $\frak m$ is a maximal ideal of $\T_N^{\new}$, then attached to
$\frak m$ is a semi-simple representation
$\rhobar_{\frak m}: G_{\Q} \rightarrow \GL_2(\T_N^{\new}/\frak m),$
uniquely determined by the properties:
\begin{enumerate}
\item $\rhobar_{\frak m}$ is unramified away from $N p$.

\item If $\ell$ is a prime not dividing $N p$ then
$\rhobar_{\frak m}(\Frob_{\ell})$ has
characteristic polynomial equal to
$X^2 - T_{\ell} X + \langle \ell \rangle \ell^{-1} \bmod \frak m
\in (\T_N^{\new}/\frak m)[X].$
\end{enumerate}

Furthermore, $\rhobar_{\frak m}$ satisfies the following condition:

\begin{enumerate}
\setcounter{enumi}{2}
\item The restriction of $\rhobar_{\frak m}$ to $D_p$ has the following shape
(with respect to a suitable choice of basis):
$$\left( \begin{matrix} \chi & * \\ 0 & \psi \end{matrix}\right),$$
where $\chi$ and $\psi$ are $(\T_N^{\new}/\frak m)^{\times}$-valued
characters of $D_p$, such that $\psi$ is unramified and
$\psi(\Frob_p) = U_p \bmod \frak m$.  
\end{enumerate}
\end{thm}
\begin{proof}
The representation $\rhobar_{\frak m}$ is constructed in the usual
way, by choosing an integral model for $\rho$
over $\normalclosure{\T}_N^{\new},$
reducing this model modulo a maximal ideal $\normalclosure{\frak m}$
lifting $\frak m$,
semi-simplifying, and then descending (if necessary) from
$\normalclosure{\T}_N^{\new}/\normalclosure{\frak m}$ to $\T_N^{\new}/\frak m$.
The stated properties follow from the corresponding properties of $\rho$.
\end{proof}

\begin{prop}\label{prop:lattice}
If $\frak m$ is a maximal ideal of $\T_N^{\new}$
for which the associated residual representation $\rhobar_{\frak m}$
is irreducible,
then $\rho_{\frak m}$ admits a model over $(\T_N^{\new})_{\frak m}$
(which we denote by the same symbol)
$$\rho_{\frak m}: G_{\Q} \rightarrow \GL_2((\T_N^{\new})_{\frak m}),$$
unique up to isomorphism.
Furthermore, in this model, the space of $I_p$-coinvariants is
free of rank one over $(\T_N^{\new})_{\frak m}$.
\end{prop}
\begin{proof}   This follows from the irreducibility of $\rhobar_{\frak m}$,
and the fact that the traces of $\rho_{\frak m}$ lie in
$(\T_N^{\new})_{\frak m}$.  (See \cite{carayol}.)
\end{proof}

Suppose that $\frak m$ is a maximal ideal of $\T_N^{\new}$ satisfying
the hypothesis of the preceding proposition, and let $M$ denote a choice
of free rank two $(\T_N^{\new})_{\frak m}$-module on which the representation
$\rho_{\frak m}$ acts.  The group $G_{\Q}$ then acts on the quotient
$M/\frak m M$ via the residual representation $\rhobar_{\frak m}$.
Since the space $M_{I_p}$ of $I_p$-coinvariants of $M$ is a free rank
one quotient of $M$, we see that $(M_{I_p})/\frak m (M_{I_p})$ is
a one dimensional quotient of the $I_p$-coinvariants
$(M/\frak m M)_{I_p}$ of $M/\frak m M$.

\begin{definition}\label{def:p-stab} 
If $k$ is a finite field of characteristic $p$,
and $\rhobar: G_{\Q} \rightarrow \GL_2(k)$ is a continuous two dimensional
Galois representation defined over $k$, acting on the two dimensional
$k$-vector space $V$, say, then a {\em $p$-stabilisation}
of $\rhobar$ consists of the choice of a one dimensional quotient of
the space $V_{I_p}$ of $I_p$-coinvariants of $V$.
\end{definition}

\begin{definition}\label{def:canonical p-stab}
The discussion preceding Definition~\ref{def:p-stab} shows that
if $\frak m$ is a maximal ideal of $\T_N^{\new}$ for which 
$\rhobar_{\frak m}$ is irreducible,
then $\rhobar_{\frak m}$ comes equipped with a natural $p$-stabilisation.
We refer to this as the {\em canonical $p$-stabilisation} of $\rhobar_{\frak m}$
attached to the maximal ideal $\frak m$.
\end{definition}

\subsection{The reduced Hida algebra}\label{subsec:reduced}


\begin{definition}\label{def:reduced}
For any level $N$, we let $\T_N'$ denote the $\Iw$-subalgebra
of $\T_N$ generated by the Hecke operators $T_{\ell}$ for $\ell$
prime to $N p$, together with the operator $U_p$.
\end{definition}

Since $\T_N'$ is a subalgebra of the finite flat $\Iw$-algebra
$\T_N$, it is certainly finite and torsion free over $\Iw$.
It turns out that $\T_N'$ is also reduced.  (Being reduced
is a standard property of Hecke algebras in which
we omit the operators indexed by the primes dividing the level.)
In fact, we can be somewhat more precise.

If $M$ is a divisor of $N$, then restricting the action of
the prime-to-$N$ Hecke operators to ordinary forms of level
dividing $M$ yields a surjective map of $\Iw$-algebras
$\T_N' \rightarrow \T_M'$.  Composing this with
the composite $\T_M' \subset \T_M \rightarrow \T_M^{\new}$
yields a map
\begin{equation}\label{eqn:comparison}
\T_N' \rightarrow \T_M^{\new}.
\end{equation}
Taking the product of these maps over all divisors $M$ of $N$,
we obtain a map
\begin{equation}\label{eqn:embedding}
\T_N'\rightarrow \prod_{M | N} \T_M^{\new}.
\end{equation}

\begin{prop}\label{prop:embedding}
~
\begin{enumerate}
\item The map~(\ref{eqn:embedding}) is injective, and induces
an isomorphism after localising at any classical height one prime ideal of $\Iw$
(and hence also after tensoring over~$\Iw$ with its fraction field $\mathcal L$).

\item If $\heightone$ is a classical height one prime ideal of $\Iw$,
then $\T_N'\otimes_{\Iw} \Iw_{\heightone}$ is a finite \'etale
extension of the discrete valuation ring $\Iw_{\heightone}$.
\end{enumerate}
\end{prop}

\begin{proof}
It is an easy consequence of the theory of newforms for ordinary
families developed in \cite{hida1} that the map~(\ref{eqn:embedding})
is injective, and that it induces an isomorphism after tensoring
with $\mathcal L$ over $\Iw$.
We wish to check that the same is true if we tensor with $\Iw_{\heightone}$
over $\Iw$ for any classical height one prime $\heightone$ of $\Iw$.
For this, it suffices to replace $\T_N'$ and each of the $\Iw$-algebras
$\T_M^{\new}$ 
by their free closures (in the terminology of \cite{hida1}).
We denote these free closures by $\free{\T}_N'$ and
$\free{\T}_M^{\new}$
respectively.

It is proved in \cite{hida1}, and recalled in the proof of
Theorem~\ref{thm:new}, that the isomorphism of part~(2) of
Theorem~\ref{thm:hida} induces an isomorphism of
$(\free{\T}_M^{\new}/\heightone \free{\T}_M^{\new})[1/p]$
with the $\Q_p$-Hecke algebra acting on the space of newforms
$$
S_k(M p^{\infty}, \O(\heightone)[1/p])^{\ord}_{\new}[\kappa_{\heightone}].
$$
Strong multiplicity one for newforms, together with the standard
duality between spaces of modular forms and Hecke algebras,
then shows that the map
$$(\free{\T}_N'/\heightone \free{\T}_N')[1/p] \rightarrow
\prod_{M | N}
(\free{\T}_M^{\new}/\heightone \free{\T}_M^{\new})[1/p]$$
is surjective.
Thus the map 
$$ \free{\T}_N' \otimes_{\Iw} \Iw_{\heightone} \rightarrow
\prod_{M | N} \free{\T}_M^{\new} \otimes_{\Iw} \Iw_{\heightone}$$
is a surjection of finite free $\Iw_{\heightone}$-modules.
Since it induces an isomorphism after extending scalars to $\mathcal L$,
it must be an isomorphism.  This completes the proof of part~(1).  
Part~(2) is an immediate consequence of part~(1),
together with part~(3) of Theorem~\ref{thm:new}.
\end{proof}

Taking the product of
the Galois representations
$$G_{\Q} \rightarrow \GL_2(\T_M^{\new}\otimes_{\Iw} \mathcal L)$$
given by Theorem~\ref{thm:galois}, as $M$ ranges over all divisors of $N$,
and taking into account part~(1) of the preceding proposition,
we obtain a Galois representation
$$\rho: G_{\Q} \rightarrow \GL_2(\T_N'\otimes_{\Iw} \mathcal L)$$
satisfying the analogue of Theorem~\ref{thm:galois}.

Just as in the preceding section, we may reduce $\rho$ modulo a
height one prime ideal $\heightone$ or a maximal ideal $\frak m$
of $\T_N'$.  We denote the corresponding residual representations by 
$\rhobar_{\heightone}$ and $\rhobar_{\frak m}$ respectively.
Similarly, we may localise $\T'_N$ at any maximal ideal
$\frak m$, and obtain a corresponding representation
$$\rho_{\frak m} : G_{\Q} \rightarrow
\GL_2((\T'_N)_{\frak m} \otimes_{\Iw} \mathcal L).$$
If $\rhobar_{\frak m}$ is irreducible, then the analogue
of Proposition~\ref{prop:lattice} holds (by the same appeal
to the results of \cite{carayol}),
and so we obtain a uniquely determined representation
$$\rho_{\frak m}: G_{\Q} \rightarrow \GL_2((\T_N')_{\frak m}).$$

\subsection{Hecke eigenvalues}

Fix a tame level $N$,
a maximal ideal $\frak m$ of $\T_N'$,
a divisor $M$ of $N$,
and a maximal ideal $\frak n$ of $\T_M^{\new}$
that pulls back to $\frak m$ under the map~(\ref{eqn:comparison}),
and hence induces a local map
$(\T_N')_{\frak m}  \rightarrow (\T_M^{\new})_{\frak n}.$ 
We will let $\T$ denote the image of this map;
by construction $(\T_M^{\new})_{\frak n}$ is obtained from
$\T$ by adjoining the Hecke operators $T_{\ell}$ for those primes $\ell$
dividing $N$, together with the diamond operators corresponding
to elements of $(\Z/M)^{\times}$. 

The residual representations $\rhobar_{\frak m}$ and $\rhobar_{\frak n}$
are isomorphic, up to a possible extension of scalars; we assume
that they are irreducible (equivalently, absolutely irreducible,
since the residue characteristic $p$ is odd).   It follows
from the results of \cite{carayol} that
the Galois representation
$\rho_{\frak n}: G_{\Q} \rightarrow \GL_2((\T_M^{\new})_{\frak n})$
is in fact defined over $\T$.
If $\ell$ does not divide $M$, then $\rho_{\frak n}$
is unramified at $\ell$, and $T_{\ell}$ is equal to the trace
of $\rho_{\frak n}(\Frob_{\ell})$.   In particular, for such
primes $\ell$, the Hecke operator $T_{\ell}$ lies in $\T$.
The map $(\T_N')_{\frak m} \rightarrow (\T_M^{\new})_{\frak n}$
thus factors through a surjection
$(\T_N')_{\frak m} \rightarrow (\T_M')_{\frak m}$
(where we denote the maximal ideal in $\T_M'$
induced by $\frak m$ by the same letter).  Thus we may suppose
that $M = N$, and we do so from now on.
The diamond operators given by $(\Z/N)^{\times}$ also lie 
in $\T$, since these can be recovered by evaluating the determinant
of $\rho_{\frak n}$ on Frobenius elements.
Thus we see that $(\T_N^{\new})_{\frak n}$ is generated by adjoining
to $\T$ the Hecke operators $T_{\ell}$, for $\ell$ dividing $N$.
In the remainder of this
section we explain the representation-theoretic interpretation
of these Hecke operators,
and also consider the related question
of whether or not the inclusion
$\T \subset (\T_N^{\new})_{\frak n}$ is an equality.


Before proceeding further, it will be helpful to
recall the Galois representation-theoretic description
of the Hecke eigenvalues of a classical newform.
Thus we let $f$ denote a classical newform 
over $\Q_p$ of some weight,
and of tame level $N$,
let $\rho_f$ denote the $p$-adic Galois representation
attached to $f$,
and let $a_{\ell}(f)$ denote its Hecke eigenvalues.
Let $V$ be a two dimensional
vector space on which the $p$-adic Galois representation $\rho_f$
acts.  If $\ell$ does not divide $N$, then $\rho_f$ is unramified
at $\ell$, and $a_{\ell}(f)$ is equal to the trace of $\rho_f(\Frob_{\ell})$.
For each prime $\ell$ dividing $N$, the space $V_{I_{\ell}}$
of $I_{\ell}$-coinvariants in $V$ is either trivial or one dimensional.
In the first case, one has that $a_{\ell}(f) = 0$.  In the second
case, one has that $a_{\ell}(f)$ is equal to the eigenvalue
of $\Frob_{\ell}$ acting on $V_{I_{\ell}}$. 
If $\O(f)$
denotes the finite extension of $\Z_p$ generated by the Fourier
coefficients of $f$, then considerations of an integral model
for $\rho_f$ show that this eigenvalue is a unit
in ${\O(f)}$.  Thus we see that, for the primes $\ell$ dividing $N$,
if $a_{\ell}(f)$ lies in
the maximal ideal of ${\O(f)}$, then it actually vanishes.

We are now in a position to prove the following lemma.

\begin{lemma}\label{lem:ramified hecke}
If $\ell$ is a prime dividing $N$, then the following are equivalent:
\begin{enumerate}
\item $T_{\ell}$ lies in the maximal ideal $\frak n$.

\item There is a classical height one prime ideal $\heightone$
of $(\T_N^{\new})_{\frak n}$
such that $T_{\ell} \equiv 0 \bmod \heightone.$

\item For every classical height one prime ideal $\heightone$ of
$(\T_N^{\new})_{\frak n}$,
we have $T_{\ell} \equiv 0 \bmod \heightone.$

\item We have $T_{\ell} = 0$ in $(\T_N^{new})_{\frak n}$.

\item Letting $\mathcal V$ denote a free rank two
$(\T_N^{\new})_{\frak n} \otimes_{\Iw}
\mathcal L$-module on which $\rho_{\frak n}$ acts, we have
${\mathcal V}_{I_{\ell}} = 0.$
\end{enumerate}
If these equivalent conditions do not hold,
then we have that ${\mathcal V}_{I_{\ell}}$ is free of rank one
over $(\T_N^{\new})_{\frak n}\otimes_{\Iw} \mathcal L$,
and $T_{\ell}$ is equal to the eigenvalue of $\Frob_{\ell}$ acting
on this free rank one module.
\end{lemma}
\begin{proof}
It is clear that $\text{(4)} \implies \text{(3)}
\implies \text{(2)} \implies \text{(1)}.$
Condition~(1) implies
that $a_{\ell}(f_{\heightone})$ is contained in the maximal
ideal of ${\O(f_{\heightone})} = \O(\heightone)$
for each classical height one prime $\heightone$ of $(\T_N^{\new})_{\frak n}$,
and thus (by the discussion preceding the statement of the lemma)
that $a_{\ell}(f_{\heightone}) = 0$ for each such prime.
Thus~(1) implies~(3).
Proposition~\ref{prop:coinvariants}, together with the discussion
preceding the statement of the lemma, shows that~(3) and~(5)
are equivalent.  
Let $(\free{\T}^{\new}_N)_{\frak n}$ denote the free closure
of $(\T^{\new}_N)_{\frak n}$ as a $\Iw$-module.  If condition~(3)
holds, then (taking into account part~(3) of Theorem~\ref{thm:new})
we find that $$T_{\ell} \in \bigcap_{\heightone' \text{ classical ht.~1 in } 
\Iw} \heightone' (\free{\T}^{\new}_N)_{\frak n},$$
and thus that $T_{\ell} = 0$.  Thus~(3) implies~(4), and this completes
the proof of the claimed equivalences.

Suppose now that these equivalent conditions do not hold.
The remarks preceding the statement of the lemma, together with
Proposition~\ref{prop:coinvariants}, show that
in this case the space ${\mathcal V}_{I_{\ell}}$ is free of rank one
over $(\T_N^{\new})_{\frak n}\otimes_{\Iw} \mathcal L$, as claimed.
Let $A_{\ell} \in (\T_N^{\new})_{\frak n}\otimes_{\Iw} \mathcal L$
denote the eigenvalue of $\Frob_{\ell}$ acting on this module.
Again, Proposition~\ref{prop:coinvariants} shows that this
module admits a model over $(\T_N^{\new})_{\heightone}$,
for each classical height one prime ideal of $(\T_N^{\new})_{\frak n}$,
and thus that in fact $A_{\ell} \in (\T_N^{\new})_{\heightone}$
for each such $\heightone$.   The discussion preceding
the statement of the lemma, together with part~(3) of Theorem~\ref{thm:new},
then shows that 
$$
T_{\ell} - A_{\ell} \in 
\bigcap_{\heightone'  \text{ classical ht.~1 in } \Iw}
\heightone' (\free{\T}^{\new}_N)_{\frak n}\otimes_{\Iw} \Iw_{\heightone'}.
$$
This implies that $T_{\ell} - A_{\ell} = 0,$ and completes the
proof of the remaining claim of the lemma.
\end{proof}

Recall that $\T$ denotes the image of $(\T'_N)_{\frak m}$ in
$(\T_N^{\new})_{\frak n}$.  We now address the issue of whether
or not $\T$ equals all of $(\T_N^{\new})_{\frak n}$.
We note that the remainder of this section is tangential to the main
arguments of the paper.

\begin{prop}\label{prop:integral surjectivity}
Assume that $\rhobar_{\frak m}$ is irreducible
and that 
the characteristic polynomial of any element in the image
of $\rhobar_{\frak m}$ splits over $\T'_N/\frak m$
if and only it splits over $\T^{\new}_N/\frak n$.
Then for any prime $\ell$ prime to $p$ we have that $T_{\ell}$
lies in $\T$,
unless $\ell$ satisfies both of the following conditions:
\begin{enumerate}
\item $\rhobar_{\frak m}$ is unramified at $\ell$,
and the characteristic polynomial
of $\rhobar_{\frak m}(\Frob_{\ell})$ has repeated roots.

\item $\mathcal V_{I_{\ell}}$ is free of rank one.
\end{enumerate}
Furthermore, if these conditions hold,
then $\ell \equiv 1 \bmod p$.
\end{prop}

\begin{remark}
Note that the condition on characteristic polynomials above
holds if the map
$\T'_N/\frak m \rightarrow \T^{\new}_N/\frak n$ is an isomorphism,
or if the characteristic polynomial of every element
in the image of $\rhobar_{\frak m}$ already splits over
$\T'_N/\frak m$.
\end{remark}

\begin{proof}
Let $\overline{V}$ denote a two dimensional vector space
on which $\rhobar_{\frak m}$ acts.  If $N(\rhobar_{\frak m})$ denotes the
tame conductor of $\rhobar_{\frak m}$, then \cite[\S 1]{livne}
implies that
$$N/N(\rhobar_{\frak m}) =
\prod_{\ell \neq p} \ell^{\dim \overline{V}_{I_{\ell}}
- \dim {\mathcal V}_{I_{\ell}}}.$$
If $\dim {\mathcal V}_{I_{\ell}} = 0,$ then Lemma~\ref{lem:ramified hecke}
shows that $T_{\ell} = 0,$ and so there is nothing to prove.
If $\dim \overline{V}_{I_{\ell}} = \dim {\mathcal V}_{I_{\ell}},$
then arguments similar to those used in the proofs of Remarks~2.9
and~2.11 of \cite{wiles} prove that $T_{\ell} \in \T$.

It remains to consider primes $\ell$ for which
${\mathcal V}_{I_{\ell}}$ is one dimensional,
while $\rhobar_{\frak m}$ is unramified
at $\ell$.  In such a case, we see that $T_{\ell}$ reduces modulo $\frak n$
to one of the two roots of
the characteristic polynomial of $\rhobar_{\frak m}(\Frob_{\ell})$.
By hypothesis
this characteristic polynomial splits over 
the residue field $\T'_N/\frak m$.   If it has distinct roots,
then Hensel's lemma shows that $T_{\ell} \in \T$.
Thus if $T_{\ell} \not\in \T$, then the characteristic
polynomial of $\rhobar_{\frak m}(\Frob_{\ell})$ must have repeated roots.


Let us finally explain how 
to conclude that $\ell \equiv 1 \bmod p$.
If $\heightone$ is any classical height one prime of $\T$,
then recall that $\rhobar_{\heightone}$ is the $p$-adic Galois representation
attached to the classical newform $f_{\heightone}$ of level $N$.
Proposition~\ref{prop:coinvariants} shows that
$\rhobar_{\heightone}$ has one dimensional $I_{\ell}$-coinvariants.
It it then well-known that there is either an isomorphism
$$(\rhobar_{\heightone \, | D_{\ell}})^{\semisimp}
\iso \psi \epsilon \oplus \psi,$$
where $\psi$ denotes an unramified character and $\epsilon$ the $p$-adic
cyclotomic character (the special case),
or
$$(\rhobar_{\heightone \, | D_{\ell}})^{\semisimp}
\iso \chi \oplus \psi$$
where $\chi$ is tamely ramified and $\psi$ is unramified (the principal
series case).  In the first case, since $\rhobar_{\frak m}(\Frob_{\ell})$
acts on $(\rhobar_{\heightone \, | D_{\ell}})^{\semisimp}$
as a scalar, we see that $\epsilon(\Frob_{\ell}) = \ell \equiv 1
\bmod p$.  In the second case, since $\rhobar_{\frak m}$ is unramified at
$\ell,$ we see that $\chi_{| I_{\ell}} \equiv 1 \bmod p$;
the group $D_{\ell}$ admits
a ramified character with this property only if $\ell \equiv 1 \bmod p$.
\end{proof}

If we are willing to invert $p$, then we obtain the following
more definitive result.

\begin{prop}\label{prop:surjectivity}
The inclusion
$\T[1/p] \subset (\T_N^{\new})_{\frak n}[1/p]$
is an equality.
\end{prop}
\begin{proof}
As in the proof of the preceding proposition,
it suffices to show that if $\ell$ is a prime at which
$\rhobar_{\frak m}$ is unramified and for which
${\mathcal V}_{I_{\ell}}$ is one dimensional,
then $T_{\ell}$ lies in $\T[1/p]$.
We may factor $\T[1/p]$ as a product $\T_{\spe} \times \T_{\ps},$
where $\T_{\spe}$ is the factor on which the diamond operators at
$\ell$ act trivially, and $\T_{\ps}$ is a the factor on which these
diamond operators act non-trivially.   (The subscripts are for
``special'' and ``principal series'', respectively.)

By specialising at classical height one primes, and arguing
as in the proof of Lemma~\ref{lem:ramified hecke}, we find
that over $\T_{\ps}$,
the associated two dimensional Galois representation has a free rank one 
space of $I_{\ell}$-invariants, and that $T_{\ell}$ is equal to
the eigenvalue of $\Frob_{\ell}$ acting on this one dimensional space
of invariants.  In particular, we find that $T_{\ell} \in \T_{\ps}$.
A similar specialisation argument shows that
$T_{\ell}^2 = \ell^{-1} \det \rho(\Frob_{\ell})$ in $\T_{\spe}$.
An application of Hensel's lemma then shows that $T_{\ell}$ lies in
the image of $\T$ in $\T_{\spe}$, and so in particular
in $\T_{\spe}$.
\end{proof}

%

\subsection{The reduced Hida algebras attached to $\rhob$}\label{subsec:family}

If $k$ is a finite field of characteristic $p$,
and $\rhobar: G_{\Q} \rightarrow \GL_2(k)$ is a continuous two dimensional
Galois representation defined over $k$, then we say that $\rhobar$ is
ordinary if $\rhobar$ satisfies condition~(3) of Theorem~\ref{thm:residual}.
We say that $\rhobar$ is $p$-distinguished if furthermore
the characters $\chi$ and $\psi$ appearing in the statement of
condition~(3) are distinct.
Clearly $\rhobar$ admits a $p$-stabilisation, in the sense of
Definition~\ref{def:p-stab},
if and only if $\rhobar$ is ordinary.  
If $\rhobar$ is furthermore $p$-distinguished,
then $\rhobar$ admits at most two choices of $p$-stabilisation
(and does admit two such choices precisely when the determinant of $\rhobar$
is unramified at $p$).

Let us now fix such a representation
$\rhobar: G_{\Q} \rightarrow \GL_2(k)$, 
and let $\overline{V}$ be a two dimensional $k$-vector space
on which $\rhobar$ acts.
We assume that $\rhobar$ is
irreducible, odd, ordinary
and $p$-distinguished, and we fix a choice of $p$-stabilisation of $\rhobar$.
We assume that $k$ is equal to the
field generated by the traces of $\rhobar$.  (If it were not,
we could replace $k$ by this field of traces, and descend $\rhobar$
to the smaller field.)
Finally, we suppose that $\rhobar$ is modular (i.e.\ that it arises
as the residual representation attached to a modular form of some
weight and level defined over $\Qbar_p$).
Our goal in this section is to define the reduced ordinary Hecke algebras
$\T_{\Sigma}(\rhobar)$ attached to an ordinary residual representation,
and to describe its basic properties.

We let $N(\rhobar)$ denote the tame conductor of $\rhobar$.
If $\ell \neq p$ is prime, then define
$$m_{\ell} = \dim_k \overline{V}_{I_{\ell}},$$
and for any finite set of primes $\Sigma$ that does not contain $p$, 
write
$$N(\Sigma) = N(\rhobar) \prod_{\ell \in \Sigma} \ell^{m_{\ell}}.$$

\begin{thm}
There is a unique maximal ideal $\frak m$
of $\T_{N(\Sigma)}'$ such that
$\rhobar_{\frak m}$, with its canonical $p$-stabilisation,
is isomorphic to $\rhobar$, with its given $p$-stabilisation.
\end{thm}
\begin{proof}  The uniqueness is clear.  The existence follows
from the assumption that $\rhobar$ is modular, and the results 
of \cite{diamond1}.
\end{proof}

\begin{prop}\label{prop:reduced and full}
If $\frak m$ denotes the maximal ideal of $\T_{N(\Sigma)}'$ of
the preceding theorem,
then there is a unique maximal ideal $\frak n$ of $\T_{N(\Sigma)}$
satisfying the following conditions:
\begin{enumerate}
\item $\frak n$ lifts $\frak m$.

\item $T_{\ell} \in \frak n$ for each $\ell \in \Sigma$.

\item The natural map of localisations $(\T_{N(\Sigma)}')_{\frak m}
\rightarrow (\T_{N(\Sigma)})_{\frak n}$ is an isomorphism of $\Iw$-algebras.
\end{enumerate}
In particular, $(\T_{N(\Sigma)}')_{\frak m}$ is a finite flat $\Iw$-algebra.
Also, the image of $T_{\ell}$ in the localisation $(\T_{N(\Sigma)})_{\frak n}$ 
in fact vanishes for each $\ell \in \Sigma$.
\end{prop}
\begin{proof}  This is a variant of~\cite[Prop.~2.1.5]{wiles},
and is proved in an analogous manner.
The second to last claim follows from the rest of the proposition,
together with part~(1) of
Theorem~\ref{thm:hida}.
\end{proof}

\begin{definition} We let $\T_{\Sigma}(\rhobar)$ (or simply $\T_{\Sigma},$
if $\rhobar$ is understood) denote the localisation of
$\T_{N(\Sigma)}'$ at the maximal ideal whose existence is 
guaranteed by Proposition~\ref{prop:reduced and full}.
We let $\rho_{\Sigma}:G_{\Q} \rightarrow \GL_2(\T_{\Sigma})$
denote the Galois representation attached to this local factor
of $\T_{N(\Sigma)}',$ as discussed in section~\ref{subsec:reduced}.
Recall that $\rho_{\Sigma}$ is characterised by the following
property:
If~$\ell$ is a prime not dividing $N(\Sigma) p$, then $\rho_{\Sigma}(\Frob_{\ell})$ has trace equal to $T_{\ell}  \in \T_{\Sigma}$.
\end{definition}

Taking into account Proposition~\ref{prop:reduced and full},
we see that $\T_{\Sigma}$
is a reduced and finite flat $\Iw$-algebra.
Note that if $\Sigma \subset \Sigma'$ then $N(\Sigma) \mid N(\Sigma'),$
and the natural surjection $\T_{N(\Sigma')}' \rightarrow \T_{N(\Sigma)}'$
induces a surjection $\T_{\Sigma'} \rightarrow \T_{\Sigma}.$
The Galois representations $\rho_{\Sigma'}$ and $\rho_{\Sigma}$ are 
evidently compatible with this surjection.

We refer to $\Spec \T_{\Sigma}$ as the universal ordinary family
of newforms, or sometimes simply ``the Hida family'',
minimally ramified outside $\Sigma$, attached to $\rhobar$ and our
chosen $p$-stabilisation.
Localising $\rho_{\Sigma}$ over $\Spec \T_{\Sigma}$, we obtain
a two dimensional vector bundle on which $G_{\Q}$ acts, which we
refer to as the universal family of Galois representations over
$\Spec \T_{\Sigma}$.
If $\Sigma \subset \Sigma'$, then the surjection $\T_{\Sigma} \rightarrow
\T_{\Sigma'}$ induces a closed embedding
$\Spec \T_{\Sigma} \rightarrow \Spec \T_{\Sigma'}$, and the universal family
of Galois representations on the target pulls back to the universal family
of Galois representations on the source.  If we consider all $\Sigma$
simultaneously, then we obtain an Ind-scheme
$$\formalilim{\Sigma} \Spec \T_{\Sigma},$$
and a family of two dimensional
Galois representations lying over it.  We refer to this Ind-scheme
as the universal ordinary family of newforms, or simply ``the Hida family'',
attached to $\rhobar$ and its chosen $p$-stabilisation.

The Hida family corresponding to $\Sigma = \emptyset$ will play a special
role; we refer to it as the minimal Hida family attached to $\rhobar$
and our chosen $p$-stabilisation.

When $\rhob$ is irreducible after restriction to
the quadratic field of discriminant $\pm p$,
the results of Wiles and Taylor--Wiles \cite{wiles, taylor-wiles},
as strengthened by Diamond \cite{diamond2}, in fact allow us to 
identify the rings $\T_{\Sigma}$, and their accompanying Galois
representations $\rho_{\Sigma}$, with certain universal deformation
rings, and their accompanying universal Galois representations,
attached to the residual representation $\rhobar$.  Specifically,
let $R_{\Sigma}$ denote the universal deformation ring parameterizing
lifts of $\rhob$ which are ordinary at $p$, and whose tame conductor
coincides with that of $\rhob$ at primes not in $\Sigma \cup \{ p \}$.
The representation $\rho_{\Sigma}$ induces a map
$R_{\Sigma} \to \T_{\Sigma}$
which by \cite{diamond2} is 
an isomorphism after tensoring with the quotient
of $\Lambda$ by any classical height one prime; it follows that in fact
$R_{\Sigma} \to \T_{\Sigma}$ is an isomorphism, as claimed.

\subsection{Branches}
\label{sec:branches}

We now prove some results concerning the irreducible components
of the Hida family attached to $\rhobar$.

\begin{definition}
If $\frak a$ is a minimal prime ideal in $\T_{\Sigma}$,
for any $\Sigma$ as above, then we will write
$\T(\frak a) := \T_{\Sigma}/\frak a$.
Note that $\T(\frak a)$ is a local domain, finite over $\Iw$.
We will write $\K(\frak a)$ to denote the fraction field
of $\T(\frak a)$.  (Thus there is an isomorphism
$\K(\frak a) \iso \T(\frak a) \otimes_{\Iw} \mathcal L$.)
We let $\rho(\frak a)$ denote the Galois representation
$$\rho(\frak a): G_{\Q} \rightarrow \GL_2(\T(\frak a))$$
induced by $\rho_{\Sigma}$.  
\end{definition}

\begin{prop}
\label{prop:components}
If $\frak a$ is a minimal prime of $\T_{\Sigma}$, for any $\Sigma$
as above, there is a unique divisor $M$ of $N(\Sigma)$ and
a unique minimal prime $\a' \subseteq \T_M^{\new}$ sitting over $\a$ such that 
$$
\xymatrix{
{\T_\Sigma} \ar[r] \ar[d] & {\T'_{N(\Sigma)}}  \ar[r] &
{\prod_{M \mid N(\Sigma)} \T_M^{\new}} \ar[d] \\
{\T_\Sigma/\a} \ar@{->}^{=}[r] & {\T(\a)} \ar[r] & {\T_M^{\new}/\a'}}
$$
commutes.
\end{prop}

\begin{proof}
Since $\T_{\Sigma}$ is finite over $\Iw,$
the minimal primes of $\T_{\Sigma}$ are in bijection with
the local components of $\T_{\Sigma} \otimes_{\Iw} \mathcal L$.
Since $\T_{\Sigma}$ is a local factor of $\T'_{N(\Sigma)}$, 
these local components are included in the local components of
$\T'_{N(\Sigma)} \otimes_{\Iw} \mathcal L$.  By part~(1) of
Proposition~\ref{prop:embedding}, we have that 
$$
\T'_{N(\Sigma)} \otimes_{\Iw} \mathcal L \cong 
\prod_{M \mid N(\Sigma)} \T_M^{\new} \otimes_{\Iw} \mathcal L.
$$
The local components of $\prod_{M \mid 
N(\Sigma)} \T_M^{\new} \otimes_{\Iw} \mathcal L$
are in one-to-one correspondence with its minimal primes.  Thus, our given
minimal prime $\a$ gives rise to a minimal prime of this ring.  However, any such
minimal prime corresponds to a minimal prime $\a'$ in $\T_M^{\new}$ for some $M \mid N(\Sigma)$,
which establishes the proposition.
\end{proof}

\begin{definition}\label{def:components}
In the context of the preceding proposition, we refer to $M$ as
the tame conductor attached to $\frak a$
(or to the irreducible component of $\Spec \T_{\Sigma}$ corresponding
to $\frak a$), and write $N(\frak a)$ for $M$. 

We write $\full{\T(\frak a)} := \T_{N(\frak a)}^{\new}/\frak a'$; 
Proposition~\ref{prop:components} then gives rise to an embedding of local domains
$\T(\frak a) \rightarrow \full{\T(\frak a)}.$
(To see that $\full{\T(\frak a)}$ is local, note that it is a complete
finite $\Iw$-algebra, and hence a product of local rings.  Being
a domain, it must be local.)
\end{definition}

If $\frak a$ is a minimal prime of $\T_{\Sigma}$ for
which the tame conductor $N(\frak a)$ is minimal (i.e.\ equal
to $N(\rhobar)$), then the embedding
$\T(\frak a) \rightarrow \full{\T(\frak a)}$ is an isomorphism.
(This follows from the arguments used in the proofs of
Remarks~2.9 and~2.11 of \cite{wiles}.)
More generally, Proposition~\ref{prop:surjectivity}
shows that the induced embedding $\T(\frak a)[1/p]
\rightarrow \full{\T(\frak a)}[1/p]$ is an isomorphism.
In particular,
the embedding $\T(\frak a) \rightarrow \full{\T(\frak a)}$ induces
an isomorphism on fraction fields, and we will use this isomorphism to identify
the fraction field of $\full{\T(\frak a)}$ with $\K(\frak a)$.
One can thus think of $\full{\T(\frak a)}$ as a partial
normalisation of $\T(\frak a)$ inside $\K(\frak a)$,
which coincides with $\T(\frak a)$ after inverting $p$.

We next observe that at the generic point of the $\frak a$-component
of $\Spec \T_{\Sigma},$ the universal Galois representation has
conductor equal to the conductor of $\frak a$.

\begin{cor}\label{cor:conductor}
The Galois representation $\rho(\frak a)$
has tame conductor
equal to $N(\frak a)$.  (Here we define the tame conductor by
regarding $\rho(\frak a)$ as a Galois representation defined over
the field of fractions $\K(\frak a)$ of $\T(\frak a)$.)
\end{cor}
\begin{proof}
This follows from Propositions~\ref{prop:conductor}
and~\ref{prop:components}.
\end{proof}

The next result deals with the classical height one primes in
$\T_{\Sigma}$.

\begin{prop}\label{prop:classical}
Let $\heightone$ be any classical height one prime ideal of
$\T_{\Sigma}$.

\begin{enumerate}
\item The ring $\T_{\Sigma}$ is \'etale over $\Iw$ (and so
regular) in a neighbourhood
of $\heightone$; 
consequently $\heightone$ contains a unique minimal
prime $\frak a$ of $\T_{\Sigma}$,
and the natural map of localisations
$(\T_{\Sigma})_{\heightone} \rightarrow \T(\frak a)_{\heightone}$
is an isomorphism.

\item
Thinking of $\heightone$ as a height one prime of $\T(\frak a)$,
the map $\T(\frak a) \rightarrow \full{\T(\frak a)}$ is an isomorphism
in a neighbourhood of $\heightone$.
Consequently,
there is a unique height one prime $\heightone'$ of $\full{\T(\frak a)}$
that pulls back to $\heightone$ under the this map,
and the map of localisations $\T(\frak a)_{\heightone} \rightarrow
\full{\T(\frak a)}_{\heightone'}$ is an isomorphism.
\end{enumerate}
\end{prop}

\begin{proof}
Both claims follow directly from
Proposition~\ref{prop:embedding}.
(Part~(2) also follows from the fact that $p \notin \heightone$,
and the equality $\T(\frak a)[1/p] = \full{\T(\frak a)}[1/p].$)
\end{proof}

In the situation
of the preceding proposition,
we write $\full{\O(\heightone)}:= \O(\heightone')$;  
this is a finite extension of $\O(\heightone)$.
Recall that we have defined
the classical newform $f_{\heightone'}$ attached to the height one prime
ideal $\heightone'$ of $\full{\T(\frak a)}$ (thought of as
a classical height one prime of $\T_{N(\frak a)}^{\new}$).  We write
$f_{\heightone} := f_{\heightone'};$
part~(2) of Theorem~\ref{thm:new} 
implies that $f_{\heightone}$ lies in
$S_k(N(\frak a) p^{\infty}, \full{\O(\heightone)})^{\ord}_{\new}
[\kappa_{\heightone}].$

\subsection{$\Lambda$-adic modular forms and Euler factors}
\label{subsec:euler}

For each minimal prime $\frak a$ of $\T_{\Sigma}$,
we can define a formal $q$-expansion along the
partial normalisation $\Spec \full{\T(\frak a)}$ of the component
$\Spec \T(\frak a)$ of $\Spec \T_{\Sigma}$ that interpolates
the $q$-expansions of the newforms $f_{\heightone}$ obtained from the
classical primes $\heightone$ in $\Spec \T(\frak a)$.
Namely, if we write $\full{T(\frak a)} =
\T_{N(\frak a)}^{\new}/\frak a'$, as in Definition~\ref{def:components},
then we define
$f(\frak a,q) \in \full{\T(\frak a)}[[q]]$ via 
$$
f(\frak a,q) = \sum_{n \geq 1} (T_n  \bmod \frak a') q^n.
$$
(Here we have written $T_n$ rather than $U_p^r T_{n'},$ when
$n$ is of the form $n = n' p^r$ with $(n',p) = 1$, for the sake
of uniformity of notation.)

An alternative way of describing this formal $q$-expansion along
$\Spec \full{\T(\frak a)}$ is to describe the corresponding
Euler factors.

\begin{definition}\label{def:euler}
Let $\frak a$ be a minimal prime of $\T_{\Sigma}$.
As above, write $\full{\T(\frak a)} \iso \T_{N(\frak a)}^{\new}/\frak a'$,
For each prime $\ell \neq p$, define the reciprocal
Euler factor $E_{\ell}(\frak a,X) \in \full{\T(\frak a)}[X]$
via the usual formula:
$$E_{\ell}(\frak a,X) := \begin{cases} \text{ 
$(1 - (T_{\ell} \bmod \frak a')X + \langle \ell \rangle \ell^{-1} X^2)$
if $\ell$ is prime to $N(\frak a)$} \\
\text{$(1 - (T_{\ell} \bmod \frak a')X)$
otherwise}.
\end{cases}$$
For the sake of completeness, 
define $E_p(\frak a,X) = (1 - (U_p \bmod \frak a') X)$.
\end{definition}

In terms of these reciprocal Euler factors,
the formal $q$-expansion $f(\frak a,q)$ is characterised by the fact 
that its ``formal Mellin transform'' is equal to the formal Dirichlet series
$$\prod_{\ell} E_{\ell}(\frak a, \ell^{-s})^{-1}.$$

We will now give a Galois theoretic description of these
reciprocal Euler factors.  


\begin{prop}\label{prop:euler}
Let $\frak a$ be a minimal prime of $\T_{\Sigma}$,
and let $\mathcal V$ denote a two dimensional $\K(\frak a)$-vector space
on which the Galois representation $\rho(\frak a)$ acts.
If $\ell \neq p$ is prime,
then the Euler factor
$E_{\ell}(\frak a,X) \in \K(\frak a)[X]$ is equal to the determinant
$\det(\Id - \Frob_{\ell} X \, | \, {\mathcal V}_{I_{\ell}}).$
(Here ${\mathcal V}_{I_{\ell}}$ denotes the space
of $I_{\ell}$-coinvariants in $\mathcal V$.)
\end{prop}
\begin{proof}
This follows from Lemma~\ref{lem:ramified hecke}.
\end{proof}

Since $\full{\T(\frak a)} \subset \T(\frak a)[1/p]$,
we may in particular regard the reciprocal Euler factors
$E_{\ell}(\frak a,X)$ and the formal $q$-expansion $f(\frak a,q)$
as varying over $\Spec \T(\frak a)[1/p]$.
Let us close this section by signaling a phenomena which will be fundamental
to all that follows: the reciprocal Euler factors $E_{\ell}(\frak a,X)$
(or equivalently the formal $q$-expansions $f(\frak a,q)$) do not extend
in a well-defined fashion over $\Spec \T_{\Sigma}[1/p]$.  In general, 
if $\heightone$ is a (necessarily non-classical) height one prime lying
in the intersection of two different components of $\Spec \T_{\Sigma}[1/p],$
say $\Spec \T(\frak a_1)[1/p]$ and $\Spec \T(\frak a_2)[1/p],$ 
then $E_{\ell}(\frak a_1,X)$ and $E_{\ell}(\frak a_2,X)$ (and hence
$f(\frak a_1,q)$ and $f(\frak a_2,q)$)
may have different specialisations in $\O(\heightone)[1/p][X]$ (respectively
$\O(\heightone)[1/p][[q]]$) for certain values
of $\ell$ dividing $N(\frak a_1)$ or $N(\frak a_2)$.

\section{Algebraic Iwasawa invariants}
\label{sec:algebraic}
 
\subsection{Selmer groups of modular forms} \label{sec:selmer}

Let $f = \sum a_{n}q^{n}$ be a $p$-ordinary and $p$-stabilized
newform of weight $k \geq 2$, tame level $N$, and character $\chi$.
Let $K$ denote the finite extension of $\Qp$ generated by the
Fourier coefficients of $f$ and let $\O$ denote the ring of integers of
$K$; we write $k$ for the residue field and fix also a uniformizer $\pi$
of $\O$.  Let
$$\rho_{f} : \GQ \to \GL_{2}(K)$$
be the corresponding Galois representation, characterized by the fact
that the characteristic polynomial under $\rho_{f}$
of an arithmetic Frobenius at a prime $\ell \nmid Np$ is
$$X^{2} - a_{\ell}X + \chi(\ell)\ell^{k-1}.$$
By \cite[Thm.\ 2.6]{GS}
the restriction of $\rho_{f}$ to $G_{p}$ is of the form
\begin{equation} \label{eq:ordes}
\rho|_{G_{p}} \cong \left( \begin{array}{cc}
\vep^{k-1}\chi\varphi^{-1} & * \\ 0 & \varphi \end{array}\right)
\end{equation}
with $\varphi : \Gp \to \O^{\times}$ the unramified character sending
an arithmetic Frobenius to $a_{p}$.  

We assume that the semisimple residual representation
$$\bar{\rho}_{f} : \GQ \to \GL_{2}(k)$$
is absolutely irreducible.  It follows that, up to conjugation by
$\GL_{2}(\O)$, there is a unique integral model
$$\rho_{f} : \GQ \to \GL_{2}(\O)$$
of $\rho_{f}$, which we now fix.  For $0 \leq i \leq p-2$ 
let $\Afi$ denote a cofree $\O$-module of corank $2$ with
$\GQ$-action via $\rho_{f} \otimes \omega^{i}$.  
We obtain from (\ref{eq:ordes}) and \cite[Prop. 12.1]{Gross}
an $\O[\Gp]$-equivariant exact sequence
\begin{equation} \label{eq:prefund}
0 \to (K/\O)(\vep^{k-1}\chi\omega^{i}\varphi^{-1}) \to \Afi \to 
(K/\O)(\omega^{i}\varphi) \to 0.
\end{equation}
We write $\Afi'$ (resp.\ $\Afi''$) for the submodule (resp.\ quotient
module) of $\Afi$ in the above sequence.

For a place $v$ of $\Qi$ define
$$\Hs^{1}(\Qiv,\Afi) = \begin{cases}
H^{1}(\Qiv,\Afi) & v \neq v_{p}; \\
\im \bigl( H^{1}(\Qivp,\Afi) \to H^{1}(\Ivp,\Afi'') \bigr) & 
v = v_{p}. \end{cases}$$
Following \cite{Greenberg1}, 
the {\it Selmer group} of $\Afi$ is defined by
\begin{align*}
\Sel(\Qi,\Afi) &= \ker \bigl( H^{1}(\Qi,\Afi) \to \prod_{v} 
\Hs^{1}(\Qiv,\Afi) \bigr) \\
&= \ker \bigl( H^{1}(\QS/\Qi,\Afi) \to \prod_{v \in \Sigma} 
\Hs^{1}(\Qiv,\Afi) \bigr)
\end{align*}
for any finite set of places $\Sigma$ containing $v_{p}$, all archimedean
places, and all places dividing $N$.  We regard
$\Sel(\Qi,\Afi)$ as a $\Lambda_{\O}$-module via the natural action
of $\Gamma$.

Recall that $\mu^{\alg}(f,\omega^{i})$ (resp.\ 
$\lalg(f,\omega^{i})$) is defined to be the
largest power of $\pi$ dividing (resp.\ the number of zeroes of) 
the characteristic power series of the $\Lambda_{\O}$-dual of $\Sel(\Qi,\Afi)$
(assuming that this Selmer group is $\Lambda_{\O}$-cotorsion).

\begin{thm} \label{thm:kato}
Let $f$ be a $p$-ordinary and 
$p$-stabilized newform with $\rhob_{f}$ absolutely irreducible.  Then
$\Sel(\Qi,\Afi)$ is co-finitely generated, $\Lambda_\O$-cotorsion,
and has no proper $\Lambda_\O$-submodules of finite index.
Furthermore, $\mu^{\alg}(f,\omega^{i})$ vanishes if and only if
$\Sel(\Qi,\Afi)[\pi]$ is finite.  If this is the case, then
$\Sel(\Qi,\Afi)$ is $\O$-divisible and
$$\lalg(f,\omega^{i}) = \dim_{k} \Sel(\Qi,\Afi)[\pi].$$
\end{thm}
\begin{proof}
It is shown in \cite[Prop.\ 6]{Greenberg1} that
Selmer groups are always co-finitely generated.  The fact that
they are also $\Lambda_\O$-cotorsion for modular forms is proven in \cite{KKT}.
The equivalence of the vanishing of $\mu^{\alg}(f,\omega^{i})$ and
the finiteness of $\Sel(\Qi,\Afi)[\pi]$ is now an immediate consequence
of the structure theory of $\Lambda_\O$-modules.

The proof of \cite[Prop.\ 4.14]{Greenberg2} easily adapts
to show that $\Sel(\Qi,\Afi)$ has no proper $\Lambda_\O$-submodules of
finite index.  (As $\Afi$ need not be self-dual, this also requires
the fact that $\Sel(\Qi,\Afi^*)$ is $\Lambda_\O$-cotorsion,
which follows since $\Afi^{*}$ is also modular.)
When $\Sel(\Qi,\Afi)[\pi]$ is finite, the maximal $\O$-divisible
$\Lambda_\O$-submodule of $\Sel(\Qi,\Afi)$ has finite index, so that it
must coincide with $\Sel(\Qi,\Afi)$; that is,
$\Sel(\Qi,\Afi)$ is divisible.  It now follows again
from the structure theory of $\Lambda_\O$-modules that as $\O$-modules
$$\Sel(\Qi,\Afi) \cong (K/\O)^{\lalg(f,\omega^{i})},$$
which proves the last statement.
\end{proof}

We close this section with a useful result on the local invariants
$H^{0}(\Gv,\Afi)$ for places $v$ dividing primes $\ell \neq p$.
Let $\cond_{\ell}(\rhob_{f})$
denote the exponent of the highest power of $\ell$ that divides
the conductor of $\rhob_{f}$.

\begin{lemma} \label{lemma:locdiv}
Let $v$ be a place of $\Qi$ dividing a prime $\ell \neq p$.
If $\cond_{\ell}(\rhob_{f}) = \ord_{\ell}(N)$,
then $H^{0}(\Gv,\Afi)$ is $\O$-divisible for all $i$.
\end{lemma}
\begin{proof}
By the invariance of the Swan conductor under
reduction (see \cite[\S 1]{livne}), we have
$$\cond_{\ell}(\rho_{f}) - \cond_{\ell}(\rhob_{f}) =
\dim_{k} A_f[\pi]^{\Il} - \dim_{K} V_{f}^{\Il}$$
where $V_{f}$ is a two dimensional $K$-vector space with $\GQ$-action
via $\rho_{f}$.  Since $\cond_{\ell}(\rho_{f})$ equals $\ord_{\ell}(N)$
by \cite{Carayol1}, we see that
the hypothesis of the lemma is equivalent
to the equality
$$\dim_{k} A_{f}[\pi]^{\Il} = \dim_{K} V_{f}^{\Il}.$$
It follows easily from this that $A_{f}^{\Il} = A_{f}^{\Iv}$ is $\O$-divisible.
As the $\Gv/\Iv$-invariants of an
$\O$-divisible $\Gv/\Iv$-module are again $\O$-divisible, we conclude that
$H^{0}(\Gv,A_{f})$ is divisible, as claimed.
Since $\omega$ is unramified at $\ell$, the above argument
works for $\Afi$ as well.
\end{proof}

\subsection{Residual Selmer groups}

Let $\Ab$ denote a two-dimensional $k$-vector space equipped with
a continuous $k$-linear action of $\GQ$ and a $k[\Gp]$-equivariant
exact sequence
\begin{equation} \label{eq:prefund2}
0 \to \Ab' \to \Ab \to \Ab'' \to 0
\end{equation}
with $\Ab'$ and $\Ab''$ one-dimensional.
We make the following assumptions on this data.
\begin{enumerate}
\item $\Ab$ is absolutely irreducible as a $k[\GQ]$-module;
\item The $k$-vector space $\Ab''$ is unramified
as a $\Gp$-module;
\item $\Ab$ is {\it $p$-distinguished} in the sense that
the representation of $\Gp$ on $\Ab$ is non-scalar;
\item $\Ab$ is {\it modular}: there exists a totally ramified
extension $K'$ of $K$, a $p$-stabilized newform
$f \in \O'[[q]]$ of weight $k\geq 2$, and a $k[\GQ]$-isomorphism
$\Ab \cong A_{f}[\pi']$
identifying $\Ab'$ with $A_{f}'[\pi']$ in the notation of the previous
section; here $\pi'$ is a uniformizer of the ring 
of integers $\O'$ of $K'$.
\end{enumerate}
We will study Selmer groups of $\Ab$ and its cyclotomic twists
$\Ab \otimes \omega^{i}$ with respect to various local conditions.

\begin{definition}
A {\it finite/singular structure} $\SS$ on $\Ab \otimes \omega^{i}$
is a choice of $k$-subspaces
\begin{equation} \label{eq:finitedef}
\HfS^{1}(\Qiv,\Ab \otimes \omega^{i}) \inj H^{1}(\Qiv,\Ab \otimes \omega^{i})
\end{equation}
for each place $v$ of $\Qi$, subject to the restrictions:
\begin{enumerate}
\item $\HfS^{1}(\Qiv,\Ab \otimes \omega^{i}) = 0$ for almost all $v$;
\item $\HfS^{1}(\Qivp,\Ab \otimes \omega^{i}) = 
\ker \bigl( H^{1}(\Qivp,\Ab \otimes \omega^{i}) \to 
H^{1}(\Ivp,\Ab'' \otimes \omega^{i}) \bigr).$
\end{enumerate}
(Note that we are not allowing any variation in the choice of condition
at $v_{p}$).  We define 
$\HsS^{1}(\Qiv,\Ab \otimes \omega^{i})$ as the cokernel of 
(\ref{eq:finitedef}).
The {\it $\SS$-Selmer group} of $\Ab \otimes \omega^{i}$ is
\begin{align*}
\Sel_{\SS}(\Qi,\Ab \otimes \omega^{i}) &= \ker 
\bigl( H^{1}(\Qi,\Ab \otimes \omega^{i}) \to 
\prod_{v} \HsS^{1}(\Qiv,\Ab \otimes \omega^{i}) \bigr) \\
&= \ker \bigl( H^{1}(\QS/\Qi,\Ab \otimes \omega^{i}) \to \prod_{v \in \Sigma}
\HsS^{1}(\Qiv,\Ab \otimes \omega^{i}) \bigr)
\end{align*}
for any finite set of places $\Sigma$ containing $p$, all
archimedean places, all places at which $\Ab$ is ramified,
and all places for which $\HfS^{1}(\Qiv,\Ab \otimes \omega^{i})$ 
does not vanish.
\end{definition}

We will be especially interested in two kinds of finite/singular
structures.  First,
the {\it minimal structure} $\SS_{\min}$ on $\Ab \otimes \omega^{i}$ is given
by
$$\HfSm^{1}(\Qiv,\Ab \otimes \omega^{i}) = 0$$
for $v \neq v_{p}$.

Next let $f$ be a newform as in (4) above.
We define the {\it induced structure} $\SS(f,i)$ on 
$\Ab \otimes \omega^{i}$ by
setting
$$\HfSf^{1}(\Qiv,\Ab \otimes \omega^{i}) = 
\ker \bigl(
H^{1}(\Qiv,\Ab \otimes \omega^{i}) \to \Hs^{1}(\Qiv,\Afi) \bigr)$$
for all places $v$ of $\Qi$.
Note that by the definition of $\Hs^{1}(\Qivp,\Afi)$ we have
\begin{multline*}
\HfSf^{1}(\Qivp,\Ab \otimes \omega^{i}) = \ker
H^{1}(\Qivp,\Ab \otimes \omega^{i}) \to
H^{1}(\Ivp,\Ab'' \otimes \omega^{i}) \\
\to
H^{1}(\Ivp,\Afi''); \end{multline*}
in fact, since $H^{0}(\Ivp,\Afi'')$ is either $A_{f}''$ (for $i=0$)
or else zero (for $i \neq 0$), the latter map is injective, so that
$$\HfSf^{1}(\Qivp,\Ab \otimes \omega^{i}) = \ker
H^{1}(\Qivp,\Ab \otimes \omega^{i}) \to
H^{1}(\Ivp,\Ab'' \otimes \omega^{i})$$
as required.

\begin{lemma} ~\label{lemma:sared}
For $v \nmid p$ we have
$$\HfSf^{1}(\Qiv,\Ab \otimes \omega^{i}) = \im \bigl(
\Afi^{\Gv}/\pi \inj H^{1}(\Qiv,\Ab \otimes \omega^{i}) \bigr)$$
In particular, $\HfSf^{1}(\Qiv,\Ab \otimes \omega^{i})=0$
if $\Ab \otimes \omega^{i}$ is unramified at $v \neq v_{p}$ so that
$\SS(f,i)$ is a finite/singular structure.
The natural map
$$\Sel_{\SS(f,i)}(\Qi,\Ab \otimes \omega^{i}) \to \Sel(\Qi,\Afi)[\pi]$$
is an isomorphism.
\end{lemma}
\begin{proof}
The lemma follows from various exact sequences in cohomology
coming from the exact sequence
$$0 \to \Ab \otimes \omega^{i} \to \Afi \overset{\pi}{\longrightarrow} 
\Afi \to 0$$
together with Lemma~\ref{lemma:locdiv}.
We leave the details to the reader.
\end{proof}

\begin{prop} \label{prop:levellowering}
Let $\Ab$ be as above.  Then there exists a newform $f$ as in (4) above
such that the finite/singular structures $\SS_{\min}$ and
$\SS(\Afi)$ coincide under the isomorphism
$\Ab \otimes \omega^{i} \cong \Afi[\pi']$ for any $i$.
\end{prop}
\begin{proof}
Let $N$ denote the tame conductor of the Galois representation
$\Ab$.  As $f$ is $p$-distinguished, by 
\cite[Thm.\ 6.4]{diamond1}
there exists a finite totally ramified
extension $K'$ of $K$ and a (not necessarily unique) $p$-stabilized
newform $f \in K'[[q]]$ of tame level $N$ and weight $2$ satisfying
the condition of (4).
(The final condition in (4) is in fact already automatic from the
more standard level lowering result \cite[Thm.\ 1.1]{diamond1}
unless
$\Ab$ is an unramified $G_{p}$-module, in which case the
Selmer case of \cite[Thm.\ 6.4]{diamond1} ensures the existence of such
an $f$.)

Since the tame conductor of $f$ equals the conductor of $\Ab$,
by Lemma~\ref{lemma:locdiv}, $\Afi^{\Gv}$ is divisible for
any place $v \neq v_{p}$.  It then follows from
Lemma~\ref{lemma:sared} that the structure
$\SS(\Afi)$ on $\Afi[\pi'] \cong \Ab$ 
agrees with $\SS_{\min}$, as desired.
\end{proof}

The next proposition, which follows from a result of Greenberg, 
is crucial to our method.

\begin{prop} \label{prop:reduction}
For $f$ be as above
there is an exact sequence
\begin{multline*}
0 \to \Sel_{\SS(f,i)}(\Qi,\Ab \otimes \omega^{i}) \to H^{1}(\QS/\Qi,\Ab \otimes \omega^{i}) 
\to \\
\prod_{v \in \Sigma} \HsSf^{1}(\Qiv,\Ab \otimes \omega^{i}) \to 0
\end{multline*}
for any finite set of places $\Sigma$ containing $v_{p}$, all
archimedean places and all places dividing the tame level of $f$.
\end{prop}
\begin{proof}
Since $\Afi$ is odd and $\Sel(\Qi,\Afi)$ is $\Lambda_{\O}$-cotorsion, 
by \cite[Prop.\ 2.1]{GV} there is an exact sequence
$$0 \to \Sel(\Qi,\Afi) \to
H^{1}(\QS/\Qi,\Afi) \to \prod_{v \in \Sigma} \Hs^{1}(\Qiv,\Afi) \to 0$$
with $\Sigma$ as above.  As 
$\Sel(\Qi,\Afi)$ is $\O$-divisible 
the $\pi$-torsion of this sequence is an exact sequence
$$0 \to \Sel(\Qi,\Afi)[\pi] \to H^{1}(\QS/\Qi,\Afi)[\pi]
\to \prod_{v \in \Sigma} \Hs^{1}(\Qiv,\Afi)[\pi] \to 0.$$
By Lemma~\ref{lemma:sared} and the definition of the induced structure
this sequence identifies with that of the
proposition.
\end{proof}

Our main algebraic theorems are consequences of the next result.

\begin{cor} \label{cor:key}
Let $\Ab$ be as above and let $f$ be a
$p$-stabilized newform of tame level $N_{f}$ such that
$A_{f}[\pi] \cong \Ab$ as in (4).  Then the sequence
\begin{equation} \label{eq:serge}
0 \to \Sel_{\SS_{\min}}(\Qi,\Ab \otimes \omega^{i}) \to
\Sel(\Qi,\Afi)[\pi] \to
\prod_{v \mid N_{f}} \Afi^{\Gv}/\pi \to 0
\end{equation}
is exact for any $i$.
\end{cor}
\begin{proof}
The exactness
of (\ref{eq:serge}) follows from the definitions and 
Lemma~\ref{lemma:sared} except
for the surjectivity of the last map.  For this,
by Proposition~\ref{prop:levellowering} there exists a totally ramified
extension $K'/K$ and a $p$-stabilized newform $f_{0}$ over $K'$ such that
$(\Ab \otimes \omega^{i},\SS_{\min}) \cong 
(\AfOi[\pi'],\SS(\AfOi))$.
The surjectivity is then a formal consequence of the exact sequence of
Proposition~\ref{prop:reduction}
applied to both $\AfOi$ and $\Afi$ with
$\Sigma = \{ v \mid N_{f}p \}$.
\end{proof}

\subsection{Algebraic Iwasawa invariants}
\label{sec:alginv}

Let $\Ab$ be as in the previous section.  We write
$\rhob : \GQ \to \GL_{2}(k)$ for the corresponding Galois representation.
We now use Corollary~\ref{cor:key} to study
the relations between 
Selmer groups of newforms in the Hida family of $\rhob$.
We first consider the behavior of the
corresponding finite/singular structures on $\Ab$.

\begin{lemma} \label{lemma:lambdabranch}
Let $f$ be a newform on the branch $\T(\a)$ of the Hida family of $\rhob$.
Then
the finite singular structure $\SS(f,i)$ depends only on the branch
$\T(\a)$ and $i$.
\end{lemma}
\begin{proof}
Fix a place $v \neq v_{p}$ dividing the rational prime $\ell$;
we must show that $\HfSf^{1}(\Qiv,\Ab \otimes \omega^{i})$ is independent of $f$ on the branch
$\T(\a)$.  Consider the exact sequence
\begin{equation} \label{eq:decomp}
0 \to \Afi^{\Gv}/\pi \to H^{1}(\Qiv,\Ab \otimes \omega^{i}) \to
H^{1}(\Qiv,\Afi)[\pi] \to 0
\end{equation}
defining the local condition.
By \cite[Prop.\ 2.4]{GV}, $H^{1}(\Qiv,\Afi)[\pi]$ has $k$-dimension 
equal to the multiplicity of
$\omega^{1-i}$ in the residual representation of the unramified
$\Gl$-representation $(V_{f})_{\Il}$.
By the proof of Proposition \ref{prop:conductor},
the dimension of $(V_{f})_{\Il}$ depends only on the branch $\T(\a)$.
In particular, if $(V_{f})_{\Il}=0$, then
$$\HfSf^{1}(\Qiv,\Ab \otimes \omega^{i}) = H^{1}(\Qiv,\Ab \otimes \omega^{i})$$
is certainly independent of $f$, while if $(V_{f})_{\Il} = V_{f}$, then
Lemma~\ref{lemma:locdiv} shows that
$$\HfSf^{1}(\Qiv,\Ab \otimes \omega^{i}) \cong (\Afi)^{\Gv}/\pi = 0$$
is again independent of $f$.  

To prove the lemma it thus suffices to
consider the case that $(V_{f})_{\Il}$ has dimension one.  
In this case
$f$ at $\ell$ is either:
\begin{enumerate}
\item principal series associated to a ramified character $\chi_{1}$ and
an unramified character $\chi_{2}$;
\item special associated to an unramified character $\chi$.
\end{enumerate}
We claim first that which of (1) or (2) holds is constant on the branch
$\T(\a)$.  Indeed, the semi-simplification $\rho(\a)|_{\Gl}^{\sss}$ is
a sum of two characters.  As the restriction of these characters to
$\Il$ has finite image, the conductor of $\rho(\a)|_{\Gl}^{\sss}$ does
not change on reduction modulo any height one prime of $\T(\a)$.
It follows that (1) holds for all forms on the branch $\T(\a)$ if and only
if $\rho(\a)|_{\Gl}^{\sss}$ is ramified; that is, which of (1) and (2)
holds is constant on the branch $\T(\a)$.

Let $A(\a)$ denote $\Hom_{\Zp}(\T(\a)^{2},\Qp/\Zp)$ endowed with
a $\GQ$-action via $\rho(\a)$.  
We claim that the map
$$H^{1}(\Gv,A(\a)_{(i)}[\wp_{f}]) \to H^{1}(\Gv,A(\a)_{(i)})[\wp_{f}]$$
is an isomorphism; here $\wp_{f}$ is the height one prime of $\T(\a)$
corresponding to $f$ and $A(\a)_{(i)}$ denotes $A(\a) \otimes \omega^{i}$.
Indeed, since which of (1) and (2) holds is constant on $\T(\a)$,
the multiplicity $\nu$ of $\omega^{1-i}$ in the $\Gv/\Iv$-representation
$A(\a)[\wp_{f}]_{\Iv}$ is equal to that in
$A(\a)_{\Iv}$.  Computing as in \cite[p.\ 38]{GV},
it follows that the map above is dual to the natural map
$$\T(\a)^{\nu}/\wp_{f} \to (\T(\a)/\wp_{f})^{\nu}$$
which is visibly an isomorphism.

It follows that
\begin{align*}
\HfSf^{1}(\Gv,\Ab \otimes \omega^{i}) &:= \ker\left(
H^{1}(\Gv,\Ab \otimes \omega^{i}) \to H^{1}(\Gv,A(\a)_{(i)}[\wp_{f}])
\right) \\
&= \ker\left(
H^{1}(\Gv,\Ab \otimes \omega^{i}) \to H^{1}(\Gv,A(\a)_{(i)})[\wp_{f}]
\right) \\
&= \ker\left(
H^{1}(\Gv,\Ab \otimes \omega^{i}) \to H^{1}(\Gv,A(\a)_{(i)})
\right)
\end{align*}
which is certainly independent of $f$.
\end{proof}

For a branch $\T(\a)$ of the Hida family of $\rhob$ and a place
$v \neq v_{p}$, we may now define
$$\dalg_{v}(\a,\omega^{i}) := \dim_{k} H^{1}_{f,\SS(f,i)}
(\Qiv,\Ab \otimes \omega^{i}) = \dim_{k} \Afi^{\Gv}/\pi$$
for any form $f$ lying on the branch $\T(\a)$.

We say that 
$$\mu^{\alg}(\rhob,\omega^{i}) = 0$$ 
if
$\Sel_{\SS_{\min}}(\Qi,\Ab \otimes \omega^{i})$ 
is finite dimensional over $k$; if this
is the case  we define
$$\lambda^{\alg}(\rhob,\omega^{i}) = \dim_{k}\Sel_{\SS_{\min}}(\Qi,\Ab \otimes
\omega^{i}).$$

\begin{thm} \label{thm:prealgebraic}
Let $\rhob$ be as above.
Let $f$ be any newform in the Hida family of $\rhob$.
Then $\mu^{\alg}(\rhob,\omega^{i})=0$ 
if and only if $\mu^{\alg}(f,\omega^{i})=0$.
\end{thm}
\begin{proof}
By Theorem~\ref{thm:kato}, $\mu^{\alg}(f,\omega^{i})$ vanishes if and only if
$\Sel(\Qi,\Afi)[\pi]$ is finite dimensional.  By
Corollary~\ref{cor:key} this is equivalent to the finite dimensionality of
$\Sel_{\SS_{\min}}(\Qi,\Ab \otimes \omega^{i})$, as claimed.
\end{proof}

\begin{thm} \label{thm:algebraic}
Let $\rhob$ be as above.
Assume that $\mu^{\alg}(\rhob,\omega^{i})=0$.
\begin{enumerate}
\item Let $f$ be a newform of level $N_{f}$ lying on the branch $\T(\a)$ of the
Hida family of $\rhob$.  Then
$$\lalg(f,\omega^{i}) = 
\lalg(\rhob,\omega^{i}) + \sum_{v \mid N_{f}} \dalg_{v}(\a,\omega^{i}).$$
In particular, $\lalg(f,\omega^{i})$ 
depends only on the branch $\T(\a)$ of $f$;
we write $\lalg(\a,\omega^{i})$ for this value.
\item Let $\T(\a_{1})$ and $\T(\a_{2})$ be two branches of the Hida family
of $\rhob$.  Then
$$\lalg(\a_{1},\omega^{i}) - \lalg(\a_{2},\omega^{i}) = 
\sum_{v \neq v_{p}} 
\dalg_{v}(\a_{1},\omega^{i}) - \dalg_{v}(\a_{2},\omega^{i}).$$
\end{enumerate}
\end{thm}
\begin{proof}
The first statement is immediate from
Corollary~\ref{cor:key} and the definition of $\dalg_{v}(\a,\omega^{i})$.  
The second statement follows from
the first.
\end{proof}

\section{Two-variable $p$-adic $L$-functions}
\label{sec:analytic L-functions}

\subsection{Constructions of $p$-adic $L$-functions}
\label{subsec:non-primitive}

If $N$ is an integer prime to $p$, 
then we let $\T_N^*$ denote the analogue of $\T_N$ constructed with
respect to all modular forms, rather than just cusp forms (so that
$\T_N$ is the quotient of $\T_N^*$ that acts on cusp forms).
Suitable analogues of the results stated for $\T_N$ in
section~\ref{subsec:hida}
hold for $\T_N^*$.  Furthermore,
if $\frak m$ is a maximal ideal of $\T_N^*$ for which the associated
residual Galois representation is irreducible, then the natural
map of localisations $(\T_N^*)_{\frak m} \rightarrow (\T_N)_{\frak m}$
will be an isomorphism.   For this reason, we will ultimately have
no need of the ring $\T_N^*$.  However, it is convenient to introduce
it, since it is this ring, rather than $\T_N$, that naturally acts
on some of the homology groups that we will now introduce.

For any level $M$, we let $X_1(M)$ denote the closed modular
curve of level $\Gamma_1(M)$, and let $C_1(M)$ denote its
set of cusps.  We have a short exact sequence of Hecke modules
\begin{equation}
\label{eqn:cohom}
0 \rightarrow H_1(X_1(M);\Z_p) \rightarrow H_1(X_1(M),C_1(M);
\Z_p) \rightarrow \tilde{H}_0(C_1(M);\Z_p) \rightarrow 0.
\end{equation}
(This is a part of the relative homology sequence of the pair $(X_1(M),
C_1(M))$; the tilde over the $H_0$ denotes reduced homology.)
If we localise this sequence at a maximal ideal in the Hecke
algebra corresponding to an irreducible residual $G_{\Q}$-representation,
then the localisation of $\tilde{H}_0(C_1(M);\Z_p)$ will vanish,
and the two $H_1$ terms will become isomorphic.

To ease notation we will write
$$
\hom{M}{\Zp} := H_1(X_1(M);\Zp)
$$
and
$$
\homc{M}{\Zp} := H_1(X_1(M); C_1(M); \Zp).
$$

We now consider (\ref{eqn:cohom})
with $M = N p^r$ ($r \geq 0$).  Passing to ordinary parts
(an exact functor) we obtain a short exact sequence
of $\T_N^*$-modules.  If we fix a maximal ideal
$\frak m$ of $\T_N^*$ for which the Galois representation $\rhobar_{\frak m}$
is irreducible,
then the above comments show that we obtain an isomorphism
$$\left(\hom{Np^r}{\Zp}^{\ord}\right)_{\frak m} \iso
\left(\homc{Np^r}{\Zp}^{\ord}\right)_{\frak m}.$$
Passing to the projective limit in $r,$ we obtain a corresponding isomorphism
of $(\T_N^*)_{\frak m} \iso (\T_N)_{\frak m}$-modules
$$\plim{r} \left ( \, \hom{Np^r}{\Zp}^{\ord} \right)_{\frak m}
\iso
\plim{r}\left (
(\homc{Np^r}{\Zp}^{\ord}\right )_{\frak m} .$$
We denote these isomorphic modules by $M_{\frak m}$.
There is an additional piece of structure on $M_{\frak m}$ that we 
should mention: the action of complex conjugation on the
modular curves $X_1(N p^r)$ induces an action of complex conjugation
on homology which commutes with the Hecke action.  Thus we obtain
an action of complex conjugation
on $M_{\frak m}$.  We let $M_{\frak m}^{\pm}$ denote the $\pm$-eigenspace
for this action.  Since $p$ is odd,
we have $M_{\frak m} = M_{\frak m}^+ \bigoplus M_{\frak m}^-$.

\begin{prop}\label{prop:rank one} If $\rhobar_{\frak m}$
is irreducible and $p$-distinguished, then each of the 
$(\T_N)_{\frak m}$-modules
$M_{\frak m}^{\pm}$ is free of rank one.
\end{prop}

\begin{proof} 
Let $J_{1}(Np)$ denote the Jacobian of $X_{1}(Np)$.
By \cite[Thm.~2.1]{wiles} and Theorem \ref{thm:hida}, the ordinary
$\m$-torsion
$J_1(Np)(\Qbar)[\m]^{\ord}$ is free of rank $2$ over $\T_N / \m$.
On the other hand, by \cite[Thm.~3.1]{hida2},
\begin{eqnarray}
\label{eqn:control}
\left( \ilim{r} ~ H^1(Np^r;\Qp/\Zp)^{\ord} \right)^\Gamma 
\cong H^1(Np;\Qp/\Zp)^{\ord},
\end{eqnarray}
so that
$$
\Hom(M_\m,\Qp/\Zp)[\m]
\cong H^1(Np;\Qp/\Zp)^{\ord}_\m[\m]
\cong H^1(Np;\Fp)^{\ord}_\m[\m].
$$
Since
$$
J_1(Np)(\Qbar)[p] \cong H^1(Np;\Fp),
$$
we conclude that $\Hom(M_\m,\Qp/\Zp)[\m]$, and thus also
its dual $M_\m/\m M_\m$, is free of rank $2$ over $\T_N / \m$.  
As $M_\m \otimes \Qp$ is free of rank $2$ over $(\T_N)_\m \otimes \Qp$,
it follows that $M_\m$ is free of rank $2$ over $(\T_N)_\m$.
Passing to $\pm$-subspaces establishes the proposition.
\end{proof}

The advantage of the second description of $M_{\frak m},$
as a limit of relative homology groups, is that the corresponding relative 
homology classes admit a description via modular symbols.
More precisely, for fixed $r$, there is a map
$$\P^1(\Q) \rightarrow \homc{Np^r}{\Zp}$$
defined by sending the element $a \in \P^1(\Q)$ to the
homology class corresponding to the image of the path $[\infty,a]$ 
on the modular curve $X_1(N p^r)$.  This map is compatible with
varying $r$, and so projecting to ordinary parts, localising at $\frak m,$
and then passing to the limit, we obtain a map
$$\P^1(\Q) \rightarrow M_{\frak m},$$
which we denote by $a \mapsto \{\infty,a\}.$
This map allows us to define an $M_{\frak m}$-valued measure
on $\Z_p^{\times}$ in the usual way.

\begin{definition}\label{def:measure}
For any open subset $a + p^r \Z_p$ of $\Z_p^{\times},$
we define
$$\mu(a + p^r \Z_p) = U_p^{-r} \{ \infty, a/p^r\} \in M_{\frak m}.$$
\end{definition}

\begin{prop} The function $\mu$ is a measure (i.e.\ it is additive).
\end{prop}
\begin{proof}
This is standard.
\end{proof}

Recall that the completed group ring $\BbbIw$
may naturally be regarded
as the space of 
$\Z_p$-valued measures on $\Z_p^{\times}$.
We may thus regard $\mu$ as defining an element 
$L(\frak m, N) \in M_{\frak m} \cotimes_{\Z_p} \BbbIw $.
(Here the tensor product
is completed with respect to the usual (i.e.\ profinite
topology) on $\BbbIw$ and the $\frak m$-adic topology
on $M_{\frak m}$.)
We may decompose $L(\frak m, N)$ under the action of complex conjugation
to get a pair of elements
$$L^{\pm}(\frak m, N) \in 
M_{\frak m}^{\pm} \cotimes_{\Z_p} \BbbIw.$$

If $\heightone$ is a height one prime ideal of $(\T_N)_{\frak m},$
then we may reduce $L(\frak m, N)^{\pm}$ modulo $\heightone$ to obtain
elements 
$L^{\pm}(\frak m, N)(\heightone) \in 
M_{\frak m}^{\pm}/\heightone M_{\frak m}^{\pm} \otimes_{\Z_p}
\BbbIw.$
We may also reduce $L(\frak m, N)^{\pm}$ modulo $\frak m$, and so obtain
elements $\Lbar^{\,\pm}(\frak m, N) \in 
M^{\pm}_{\frak m}/\frak m M^{\pm}_{\frak m}
\otimes_{\Z/p} \BbbIw.$

We now assume that the hypotheses of 
Proposition~\ref{prop:rank one} are satisfied. That result shows
that we may choose an isomorphism
\begin{equation}\label{eqn:first choice}
(\T_N)_{\frak m} \iso M_{\frak m}^{\pm},
\end{equation}
and so regard $L^{\pm}(\frak m, N)$ as an element of
$(\T_N)_{\frak m} \cotimes_{\Z_p} \BbbIw
\iso (\T_N)_{\frak m}[[\Z_p^{\times}]].$
For any height one prime $\heightone$ in $(\T_N)_{\frak m}$
we may regard
$L^{\pm}(\frak m, N)(\heightone)$ as an element of $\O(\heightone)
\otimes_{\Z_p} \BbbIw \iso \O(\heightone)[[\Z_p^{\times}]].$
Finally, $(\T_N)_{\frak m}/\frak m \iso~k,$ and so we may
regard
$\Lbar^{\,\pm}(\frak m, N)$ as an element of
$k \otimes_{\Z_p} \BbbIw \iso k[[\Z_p^{\times}]].$

\begin{prop}\label{prop:two vs one}
If $\heightone$ is a classical height one prime ideal
in $(\T_N)_{\frak m}$,
then $$L^{\pm}(\frak m, N)(\heightone) \in \O({\heightone})[[\Zp^\times]]$$
is the usual analytic $p$-adic $L$-function attached to
the corresponding normalised eigenform $$f_{\heightone} \in
S_k(N p^{\infty}, \O(\heightone))^{\ord}[\kappa_{\heightone}]$$
(computed with respect to a canonical period).
\end{prop}
\begin{proof}
This is the standard comparison between the specialisation
of a two-variable 
$L$-function arising from a Hida family 
and a one-variable $p$-adic 
$L$-function, and will be clear once we give a one-variable construction
analogous to our two-variable construction.
Let $\T_{r,k}$ denote the ordinary part of the Hecke algebra acting
on $S_k(Np^r)$ and let
$$
M_{r,k} = H_1\bigl(Np^r,\{\text{cusps}\},L_{k-2}(\Zp)\bigr).
$$
Choose $r$ and $k$ so that the
modular form $f_\p$ corresponds to some prime ideal $\p_{r,k}$ of $\T_{r,k}$
and let $\m_{r,k}$ be the maximal ideal of $\T_{r,k}$ containing
$\p_{r,k}$.  There is an
$(M_{r,k})_{\m_{r,k}}$-valued measure defined by
$$
\mu_{r,k}(a+p^n\Zp) = U_p^{-n} \{\infty,a/p^n\};
$$
we write
$$L_p^\pm(\m_{r,k},N) \in (M_{r,k})^\pm_{\m_{r,k}}[[\Zp^\times]]$$
for its associated power series.

We claim 
that $(M_{r,k})^\pm_{\m_{r,k}}$ is a free $(\T_{r,k})_{\m_{r,k}}$-module 
of rank one.
First note that by Theorem \ref{thm:hida}
$$
\T / \w_{r,k} \T \cong \T_{r,k}
$$
where $\T = \T_N$ and
$\w_{r,k}$ is the product of all classical primes $\p \subset \Lambda$
of weight $k$ and such that $\kappa_\p$ restricted to $1+p^r \Zp$ is trivial.
Furthermore, by \cite[Theorem 1.9]{hida3},
$$
M_\m / \w_{r,k} M_\m \cong (M_{r,k})_{\m_{r,k}}.
$$
Thus, since $M_\m^\pm$ is a free $\T_\m$-module of rank one (Proposition
\ref{prop:rank one}), it follows that
$(M_{r,k})^\pm_{\m_{r,k}}$ is a free $(\T_{r,k})_{\m_{r,k}}$-module of
rank one as well.

Our fixed identification $M_\m^\pm \cong \T_\m$ now gives an 
identification
$(M_{r,k})^\pm_{\m_{r,k}} \cong (\T_{r,k})_{\m_{r,k}}$ and we can
view $L_p^\pm(\m_{r,k},N)$ as an element of
$(\T_{r,k})_{\m_{r,k}}[[\Zp^\times]]$.
The reduction of $L_p^\pm(\m_{r,k},N)$ mod $\p_{r,k}$, viewed in
$\T_{r,k}/\p_{r,k}[[\Zp^\times]] \cong \O(\p)[[\Zp^\times]]$,
is simply the $p$-adic $L$-function of $f_\p$ computed with respect 
to a canonical period.
(This last claim follows from the definition of canonical period and the fact
that this power series has the correct interpolation property.)

By construction
\begin{eqnarray*}
L^\pm(\m,N) \equiv L^\pm(\m_{r,k},N) \pmod{\w_{r,k}} 
\end{eqnarray*}
and thus
$$
L^\pm(\m,N)(\p) \equiv L^\pm(\m_{r,k},N) \pmod{\p_{r,k}}.
$$
The right hand side of the above equation is 
the $p$-adic $L$-function of $f_\p$; the proposition follows.
(Note that the usual ambiguity of the canonical period
coming from the choice of an isomorphism $(M_{r,k})^\pm_{\m_{r,k}}
\cong (\T_{r,k})_{\m_{r,k}}$ is controlled along the Hida
family by the single identification $M^\pm_\m \cong \T_\m$.)
\end{proof}

Let us remark that in general, the normalised eigenform
$f_{\heightone}$ appearing in the statement of Proposition \ref{prop:two vs one}
need not be a newform, and thus $L^{\pm}(\frak m, N)(\heightone)$ may not
be a ``primitive'' $p$-adic $L$-function.   In the following sections,
we will make a careful study of the relationship between such 
non-primitive $p$-adic $L$-functions, and the corresponding primitive
$p$-adic $L$-functions (i.e.\ the $p$-adic $L$-functions attached
to the corresponding $p$-stabilised $p$-ordinary newforms).

If we fix a tame character $\omega^i$ for some $0 \leq i \leq p-2$,
then we may project $L^{\pm}(\frak m, N)$ onto the ``$\omega^i$-part'' 
of $\BbbIw$, and obtain 
an element of $(\T_N)_{\frak m}\cotimes_{\Z_p} \Iw_{(i)},$
which we denote by $L(\frak m, N, \omega^i).$  
(Here we choose the sign $\pm$ to equal $(-1)^i,$ since it
is well-known, and easily checked, that otherwise the projection is trivial.)

\subsection{Two variable $L$-functions on branches of the Hida family}
\label{subsec:primitive}

Fix an irreducible, $p$-ordinary, $p$-distinguished, modular Galois
representation $\rhobar: G_{\Q} \rightarrow \GL_2(k)$,
equipped with a chosen $p$-stabilisation (which we will
suppress in our discussion and notation),
as in section~\ref{subsec:family}.
Our goal in this section is to define a $p$-adic L-function
varying over each component of the Hida family of $\rhobar,$
which at any classical height one prime $\heightone$
specialises to the $p$-adic $L$-function of the newform
$f_{\heightone}$.

We begin by fixing a finite set of primes $\Sigma$,
and considering the Hecke algebra $\T_{\Sigma}$ associated to $\rhob$ as in
section~\ref{subsec:family}.  Proposition~\ref{prop:reduced and full}
yields an isomorphism
$\T_{\Sigma} \iso (\T_{N(\Sigma)})_{\frak n}$ for a
certain maximal ideal $\frak n$ of $\T_{N(\Sigma)}$.
The construction of the preceding section 
defines elements  $L(\frak n, N(\Sigma), \omega^i) 
\in (\T_{N(\Sigma)})_{\frak n}\cotimes_{\Z_p} \Iw_{(i)}$ 
(well-defined up to multiplication by an element of
$(\T_{N(\Sigma)})_{\frak n}^{\times}$).  

\begin{definition} 
~
\begin{enumerate}
\item Set $\npL{\Sigma}{\rhobar}$
equal to the element of $\T_{\Sigma} \cotimes_{\Z_p} \Iw_{(i)}$
arising from $L(\frak n, N(\Sigma), \omega^i)$ via the isomorphism
$\T_{\Sigma} \iso (\T_{N(\Sigma)})_{\frak n}$.  (This is
well-defined up to a unit of $\T_{\Sigma}$.)


\item 
Set $\npL{\Sigma}{\rhobar}(\heightone)$ equal to the element
of $\O(\heightone) \otimes_{\Z_p} \Iw_{(i)}$
obtained as the reduction of $\npL{\Sigma}{\rhobar}$ modulo
$\heightone$ for any height one prime ideal $\heightone$
of $\T_\Sigma$.  (This is well-defined up to a unit in $\O(\heightone)$.)


\item 
Set $\modpL{\Sigma}$ equal to the element
of $k \cotimes_{\Z_p} \Iw_{(i)}$ 
obtained as the reduction of
$\npL{\Sigma}{\rhobar}$ modulo the maximal ideal
of $\T_{\Sigma}$.  (This is well-defined up to a unit of $k$.)


\end{enumerate}
\end{definition}

If $\heightone$ is a classical height one
prime ideal of $\T_{\Sigma}$, then
$\heightone$ corresponds to a classical
height one prime ideal $\heightone'$ of $\T_{N(\Sigma)}$
(via the isomorphism of Proposition~\ref{prop:reduced and full}),
and we write $g_{\heightone} := f_{\heightone'}$
to denote the normalised eigenform in
$S_k(N(\Sigma) p^{\infty},\O(\heightone))^{\ord}[\kappa_{\heightone}]$
corresponding to $\heightone'$ via part~(2) of Theorem~\ref{thm:hida}.
If $\frak a$ denotes the (unique) minimal prime of $\T_{\Sigma}$ contained in
$\heightone$, then
we may also form the normalised eigenform
$f_{\heightone} \in
S_k(N(\frak a) p^{\infty},\full{\O(\heightone)})^{\ord}_{\new}
[\kappa_{\heightone}]$.
Note that in general $f_{\heightone}$ and $g_{\heightone}$
are not equal.  Indeed, they are eigenforms for different Hecke algebras:
$f_{\heightone}$ is an eigenform for $\T_{N(\frak a)}$,
while $g_{\heightone}$ is an eigenform for $\T_{N(\Sigma)}$.
Proposition~\ref{prop:reduced and full} makes it clear how they are related:
$g_{\heightone}$ is the normalised oldform obtained from $f_{\heightone}$
by ``removing the Euler factor'' at each of the primes $\ell \in \Sigma$.
Proposition~\ref{prop:two vs one} shows that
$\npL{\Sigma}{\rhobar}(\heightone)$ is the $p$-adic $L$-function
attached to $g_{\heightone}$. 

We now turn to constructing
a $p$-adic $L$-function
$\pL{\rhobar}{\frak a}
\in \full{\T(\frak a)} \cotimes_{\Z_p} \Iw_{(i)}$
for each minimal prime $\frak a$ of $\T_{\Sigma}$,
whose specialisation at each classical height one prime $\heightone$
of $\T(\frak a)$ will be equal to the $p$-adic $L$-function
of $f_{\heightone}$.  
For this, we recall the natural isomorphism
$\full{\T(\frak a)} \iso \T_{N(\frak a)}^{\new}/\frak a'$
of Definition~\ref{def:components},  
which gives rise to the composite surjection
\begin{equation}\label{eqn:surjection}
\T_{N(\frak a)} \rightarrow \T_{N(\frak a)}^{\new}
\rightarrow \T_{N(\frak a)}^{\new}/\frak a' \iso \full{\T(\frak a)}.
\end{equation}
If we let $\frak m$ denote the maximal ideal of $\T_{N(\frak a)}$
obtained as the preimage of the maximal ideal of $\full{\T(\frak a)}$
under~(\ref{eqn:surjection}),
then the construction of the preceding section
yields an $L$-function
$L(\frak m, N(\frak a), \omega^i)
\in \T_{N(\frak a)} \cotimes_{\Z_p} \Iw_{(i)}$.
The surjection~(\ref{eqn:surjection}) induces a corresponding
surjection
$$\T_{N(\frak a)}\cotimes_{\Z_p} \Iw_{(i)} \rightarrow
\full{\T(\frak a)}\cotimes_{\Z_p} \Iw_{(i)}.$$

\begin{definition}
\label{def:pLa}
We let $\pL{\rhobar}{\frak a}$ denote the image
of $L(\frak m, N(\frak a), \omega^i)$ under the preceding surjection.
\end{definition}

Proposition~\ref{prop:two vs one} shows that the specialisation
of $\pL{\rhobar}{\frak a}$ at any classical height one prime ideal $\heightone$
of $\T(\frak a)$ is equal to the $p$-adic $L$-function $L_p(f_{\heightone})$
of the associated newform $f_{\heightone}$.

\subsection{Comparisons}

Our next task is to compare 
$\pL{\rhobar}{\frak a}$ and
$\npL{\Sigma}{\rhobar}\bmod \frak a$.
Each of these two $p$-adic $L$-functions
lies in  $\full{\T(\frak a)}\cotimes_{\Z_p} \Iw_{(i)}$.
The interpolation property satisfied by
classical $p$-adic $L$-functions shows that for any height one prime
$\heightone$ of $\T(\frak a)$, the $L$-functions
$L_p(f_{\heightone})$ and $L_p(g_{\heightone})$ agree up to multiplication
by a certain product of Euler factors (reflecting the Euler factors removed
from $f_{\heightone}$ to obtain $g_{\heightone}$) and an element
of $\O(\heightone)[1/p]^{\times}$ (reflecting the ratio of the
canonical period of $f_{\heightone}$ and the canonical period of
$g_{\heightone}$).  We will show that in fact the ratio of these canonical
periods lies
in $(\full{\O(\heightone)})^{\times}$, and hence that
$L_p(f_{\heightone})$ and $L_p(g_{\heightone})$ agree up to multiplication
by a product of Euler factors, and the inevitable ambiguity of a unit
of $\full{\O(\heightone)}$.
Furthermore, we will show that
this occurs uniformly in Hida families;
in other words, that $\npL{\Sigma}{\rhobar} \bmod \frak a$
and $\pL{\rhobar}{\frak a}$
agree up to multiplication by a product of Euler factors,
and a unit in $\full{\T(\frak a)}$.

We begin by defining the relevant Euler factors.
Recall the Euler factors
$E_{\ell}(\frak a,X)$ in $\full{\T(\frak a)}[X]$
defined in section~\ref{subsec:euler}.
Since $\ell$ is prime to $p$, we may regard $\ell$ as an element
of $\Z_p^{\times}$, and so obtain a corresponding
unit element $\langle \ell \rangle \in \BbbIw$
We may substitute $\langle \ell \rangle^{-1}$ in place
of $X$ in the Euler factor $E_{\ell}(\frak a,X),$
and so obtain an element
$$E_{\ell}(\frak a,\langle \ell \rangle^{-1}) \in
\full{\T(\frak a)}[[\Z_p^{\times}]]
\iso \full{\T(\frak a)} \cotimes_{\Z_p} \BbbIw.$$ 
Recall that for $0\leq i\leq p-2$, we write $\langle \ell \rangle_i$
to denote the projection of $\langle \ell \rangle \in \BbbIw$
under the surjection $\BbbIw \rightarrow \Iw_{(i)}$.  
We write $E_{\ell}(\frak a,\langle \ell \rangle_{i}^{-1})$
to denote the corresponding element of
$\full{\T(\frak a)}\cotimes_{\Z_p} \Iw_{(i)}.$

\begin{definition} If $\frak a$ is a minimal prime of $\T_{\Sigma}$,
then we write $$E_{\Sigma}(\frak a) := \prod_{\ell \in \Sigma}
E_{\ell}(\frak a,\langle \ell \rangle^{-1}) \in 
\full{T(\frak a)}\cotimes_{\Z_p} \BbbIw,$$
and 
$$E_{\Sigma}(\frak a,\omega^i) := \prod_{\ell \in \Sigma}
E_{\ell}(\frak a,\langle \ell \rangle_{i}^{-1}) \in 
\full{T(\frak a)}\cotimes_{\Z_p} \Iw_{(i)}.$$
\end{definition}

Since $\npL{\Sigma}{\rhobar}$ and $\pL{\rhobar}{\frak a}$ are
constructed using modular symbols of levels $N(\Sigma)$
and $N({\frak a})$ respectively, in order to compare them,
we will need to be able to compare the corresponding Hecke algebras.
Recall that by construction $\T_{\Sigma} = (\T_{N(\Sigma)}')_{\frak m'}$
for a certain maximal ideal $\frak m'$ of $\T_{N(\Sigma)}'$, and
that $\frak n$ is a certain maximal ideal of
$\T_{N(\Sigma)}$ lying over $\frak m'$ with the property that
the map
$(\T_{N(\Sigma)}')_{\frak m'} \rightarrow (\T_{N(\Sigma)})_{\frak n}$
is an isomorphism.  We let $\frak m$ denote the preimage in
$\T_{N(\frak a)}$, under the surjection~(\ref{eqn:surjection}),
of the maximal ideal of $\full{\T(\frak a)}$.
Altogether, we have the following diagram of
maps between the various Hecke algebras:

\begin{equation}\label{eqn:maps}
\xymatrix{
(\T_{N(\Sigma)}')_{\frak m'} \ar[d]\ar[r]^{\iso} & (\T_{N(\Sigma)})_{\frak n}
\\
(\T_{N(\frak a)})_{\frak m} \ar[r] & \full{\T(\frak a)}.}
\end{equation}

By inverting the upper horizontal isomorphism, we obtain
a map $(\T_{N(\Sigma)})_{\frak n} \rightarrow \full{\T(\frak a)}$.
This is a homomorphism of $\T_{N(\Sigma)}'$-algebras, when we equip
the source and target with the $\T_{N(\Sigma)}'$-algebra structure
provided by the diagram~(\ref{eqn:maps}).

We write $M(N(\Sigma))_{\frak n}$ and
$M(N(\frak a))_{\frak m}$ to
denote the modules of $p$-adic modular symbols constructed in the preceding
section, for the indicated choice of tame level, and localised
at the indicated maximal ideal. 
The following result provides
the key to comparing the $L$-functions $\npL{\Sigma}{\rhobar}\bmod \frak a$
and $\pL{\rhobar}{\frak a}$.  We postpone its proof to the
section~\ref{subsec:ihara}.

\begin{thm}\label{thm:ihara}
There is an isomorphism of $\full{\T(\frak a)}$-modules
$$\full{\T(\frak a)}\otimes_{(\T_{N(\Sigma)})_{\frak n}} M(N(\Sigma))_{\frak n}
\iso \full{\T(\frak a)}\otimes_{(\T_{N(\frak a)})_{\frak m}}
M(N(\frak a))_{\frak m},$$
compatible with the action of complex conjugation,
and having the property that the induced isomorphism
$$\full{\T(\frak a)}\otimes_{(\T_{N(\Sigma)})_{\frak n}} M(N(\Sigma))_{\frak n}
\cotimes_{\Z_p} \BbbIw
\iso
\full{\T(\frak a)}\otimes_{(\T_{N(\frak a)})_{\frak m}}
M(N(\frak a))_{\frak m}
\cotimes_{\Z_p} \BbbIw$$
maps the element $1\otimes L(\frak n, N(\Sigma)) $
of the source to the product
$\left( 1 \otimes L(\frak m, N(\frak a))  E_{\Sigma}(\frak a) \right)$
in the target.
\end{thm}

\begin{cor}\label{cor:L comparison}
There is a unit $u \in \full{\T(\frak a)}$ such that 
$$
\npL{\Sigma}{\rhobar} 
\equiv u \cdot \pL{\rhobar}{\frak a} E_{\Sigma}(\frak a, \omega^i) \pmod{\a}.
$$
\end{cor}
\begin{proof}
This follows immediately from the preceding theorem.
\end{proof}

\subsection{Iwasawa invariants}

If we choose an isomorphism
\begin{equation}\label{eqn:second choice} \Iw \cong \Z_p[[T]],
\end{equation}
and hence an isomorphism
$(\T_N)_{\frak m} \cotimes_{\Z_p} \Iw_{(i)} \iso (\T_N)_{\frak m}[[T]],$
then we may regard $L(\frak m, N, \omega^i)$ as an element
of $\T_{\frak m}[[T]]$.  We are of course interested in the Iwasawa invariants
of such power series.  We begin this section by reviewing the definitions and
basic properties of such invariants.

\begin{definition}\label{def:content}
If $R$ is a ring and $f(T) \in R[[T]]$ is a one-variable
power series with coefficients in $R$, 
then we define the {\em content} of $f(T)$ to be
the ideal $\content{f(T)} \subset R$ generated
by the coefficients of $f(T)$.
\end{definition}

The proof of the next lemma is straightforward.

\begin{lemma}\label{lem:content}
Let $f(T) \in R[[T]]$.

\begin{enumerate}
\item If $\alpha: R[[T]] \rightarrow R[[T]]$ is any automorphism
of $R[[T]]$, then $\content{f(T)} = \content{\alpha(f(T))}.$

\item If $u(T) \in R[[T]]^{\times}$,
then $\content{u(T) f(T)} = \content{f(T)}$.

\item If $f(T) = g(T) h(T)$ is an equation in $A[[T]]$,
then if any two of the three power series $f(T)$, $g(T)$,
and $h(T)$ have unit content, so does the third.

\item If $\phi: R \rightarrow S$ is a morphism of rings,
and if $\bar{f}(T) \in S[[T]]$ denotes the image of $f(T)$
under the induced map $R[[T]] \rightarrow S[[T]]$,
then $\content{\bar{f}(T)} = \phi(\content{f(T)}).$
\end{enumerate}
\end{lemma}

In particular, regarding $L(\frak m, N, \omega^i)$ as a power series,
we may consider its content
$\content{L(\frak m, N, \omega^i)}$, an ideal of $(\T_N)_{\frak m}$.
Parts~(1) and~(2) of the preceding lemma show that this
ideal is in fact independent of the choice of the
isomorphisms~(\ref{eqn:first choice}) and~(\ref{eqn:second choice}).

For each height one prime $\heightone$, we may likewise
form an element $L(\frak m, N, \omega^i)(\heightone) \in \O(\heightone)[[T]]$
and the corresponding ideal
$\content{L(\frak m, N, \omega^i)(\heightone)}$ in $\O(\heightone)$.  Then
part~(4) of Lemma \ref{lem:content}
shows that this ideal
can equally well be regarded as the image of $\content{L(\frak m, N, \omega^i)}$
under the surjection $(\T_N)_{\frak m} \rightarrow 
\O(\heightone).$

Finally, we may construct an element $\Lbar(\frak m, N, \omega^i)
\in k[[T]],$ and consider the corresponding
ideal $\content{\Lbar(\frak m, N, \omega^i)}$ of $k$
generated by its coefficients.  Again, we may also regard
$\content{\Lbar(\frak m, N, \omega^i)}$ as the image of
$\content{L(\frak m, N, \omega^i)}$ under the surjection
$(\T_N)_{\frak m} \rightarrow (\T_N)_{\frak m}/\frak m \iso k.$
Note that since $k$ is field, the ideal $\content{\Lbar(\frak m, N, \omega^i)}$
is either zero, or all of $k$.

The next result is an immediate
consequence of the fact that $(\T_N)_{\frak m}$ is a local ring.

\begin{prop}\label{prop:mu vanishing}
Fix $i$, $0 \leq i \leq p-2$.  The following are equivalent:

(1) $\content{L(\frak m, N, \omega^i)} = (\T_N)_{\frak m}.$

(2) $\content{L(\frak m, N, \omega^i)(\heightone)} = \O(\heightone)$
for some height one prime
$\heightone$ of $(\T_N)_{\frak m}.$

(3) $\content{(\frak m, N, \omega^i)(\heightone)} = \O(\heightone)$
for every height one prime
$\heightone$ of $(\T_N)_{\frak m}.$

(4) $\content{\Lbar(\frak m, N, \omega^i)}$ is non-zero (and hence equal to
$k$).
\end{prop}

We remark that if $\heightone$ is a classical height one prime ideal in
$(\T_N)_{\frak m}$,
then the length of the quotient
$\O(\heightone)/\content{L(\frak m, N, \omega^i)(\heightone)}$ 
is related to the usual $\mu$-invariant of the $p$-adic $L$-function
$L(\frak m, N, \omega^i)(\heightone) \in \O(\heightone)[[T]].$
In particular, this $\mu$-invariant vanishes if and only
if $\content{L(\frak m, N, \omega^i)(\heightone)}$ is the unit ideal
of $\O(\heightone)$.

%

Analogously, we define 
the ideals 
$$
\content{\npL{\Sigma}{\rhobar}},
\content{\npL{\Sigma}{\rhobar}(\heightone)},
\text{~and~} \content{\modpL{\Sigma}}
$$
of $\T_{\Sigma},$ $\O(\heightone)$, and $k$
respectively as the content of the elements
$$
\npL{\Sigma}{\rhobar}, \npL{\Sigma}{\rhobar}(\heightone),
\text{~and~} \modpL{\Sigma}.
$$  
All of the above results apply equally 
well to these ideals since they are constructed out of
$\content{L(\frak m, N, \omega^i)}$ for some choice of $\frak m$ and $N$.

We now compute the
content of the reciprocal Euler factor $E_{\Sigma}(\frak a,\omega^i)$.

\begin{lemma}\label{lem:euler content}
The element $E_{\Sigma}(\frak a,\omega^i)$ of
$\full{T(\frak a)}\cotimes_{\Z_p} \Iw_{(i)}$
has unit content.
\end{lemma}
\begin{proof}
Part~(3) of Lemma~\ref{lem:content}
shows that it suffices to prove that
each of the reciprocal Euler factors
$E_{\ell}(\frak a,\langle \ell \rangle_{i}^{-1})$
has unit content.
If $\gamma$ is a topological generator
of $\Gamma$, then we may write $\ell^{-1}\omega(\ell)
= \gamma^{u p^n},$
for some $u \in \Z_p^{\times}$ and some integer $n \geq 0$.
If we choose our isomorphism~(\ref{eqn:second choice})
so that $\gamma \mapsto 1 + T$,
then we see that
$$E_{\ell}(\frak a,\langle \ell \rangle_i)^{-1}
= \begin{cases} \text{$(1 - (T_{\ell} \bmod \frak a') \omega^{-i}(\ell) (1 + T)^{u p^n}
+ \langle \ell \rangle \omega^{-2i}(\ell) \ell^{-1} (1+T)^{2u p^n})$} \\
\quad \text{ if $\ell$ is prime to $N(\frak a)$} \\
\text{$(1 - (T_{\ell} \bmod \frak a') \omega^{-i}(\ell) (1 + T)^{u p^n})$
otherwise}.
\end{cases}$$
In the second case
either the constant term is a unit (if $T_{\ell}$ is not a unit)
or else the coefficient of $T^{p^n}$ is a unit (if $T_{\ell}$ is a unit).
In the first case, we may compute the content after making the
substitution $(1+T)^{u} \mapsto 1+T$, in which case one sees immediately
that the coefficient of $T^{2p^{n}}$ is a unit.
\end{proof}

The following theorem establishes that if the $\mu$-invariant vanishes
for one form in a Hida family then it vanishes for every form in that family.

\begin{thm}\label{thm:mu vanishing}
The following are equivalent:
\begin{enumerate}
\item There is one ordinary newform $f$ in the Hida family of $\rhobar$
for which the $p$-adic $L$-function
$L_p(f,\omega^i)$ has vanishing $\mu$-invariant.

\item For every ordinary newform $f$ in the Hida family of $\rhobar$,
the $p$-adic $L$-function $L_p(f,\omega^i)$ has vanishing $\mu$-invariant.

\item For one irreducible component $\T(\frak a)$ of the Hida family
of $\rhobar$, the $p$-adic $L$-function
$\pL{\rhobar}{\frak a}$ has unit content.

\item For every irreducible component $\T(\frak a)$ of the Hida family of
$\rhobar$, the $p$-adic $L$-function $\pL{\rhobar}{\frak a}$ has unit content.
\end{enumerate}
\end{thm}
\begin{proof}
Let $\T(\frak a)$ be an irreducible component
of the Hida family of $\rhobar$.  Since $\full{\T(\frak a)}$
is a local ring, we deduce from Corollary~\ref{cor:L comparison}
that the following are equivalent:
$\pL{\rhobar}{\frak a}$ has unit content;
the $p$-adic $L$-function
$L_p(f_{\heightone},\omega^i)$ has unit content (i.e.\ trivial
$\mu$-invariant) for one classical height one prime $\heightone$
of $\T(\frak a)$;
the $p$-adic $L$-function
$L_p(f_{\heightone},\omega^i)$ has unit content 
for every classical height one prime $\heightone$
of $\T(\frak a)$.
(Compare Proposition~\ref{prop:mu vanishing}.)
To complete the proof of the proposition, it suffices to show
that if $\T(\frak a_1)$ and $\T(\frak a_2)$ are two irreducible
components of the Hida family of $\rhobar$,
then $\pL{\rhobar}{\frak a_1}$ has unit content if and only if
the same holds true for $\pL{\rhobar}{\frak a_2}.$
We may choose $\Sigma$ such that each of $\T(\frak a_1)$
and $\T(\frak a_2)$
is an irreducible component of $\T_{\Sigma}$.
Lemmas~\ref{lem:content}~(3)
and~\ref{lem:euler content} then show that each of
$\pL{\rhobar}{\frak a_1}$ and $\pL{\rhobar}{\frak a_2}$ have unit content
if and only if the same is true of $\npL{\Sigma}{\rhobar}$.
\end{proof}

If the equivalent conditions of the preceding theorem hold,
then we write $$\mu^{\an}(\rhobar,\omega^i) = 0.$$  (As usual, 
our notation suppresses the choice of $p$-stabilisation of $\rhobar$.)
In the case when the $\mu$-invariant vanishes we can further study 
the $\lambda$-invariant of these power series.
(Note that for a local ring  which is not a discrete valuation ring
one can only define the $\lambda$-invariant for power series
of unit content.)

\begin{definition}\label{def:lambda}
If $A$ is a local ring,
and $f(T) \in A[[T]]$ is a power series having unit
content,
then we define the $\lambda$-invariant $\lambda(f(T))$
to be the smallest degree in which $f(T)$ has a unit coefficient.
\end{definition}

We remark that if $\phi: A \rightarrow B$ is a local morphism
of complete local rings, and if $f(T)$ is an element of $A[[T]]$
having unit content, then $\lambda(f(T)) = \lambda(\phi(f(T))).$

For $\T(\frak a)$ an irreducible component of the Hida family of
$\rhobar$ with $\mu^{\an}(\rhobar,\omega^i) = 0$, by Theorem \ref{thm:mu vanishing},
we have that $\pL{\rhobar}{\a}$ has unit content.
We can therefore define the analytic $\lambda$-invariant of a branch by
$$
\lan(\rhobar, \frak a, \omega^i) = \lambda(\pL{\rhobar}{\frak a}).
$$
Our main results on analytic $\lambda$-invariants in Hida families are
as follows.

\begin{thm}
\label{thm:lambda}
Assume that $\mu^{\an}(\rhobar,\omega^i) = 0$.
\begin{enumerate}
\item  For any given irreducible component
$\T(\frak a)$ of the Hida family of $\rhobar$,
the $\lambda$-invariant of $L_p(f_{\heightone},\omega^i)$ takes on the constant
value of $\lan(\a,\omega^i)$
as $\heightone$ varies over all
classical height one primes of $\T(\frak a)$.

\item For any two irreducible components $\T(\frak a_1)$, $\T(\frak a_2)$ of
the Hida family of $\rhobar$, we have that
$$
\lan(\a_1,\omega^i) - \lan(\a_2,\omega^i) = \sum_{\ell \neq p}
( e_\ell(\a_2,\omega^i) - e_\ell(\a_1,\omega^i) )
$$
where $e_\ell(\a,\omega^i) = \lambda(E_\ell(\a,\langle \ell \rangle_i))$.
\end{enumerate}
\end{thm}

\begin{proof}
The first part follows from the remarks following
Definition~\ref{def:pLa} and Definition~\ref{def:lambda}.
For the second part, choose $\Sigma$ large enough so that both $\T(\a_1)$ and $\T(\a_2)$
are irreducible components of $\T_\Sigma$.  Then, by Corollary \ref{cor:L comparison}, 
we have that
$$
\lambda(\pL{\Sigma}{\rhobar}) = \lan(\a_j,\omega_i) + \sum_{\ell \in \Sigma} e_\ell(\a_j,\omega^i)
$$
for $j=1,2$.  Our formula then follows since $e_\ell(\a_1,\omega^i) = e_\ell(\a_2,\omega^i)$ for 
$\ell \notin \Sigma$.
\end{proof}

We remark that in section \ref{sec:examples} we will see that the formulas
of the preceding theorem compare well with the formulas of Theorem \ref{thm:algebraic}.

\subsection{Proof of Theorem \ref{thm:ihara}}
\label{subsec:ihara}

In this section, we present the proof of Theorem~\ref{thm:ihara}.
We will utilise the notation introduced prior to the statement
of the theorem, 
and begin the proof by introducing some additional notation.

If $M\geq 1$ is an integer,
$d$ is a divisor of $M$, and $d'$ is a divisor of $d$,
then we let $B_{d,d'}: X_1(M) \rightarrow X_1(M/d)$
denote the map induced by the map $\tau \mapsto d' \tau$ on the upper 
half-plane.  (It will be convenient in the following discussion
to be able to omit the level
$M$ from the notation, and hence we do so.)
This map induces a corresponding map
$$
(B_{d,d'})_* : \homc{M}{\Zp} \rightarrow \homc{M/d}{\Zp}
$$ 
on relative homology groups.

If $\ell$ is any prime distinct from $p$, then we let $e_{\ell}$ denote
the largest power of $\ell$ dividing $N(\Sigma)/N(\frak a)$.
We have the inequality $0 \leq e_{\ell} \leq 2.$   Also,
$e_{\ell} = 0$ unless $\ell \in \Sigma$, while $e_{\ell} = 2$
only if $\ell \in \Sigma$ and $\ell$ is prime to $N(\frak a)$.
For each such prime $\ell$, we write
$$\epsilon(\ell) := \begin{cases} 1 &\text{ if } e_{\ell} = 0 \\
((B_{\ell, 1})_* - \ell^{-1} T_{\ell} \, (B_{\ell,\ell})_*) &\text{ if } e_{\ell} = 1\\
((B_{\ell^2, 1})_* -  \ell^{-1} T_{\ell} \, (B_{\ell^2,\ell})_* + \ell^{-3} \langle \ell \rangle
(B_{\ell^2,\ell^2})_*) &\text{ if } e_{\ell} = 2. \end{cases}$$
Now write $\Sigma = \{\ell_1,\ldots,\ell_n\}$,
and for any $r \geq 1$ define 
$$
\epsilon_r :
\homc{N(\Sigma) p^r}{\Z_p}^{\ord}
\rightarrow \homc{N(\frak a) p^r}{\Z_p}^{\ord}$$
by
$\epsilon_r = \epsilon(\ell_n) \circ \cdots \circ \epsilon(\ell_1)$.

Let us explain this formula.  For any $i = 1 ,\ldots,n$,
write $N_i = N(\Sigma)/\ell_1^{e_{\ell_1}}\cdots \ell_i^{e_{\ell_i}}.$
In the formula for $\epsilon_r,$
the map $\epsilon(\ell_i)$
is taken to be the map $$\homc{N_{i-1} p^r}{\Z_p}^{\ord}
\rightarrow \homc{N_i p^r}{\Z_p}^{\ord}$$
given by the stated formula for $\epsilon(\ell_i)$.
(The symbol $T_{\ell_i}$ in the formula for $\epsilon(\ell_i)$
is understood to stand for the corresponding Hecke operator
acting in level $N_i p^r$.)
It is easily verified that the map $\epsilon_r$
is in fact independent of the choice of ordering
of the elements of $\Sigma$.

For any tame level $M$ we let
$(\T^*_{M})'$ denote the $\Iw$-subalgebra of the ordinary Hecke algebra
$\T^*_{M}$ generated by the Hecke operators prime to $M$.
If we regard the source and target of $\epsilon_r$ as
$(\T^*_{N(\Sigma)})'$-modules
via the inclusion $(\T^*_{N(\Sigma)})' \subset \T^*_{N(\Sigma)}$
and the natural map $(\T^*_{N(\Sigma)})' \rightarrow
(\T^*_{N(\frak a)})'\subset \T^*_{N(\frak a)}$,
then $\epsilon_r$ is immediately seen to be $(\T^*_{N(\Sigma)})'$-linear.

As $r$ varies, the sources and targets of the maps $\epsilon_r$ each form
a projective system, and the maps $\epsilon_r$ are evidently compatible
with the projection maps on source and target.  Thus,
passing to the limit in $r$, we obtain a $(\T^*_{N(\Sigma)})'$-linear map
$$
\epsilon_{\infty}:\plim{r} \, \homc{N(\Sigma) p^r}{\Z_p}^{\ord}
\rightarrow \plim{r} \, \homc{N(\frak a)p^r}{\Z_p}^{\ord}.$$
Recall that the source and target of this map are denoted by 
$M(N(\Sigma))$ and $M(N(\frak a))$ respectively.

We may regard each of the maximal
ideals $\frak m'$, $\frak n$, and $\frak m$ equally well as maximal ideals
of $(\T_{N(\Sigma)}^*)',$ $\T_{N(\Sigma)}^*$, and $\T_{N(\frak a)}^*.$ 
If we localise $\epsilon_{\infty}$ with respect to $\frak m'$,
we obtain a map
\begin{equation}\label{eqn:intermediate}
M(N(\Sigma))_{\frak m'} \rightarrow M(N(\frak a))_{\frak m'}.
\end{equation}
Now $\frak n$ and $\frak m$ each pull back to $\frak m'$ under the natural
maps $\T_{N(\Sigma)}' \rightarrow \T_{N(\Sigma)}$ and
$\T_{N(\Sigma)}' \rightarrow \T_{N(\frak a)}$, and
so the localisations 
$(\T_{N(\Sigma)})_{\frak n}$ and $(\T_{N(\frak a)})_{\frak m}$
are local factors of the complete semi-local rings
$(\T_{N(\Sigma)}')_{\frak m'} \otimes_{\T_{N(\Sigma)}'} \T_{N(\Sigma)}$
and
$(\T_{N(\Sigma)}')_{\frak m'} \otimes_{\T_{N(\Sigma)}'} \T_{N(\frak a)}$.
Thus the localisations $M(N(\Sigma))_{\frak n}$ and $M(N_{\frak a})_{\frak m}$
appear naturally as direct factors of $M(N(\Sigma))_{\frak m'}$ and
$M(N_{\frak a})_{\frak m'}$ respectively,
and so the map~(\ref{eqn:intermediate}) induces a map
$$M(N(\Sigma))_{\frak n} \rightarrow M(N(\frak a))_{\frak m}.$$
Tensoring the source of this map with $\full{\T(\frak a)}$ over
$(\T_{N(\Sigma)})_{\frak n}$, and the target
with $\full{\T(\frak a)}$ over $(\T_{N(\frak a)})_{\frak m}$,
we obtain a $\full{\T(\frak a)}$-linear map
\begin{equation}
\label{eqn:ihara}
\full{\T(\frak a)}\otimes_{(\T_{N(\Sigma)})_{\frak n}} M(N(\Sigma))_{\frak n}
\rightarrow \full{\T(\frak a)}\otimes_{(\T_{N(\frak a)})_{\frak m}}
M(N(\frak a))_{\frak m}.
\end{equation}
We claim that this map satisfies the requirements of Theorem~\ref{thm:ihara}.

We first observe that this map satisfies the claimed
property with regard to the $L$-functions.  This follows
from the an explicit calculation of the effect of the maps
$\epsilon_{\ell}$ on modular symbols.
In fact, one easily shows that for any character $\chi$ of conductor
$p^n$ and $\displaystyle \alpha = 
\sum_{a \in (\Z/p^n)^{\times}} \chi(a) \left\{a/p^n,\infty\right\}$,
we have that
$$
\epsilon_r\left( \alpha \right) = 
\prod_{i \text{ s.t.\ } e_{\ell_i} = 1}(1 - \chi^{-1}(\ell) \ell^{-1} T_{\ell}  )
\prod_{i \text{ s.t.\ } e_{\ell_i} = 2} (1 - \chi^{-1}(\ell) \ell^{-1} T_{\ell}  
+ \chi^{-2}(\ell) \ell^{-3} \langle \ell \rangle   ) 
\cdot \alpha 
$$
(In the left hand side of this equation, the modular symbols
$\left\{a/p^n,\infty\right\}$ are regarded as lying in
$\homc{N(\Sigma) p^r}{\Z_p}.$
On the right hand side, they are regarded as lying in
$\homc{N(\frak a) p^r}{\Z_p}.$)
Passing to the limit in $r$, 
and taking into account the fact that $\chi$ is an arbitrary
Dirichlet character of $p$-power conductor,
we conclude the the isomorphism~(\ref{eqn:ihara}) has the
required effect on $L$-functions.

We now turn to showing that~(\ref{eqn:ihara}) is an isomorphism.
Note that Proposition~\ref{prop:rank one} shows that
both source and target are free of rank two over
$\full{\T(\frak a)}$.  Thus to see that this map is an isomorphism,
it suffices to check that it induces a surjection after being reduced
modulo the maximal ideal of $\T(\frak a)^\circ$.  To do this, we first 
mod out by a classical prime of weight two and then by the full maximal ideal.

Let $\heightone$ denote the classical height one prime in $\Iw$
of weight two for which $\kappa_{\heightone}$ is trivial.
If we tensor each side of~(\ref{eqn:ihara}) by $\Iw/\heightone$
over $\Iw$, we obtain the map 
\begin{equation}\label{eqn:weight two map}
\full{\T(\frak a)}\otimes_{(\T_{N(\Sigma)})_{\frak n}}
\hom{N(\Sigma) p}{\Z_p}^{\ord}_{\frak n} 
\rightarrow \full{\T(\frak a)}\otimes_{(\T_{N(\frak a)})_{\frak m}}
\hom{N(\frak a) p}{\Z_p}^{\ord}_{\frak m}
\end{equation}
induced by localising the source and target of $\epsilon_1$ 
at $\frak n$ and $\frak m$ respectively, and then extending scalars
to $\full{\T(\frak a)}$. 
(Here, as in Proposition \ref{prop:rank one}, we are using 
\cite[Thm.~3.1]{hida2}.)

Thus, 
the reduction modulo
the maximal ideal of $\full{\T(\frak a)}$
of~(\ref{eqn:ihara}) 
coincides with the map
\begin{multline}
\label{eqn:reduced map}
\left( \T_{N(\frak a)}/\frak m \right)\otimes_{\T_{N(\Sigma)}/\frak n}
\left( \hom{N(\Sigma) p}{\Z_p}^{\ord}/\frak n
\hom{N(\Sigma) p}{\Z_p}^{\ord} \right) \\
\rightarrow 
\hom{N(\frak a) p}{\Z_p}^{\ord}/ \frak m
\hom{N(\frak a) p}{\Z_p}^{\ord}
\end{multline}
of $\T_{N(\frak a)}/\frak m$-vector spaces
induced by~(\ref{eqn:weight two map}).

Rather than showing directly that~(\ref{eqn:reduced map}) is surjective,
we will show that the corresponding dual map
\begin{equation}\label{eqn:dual}
\cohom{N(\frak a) p}{\F_p}^{\ord}[\frak m]
\rightarrow
\left( \T_{N(\frak a)}/\frak m \right)\otimes_{\T_{N(\Sigma)}/\frak n}
\left( \cohom{N(\Sigma) p}{\F_p}^{\ord}[\frak n]\right)
\end{equation}
is injective.  (In writing the dual of~(\ref{eqn:reduced map})
in this form, we have implicitly fixed an isomorphism of one
dimensional $\T_{N(\frak a)}/\frak m$-vector spaces
between
$\T_{N(\frak a)}/\frak m$ and its
$\T_{N(\Sigma)}/\frak m'$-linear dual.  We will suppress this
choice of isomorphism here and below.)
This map may be written as a composite
\begin{equation}
\label{eqn:factored}
\begin{aligned}
\cohom{N(\frak a) p}{\F_p}^{\ord}[\frak m] 
&\rightarrow
\left( \T_{N(\frak a)}/\frak m \right) \otimes_{\T_{N(\Sigma)}/\frak n}
\left( \cohom{N(\frak a) p}{\F_p}^{\ord}[\frak m'] \right) \\
&\rightarrow 
\left( \T_{N(\frak a)}/\frak m \right)\otimes_{\T_{N(\Sigma)}/\frak n}
\left( \cohom{N(\Sigma) p}{\F_p}^{\ord}[\frak m']\right)\\
&\rightarrow
\left( \T_{N(\frak a)}/\frak m \right)\otimes_{\T_{N(\Sigma)}/\frak n}
\left( \cohom{N(\Sigma) p}{\F_p}^{\ord}[\frak n]\right).
\end{aligned}
\end{equation}
To explain this, we first recall that since
the localisation $(\T_{N(\Sigma)})_{\frak n}$ is a local factor
of the tensor product
$(\T_{N(\Sigma)}')_{\frak m'} \otimes_{\T_{N(\Sigma)}'} \T_{N(\Sigma)}$
for which the natural map $(\T_{N(\Sigma)}')_{\frak m'} 
\rightarrow (\T_{N(\Sigma)})_{\frak n}$ is an isomorphism,
the residue field
$\T_{N(\Sigma)}/\frak n$ is a local factor of the Artin local ring
$\T_{N(\Sigma)}/\frak m' \T_{N(\Sigma)};$
write
$$\T_{N(\Sigma)}/\frak m' \T_{N(\Sigma)} \iso \T_{N(\Sigma)}/\frak n
\times A,$$
where $A$ denotes the product of the remaining local factors.
This decomposition induces a corresponding decomposition
\begin{equation}\label{eqn:decomp}
B \iso B[\frak n] \times A\otimes_{\T_{N(\Sigma)}/\frak m' \T_{N(\Sigma)}} B
\end{equation}
for any $\T_{N(\Sigma)}/\frak m' \T_{N(\Sigma)}$-module $B$.
The third arrow of~(\ref{eqn:factored}) is precisely projection
onto the the first of the two factors in~(\ref{eqn:decomp}),
with 
$$
B =
\left( \T_{N(\frak a)}/\frak m \right)\otimes_{\T_{N(\Sigma)}/\frak n}
\left( \cohom{N(\Sigma) p}{\F_p}^{\ord}[\frak m']\right).
$$
The second arrow is induced by dualising the reduction modulo
$\frak m'$ of the map~(\ref{eqn:intermediate}),
and the first arrow is induced from the obvious inclusion,
or if the reader prefers, is obtained by dualising the surjection
\begin{multline}
\left( \T_{N(\frak a)}/\frak m \right)\otimes_{\T_{N(\Sigma)}/\frak n}
\left( \hom{N(\frak a) p}{\Z_p}^{\ord}/\frak m'
\hom{N(\frak a) p}{\Z_p}^{\ord} \right) \\
\rightarrow 
\hom{N(\frak a) p}{\Z_p}^{\ord}/ \frak m
\hom{N(\frak a) p}{\Z_p}^{\ord}.
\end{multline}

\begin{lemma}\label{lem:vanishing}
The map~(\ref{eqn:dual}) is injective if and only
if the composite of the first two arrows of~(\ref{eqn:factored})
is injective.
\end{lemma}
\begin{proof} 
The only if statement is clear.
In order to prove the other statement, we first
observe that the second arrow of~(\ref{eqn:factored}) is given by
the cohomological version of the map $\epsilon_1$.  
More precisely, if for each $\ell \in \Sigma$ we define
$$\epsilon_{\ell}^* := \begin{cases} 1 &\text{ if } e_{\ell} = 0 \\
B_{\ell, 1}^* -  B_{\ell,\ell}^* \, \ell^{-1} T_{\ell} &\text{ if } e_{\ell} = 1\\
B_{\ell^2, 1}^* - B_{\ell^2,\ell}^* \, \ell^{-1} T_{\ell} +
B_{\ell^2,\ell^2}^* \ell^{-3} \langle \ell \rangle  &\text{ if } e_{\ell} = 2
\end{cases}$$
(where $B_{d,d'}^*$ denotes the map on cohomology induced by
the degeneracy map $B_{d,d'}$), 
and define
$$
\epsilon^* = \epsilon_{\ell_1}^* \circ \cdots \circ \epsilon_{\ell_n}^*:
\cohom{N(\frak a) p}{\F_p}^{\ord}
\rightarrow \cohom{N(\Sigma) p}{\F_p}^{\ord}.
$$
Then the second arrow of~(\ref{eqn:factored}) is
obtained from the map~$\epsilon^*$ by passing to the kernel
of $\frak m'$ in the source and target, and then extending scalars
from $\T_{N(\Sigma)}/\frak n$ to $\T_{N(\frak a)}/\frak m$.

One then checks that any element in the image of
$\epsilon^*$ is annihilated by the Hecke operators
$T_{\ell}$, for those $\ell \in \Sigma$. 
 On the other hand,
Proposition~\ref{prop:reduced and full} shows that
the maximal ideal $\frak n$ is uniquely characterised by
the property of containing these operators.  
Thus the image of the second arrow of~(\ref{eqn:factored}) 
lies in the local factor
$$\left( \T_{N(\frak a)}/\frak m \right)\otimes_{\T_{N(\Sigma)}/\frak n}
\left( \hom{N(\Sigma) p}{\F_p}^{\ord}[\frak n]\right)$$
of
$$\left( \T_{N(\frak a)}/\frak m \right)\otimes_{\T_{N(\Sigma)}/\frak n}
\left( \hom{N(\Sigma) p}{\F_p}^{\ord}[\frak m']\right).$$
This proves the lemma.
\end{proof}

By the preceding lemma, we are reduced to proving
that the composite of the first two arrows
of~(\ref{eqn:factored}) is injective.  
It will be notationally easier to deal with each of the maps
$\epsilon_{\ell_i}$ separately,
and so we put ourselves in the
following more general situation.  We consider a natural number
$M$ prime to $p$ and a prime $\ell$ distinct from $p$.
We define
$n_{\ell} = \text{ $1$ or $2$ }$ according to whether or not
$\ell$ divides $M$,
and write $N:= \ell^{n_{\ell}} M $.
We let $\frak m$ denote a maximal ideal of $\T_M$ for which
$\rhobar_{\frak m}$ is irreducible and $p$-distinguished, 
and let
$\frak m'$ denote the pull back of $\frak m$
under the natural map $\T_N' \rightarrow \T_M$.
The map $\epsilon_{\ell}^*$ defined in the proof of Lemma~\ref{lem:vanishing}
induces a $\T_M/\frak m$-linear map
\begin{equation}\label{eqn:general}
\cohom{M p}{\F_p}^{\ord} [\frak m] \rightarrow
\left (\T_M/\frak m\right )
\otimes_{\T_N'/\frak m'} \left( \cohom{N p}{\F_p}^{\ord} [\frak m']\right).
\end{equation}

\begin{lemma}\label{lem:injectivity}
The map~(\ref{eqn:general}) is injective.
\end{lemma}

The proof of Lemma~\ref{lem:vanishing} shows
that elements in the image of~(\ref{eqn:general}) are annihilated by
$T_{\ell}$, and hence are in fact annihilated by
a maximal ideal of the full Hecke algebra $\T_N$.
Thus we may apply Lemma~\ref{lem:injectivity}
inductively to establish the injectivity of
the composite of the first two arrows of~(\ref{eqn:factored}),
and thereby complete the proof of the theorem.

\begin{proof}[Proof of Lemma~\ref{lem:injectivity}]
It will be convenient to bring the residue field $\T_M/\frak m$
inside the coefficients of cohomology.
To do this, we choose a finite field $k$ containing $\T_M/\frak m$,
and tensor the source and target of~(\ref{eqn:general})
with $k$ over $\F_p$,
to obtain a map of
$k\otimes_{\F_p} \T_M/\frak m$-modules.
The chosen inclusion of $\T_M/\frak m$ determines a projection
\begin{equation}\label{eqn:proj}
k\otimes_{\F_p} \T_M/\frak m \rightarrow k,
\end{equation}
which realises $k$ as a local factor of $k\otimes_{\F_p} \T_M/\frak m$.
Projecting onto this local factor,
we recover our original map~(\ref{eqn:general}),
but rewritten as
\begin{equation}\label{eqn:rewritten}
\cohom{M p}{k}^{\ord}[\frak m_k]
\rightarrow 
\cohom{N p}{k}^{\ord}[\frak m_k'].
\end{equation}
Here we regard the source as a $W(k) \otimes_{\Z_p} \T_{M}$-module,
and the target as a $W(k) \otimes_{\Z_p} \T_{N}'$-module.
Also, we have written $\frak m_k$ to denote the
the maximal ideal of $W(k) \otimes_{\Z_p} \T_M$
that is the kernel of the composite
$$W(k) \otimes_{\Z_p} \T_M
\rightarrow k \otimes_{\F_p} \T_M/\frak m \rightarrow
k$$ 
(where the second arrow is given by
the projection~(\ref{eqn:proj})), 
and $\frak m_k'$ to denote 
the maximal ideal of $W(k) \otimes_{\Z_p} \T_{N}'$
that is the kernel of the composite
$$W(k) \otimes_{\Z_p} \T_{N}'
\rightarrow k \otimes_{\F_p} \T_{N}'/\frak m' \rightarrow
k \otimes_{\F_p} \T_M/\frak m \rightarrow k$$
(where the second arrow is obtained by tensoring the
injection $\T_{N}'/\frak m' \rightarrow \T_{M}/\frak m$
by $k$ over $\F_p$, and the third arrow is given by
the projection~(\ref{eqn:proj})). 
We must show that the map~(\ref{eqn:rewritten}) is injective.

Our argument will rely on the results of \cite[\S 2.2]{wiles},  
which extend a well-known result of Ihara.
Recall that the standard practice for analysing the
map~(\ref{eqn:rewritten}) is to factor the map $\epsilon_{\ell}^*$,
and so to write~(\ref{eqn:rewritten}) as a composite
\begin{equation*}
\cohom{M p}{k}^{\ord}[\frak m_k]
\buildrel \alpha_{\ell} \over \longrightarrow
(\cohom{M p}{k}^{\ord}[\frak m_k])^{n_{\ell} + 1} 
\buildrel \beta_{\ell} \over \longrightarrow
\cohom{N p}{k}^{\ord}[\frak m_k'],
\end{equation*}
where 
$$\alpha_{\ell} = \begin{cases} (1, - \ell^{-1} T_{\ell}) &\text{ if } n_{\ell} = 1
\\
(1, - \ell^{-1} T_{\ell},  \ell^{-3} \langle \ell \rangle )  &\text{ if }
n_{\ell} = 2 \end{cases}$$
and
$$\beta_{\ell}  = \begin{cases}  B_{\ell, 1}^* \pi_1
+ B_{\ell,\ell}^* \pi_2 &\text{ if }
n_{\ell} = 1 \\
B_{\ell^2,1}^* \pi_1 + B_{\ell^2, \ell}^* \pi_2 + B_{\ell^2,\ell^2}^* \pi_3 
&\text{ if } n_{\ell} = 2,\end{cases}$$
with $\pi_i$ denoting projection onto the $i$th factor of the product
$$(\cohom{N p}{k}^{\ord}[\frak m_k])^{n_{\ell} + 1}.$$

The map $\alpha_{\ell}$ is manifestly injective,
and in the case $n_{\ell} = 2$, Wiles has shown
that the map $\beta_{\ell}$ is also injective.
(See in particular the discussion at the top
of \cite[p.~497]{wiles},
or the discussion of Wiles' results provided by \cite[\S \S 4.3]{diamond-ribet}
and \cite[\S\S 4.4, 4.5]{darmon-diamond-taylor}.)
Thus when $n_{\ell} = 2$, the lemma is proved.

However, if $n_{\ell} = 1$, then the situation is more complicated.
Lemma~2.5 of \cite{wiles} provides an exact sequence
\begin{equation*}
\label{eqn:1}
\cohom{(M/\ell) p}{k}^{\ord}_{\frak m_k'}
\buildrel \gamma_{\ell} \over \longrightarrow
(\cohom{M p}{k}^{\ord}_{\frak m_k'})^2
\buildrel \beta_{\ell} \over \longrightarrow
\cohom{N p}{k}^{\ord}_{\frak m_k'},
\end{equation*}
where $\gamma_{\ell} = (B_{\ell,\ell}^*, - B_{\ell,1}^*).$
(Here the subscript $\frak m_k'$ denotes that we have localised
at this maximal ideal.)

Let $x \in \cohom{M p}{k}[\frak m]$ be
an element in the kernel of~(\ref{eqn:rewritten}). 
Then (\ref{eqn:1}) shows that we may find
$y \in \cohom{(M/\ell) p}{k}$ such that
$ (x, - T_{\ell} x) = (B_{\ell,\ell}^* y, - B_{\ell,1}^*y).$
In particular, $B_{\ell,\ell}^*y = x$ is annihilated by $\frak m$,
and so is an eigenvector for the full Hecke algebra $\T_M$.
The following result shows that $B_{\ell,\ell}^*y = 0,$
and hence that $x = 0$. This completes
the proof of the lemma.
\end{proof}

\begin{lemma}
Let $D$ be a natural number prime to $p$, let $\ell$
be a prime distinct from $p$, and consider the map
$B_{\ell,\ell}^*:
\cohom{D p}{\Fpbar}^{\ord} \rightarrow 
\cohom{D \ell p}{\Fpbar}^{\ord}$.
If $y$ is a class in the domain, with the property
that $B_{\ell,\ell}^*y$ is an eigenvector for the full
Hecke ring $\T_{D \ell }$, corresponding to a maximal
ideal for which the attached Galois representation
into $\GL_2(\Fpbar)$ is irreducible and $p$-distinguished,
then $B_{\ell,\ell}^* y = 0$.
\end{lemma}
\begin{proof}
We will prove the lemma by comparing the map $B_{\ell,\ell}^*$
on cohomology classes with the corresponding map on modular forms
$\bmod p$.  For such modular forms, the analogue of the lemma follows
immediately from a consideration of $q$-expansions.
For the comparison with modular forms $\bmod p$, we follow the
discussion in the proof of \cite[Thm.~2.1]{wiles}.

Let us make some initial reductions in the situation of the lemma.
We let $\frak m$ denote the maximal ideal in $\Fpbar\otimes_{\F_p} \T_{D\ell}$
describing the action of the Hecke operators on the eigenvector
$B_{\ell,\ell}^*y$.  
Let $\frak m'$ denote the intersection of $\frak m$ with
the subring $\Fpbar\otimes_{\F_p} \T_{D \ell}'$ of
$\Fpbar \otimes_{\F_p} \T_{D \ell}$.
We may regard each of
$\cohom{D p}{\Fpbar}^{\ord}$
and
$\cohom{D \ell p}{\Fpbar}^{\ord}$
as $\Fpbar\otimes_{\F_p} \T_{D \ell}'$-modules (the former via the natural map
$\Fpbar\otimes_{\F_p} \T_{D\ell}' \rightarrow \Fpbar\otimes_{\F_p} \T_D'$),
and the map $B_{\ell,\ell}^*$
is $\Fpbar\otimes_{\F_p} \T_{D\ell}'$-linear.
If $y_{\frak m'}$ denotes the projection of $y$ onto the
localisation 
$\cohom{D p}{\Fpbar}^{\ord}_{\frak m'}$
of $\cohom{D p}{\Fpbar}^{\ord}$ at $\frak m'$,
then $B_{\ell,\ell}^* y_{\frak m'} = B_{\ell,\ell}^* y$
(since $B_{\ell,\ell}^* y$ is assumed to be annihilated by
$\frak m$, and so in particular by $\frak m'$).
Also, there is a natural isomorphism
$$\cohom{D p}{\Fpbar}^{\ord}_{\frak m'} \iso
\prod_{\frak n \supset \frak m' (\Fpbar\otimes_{\F_p}\T_D)}
\cohom{D p}{\Fpbar}^{\ord}_{\frak n},$$
where $\frak n$ ranges over all maximal ideals of $\Fpbar\otimes_{\F_p} \T_D$
containing
$\frak m' (\Fpbar \otimes_{\F_p} \T_D)$.  
Thus to show that $B_{\ell,\ell}^* y =0,$
it suffices to show that $B_{\ell,\ell}^* y_{\frak n} = 0$
for each such maximal ideal $\frak n$, where
$y_{\frak n}$ denotes the projection of $y$ onto
$\cohom{D p}{\Fpbar}^{\ord}_{\frak n}.$
Thus for the duration of the proof  we assume that
$y \in \cohom{D p}{\Fpbar}^{\ord}_{\frak n}$,
for some maximal ideal $\frak n$ of $\Fpbar\otimes_{\F_p} \T_D$
lying over $\frak m'$.   

As in the proof of \cite[Thm.~2.1]{wiles},
we consider two cases: that in which the
mod $p$ diamond operators have non-trivial image in
$(\Fpbar \otimes_{\F_p} \T_{D \ell})/\frak m$,
and that in which they have trivial image.

In order to treat the first case,
we begin by noting that the commutative diagram
\begin{equation}\label{eqn:diagram}
\xymatrix{ 
\cohom{D p}{\Fpbar}^{\ord}_{\frak n} \ar[r]^-{B_{\ell,\ell}^*} \ar[d]^{\iso} &
\cohom{D \ell p}{\Fpbar}^{\ord}_{\frak m} \ar[d]^{\iso} \\
(\Fpbar\otimes_{\F_p} J_1(D p)[p]^{\ord})_{\frak n}
\ar[r]^-{B_{\ell,\ell}^*} &
(\Fpbar\otimes_{\F_p} J_1(D \ell p)[p]^{\ord})_{\frak m}  ,}
\end{equation}
in which the lower horizontal arrow is induced by the map
of Jacobians arising by Picard functoriality applied to
the degeneracy map $B_{\ell,\ell}: X_1(D \ell p) \rightarrow X_1(D p),$
allows us to replace
curves with coefficients in $\Fpbar$ with a consideration of $p$-torsion
in the corresponding Jacobians (tensored by $\Fpbar$ over $\F_p$).
(The superscript $\ord$ on $J_1(D p)[p]$ and $J_1(D \ell p)[p \ell]$
denotes the localisation of these $p$-torsion modules at the $p$-ordinary
part of the Hecke algebra, or equivalently, the image of the
$p$-torsion modules under Hida's idempotent $e^{\ord}$.)

Let $\rhobar: G_{\Q} \rightarrow \GL_2(\Fpbar)$ denote the
residual Galois representation attached to each of $\frak m$ and
$\frak n$.
From \cite[Thm.~2.1]{wiles}
and its proof we conclude that there are isomorphisms of Galois modules
$$(\Fpbar\otimes_{\F_p} \T_D)_{\frak n}/
p (\Fpbar\otimes_{\F_p} \T_D)_{\frak n}
\otimes_{\Fpbar} \rhobar \iso 
(\Fpbar\otimes_{\F_p} J_1(D p)[p]^{\ord})_{\frak n}$$
and
$$(\Fpbar\otimes_{\F_p} \T_{D\ell})_{\frak m}/
p (\Fpbar\otimes_{\F_p} \T_{D\ell})_{\frak m} \otimes_{\Fpbar} \rhobar \iso 
(\Fpbar\otimes_{\F_p} J_1(D \ell p)[p]^{\ord})_{\frak m}.$$
Also, passing to $p$-torsion in the short exact sequence
\cite[(2.2)]{wiles} for each of the maximal ideals $\frak m$ and
$\frak n$ yields short exact sequences
\begin{multline*}
0 \rightarrow
(\Fpbar\otimes_{\F_p} J_1(D p)[p]^{\ord})^0_{\frak n}\rightarrow
(\Fpbar\otimes_{\F_p} J_1(D p)[p]^{\ord})_{\frak n}  \rightarrow \\
(\Fpbar\otimes_{\F_p} J_1(D p)[p]^{\ord})_{\frak n}^E \rightarrow 0
\end{multline*}
and
\begin{multline*}
0 \rightarrow
(\Fpbar\otimes_{\F_p} J_1(D\ell p)[p]^{\ord})^0_{\frak m}\rightarrow
(\Fpbar\otimes_{\F_p} J_1(D\ell p)[p]^{\ord})_{\frak m}  \rightarrow  \\ 
(\Fpbar\otimes_{\F_p} J_1(D\ell p)[p]^{\ord})_{\frak m}^E \rightarrow 0.
\end{multline*}
The lower horizontal arrow $B_{\ell,\ell}^*$ of~(\ref{eqn:diagram})
is $G_{\Q}$-equivariant, and also
induces a morphism between these short-exact sequences.

If $B_{\ell,\ell}^* y$ is non-zero, then,
since $\rhobar$ is an irreducible $G_{\Q}$-module,
we may find $\sigma \in G_{\Q}$ such that $\sigma(B_{\ell,\ell}^* y)$
has non-zero image in
$(\Fpbar\otimes_{\F_p} J_1(D\ell p)[p]^{\ord})_{\frak m}^E$. 
Replacing $y$ by the image of $\sigma(y)$ in
$(\Fpbar\otimes_{\F_p} J_1(D p)[p]^{\ord})_{\frak n}^E,$ 
we are reduced to proving the following claim.

\

\noindent
{\bf Claim:} If $y \in 
(\Fpbar\otimes_{\F_p} J_1(D p)[p]^{\ord})_{\frak n}^E$ is such that
$B_{\ell,\ell}^*y \in
(\Fpbar\otimes_{\F_p} J_1(D\ell p)[p]^{\ord})_{\frak n}^E$
is annihilated by $\frak m$, then $B_{\ell,\ell}^* y$ vanishes.

\

Passing to the special fibre of the N\'eron model of each of $J_1(D p)$
and $J_1( D \ell p)$ over $\Z_p[\zeta_p]$ (where $\zeta_p$ denotes
a primitive $p$th root of unity), the discussion in the proof
of \cite[Thm.~2.1]{wiles} yields a commutative diagram in which
the horizontal arrows are isomorphisms:
$$
\xymatrix{
(\Fpbar\otimes_{\F_p} J_1(D p)[p]^{\ord})_{\frak n}^E
\ar[d]^-{B_{\ell,\ell}^*} \ar[r]^-{~} &
\left( H^0(\Sigma_1(D)^{\mu},\Omega^1) \oplus
H^0(\Sigma_1(D)^{\et}, \Omega^1)\right)^{\ord}_{\frak n}
 \ar[d]^-{B_{\ell,\ell}^*}\\
(\Fpbar\otimes_{\F_p} J_1(D \ell p)[p]^{\ord})_{\frak m}^E
\ar[r]^-{~} & 
\left( H^0(\Sigma_1(D\ell)^{\mu},\Omega^1) \oplus
H^0(\Sigma_1(D\ell)^{\et}, \Omega^1)\right)^{\ord}_{\frak m}.}
$$
(For any $M$ prime to $p$, we denote by $\Sigma_1(M)^{\mu}$ and
$\Sigma_1(M)^{\et}$
the base-change from $\F_p$ to $\Fpbar$
of the two components of the special fibre of the canonical model
of $X_1(M p)$ over $\Z_p[\zeta_p]$.)
Thus to prove the claim, it suffices to show
that if $z$ is a non-zero element of
$\left(H^0(\Sigma_1(D)^{\mu},\Omega^1) \oplus
H^0(\Sigma_1(D)^{\et}, \Omega^1)\right)^{\ord}_{\frak n}$
then $B_{\ell,\ell}^*$ cannot lie in the $\frak m$-eigenspace of 
$\left(H^0(\Sigma_1(D\ell)^{\mu},\Omega^1)\oplus
H^0(\Sigma_1(D\ell)^{\et}, \Omega^1)\right)^{\ord}_{\frak m}.$
This follows from the $q$-expansion principal; more precisely,
the $q$-expansion of $B_{\ell,\ell}^* z$ at the cusp $\infty$
involves only powers of $q^{\ell}$, and so cannot be an eigenform
for the full Hecke algebra $T_{D \ell}$.

The second case (when the mod $p$ diamond operators are trivial
modulo $\frak m$) is treated similarly, using the corresponding
results from the proof of \cite[Thm.~2.1]{wiles}.
\end{proof}

\section{Applications to the main conjecture and examples}
\label{sec:examples}

\subsection{The main conjecture}

Let $\O$ be the ring of integers of a finite extension $K$ of $\Qp$ and
let $f \in \O[[q]]$ be a $p$-ordinary eigenform such that the residual 
representation $\rhobar_f$ is irreducible.   Let 
$L_p^{\alg}(f,\omega^i) \in \Lambda_\O$ denote a generator of
the characteristic power series of the $\Lambda_{\O}$-dual of the Selmer group
$\Sel(\Q_\infty,\Afi)$ of Section~\ref{sec:selmer}.
Let $L_p^{\an}(f,\omega^i) \in \Lambda_{\O}$ 
denote the usual $p$-adic $L$-function of 
$f \otimes \omega^i$
(computed with respect to some canonical period). 
The main conjecture of Iwasawa theory in this context is the following;
it is independent of the particular choice of coefficient field $K$ over
which $f$ is defined.

\begin{conj}
There is a unit $u \in \Lambda_{\O}^{\times}$ such that
$$
L_p^{\alg}(f,\omega^i) \cdot u = L_p^{\an}(f,\omega^i).
$$
\end{conj}

We remark that the Selmer group in the above conjecture is not the 
Bloch--Kato Selmer group in which the local condition at $p$ is
defined via crystalline periods, but instead is
the Greenberg Selmer group in which the local
condition at $p$ is defined via the ordinary filtration as in
Section~\ref{sec:selmer}.  The former group is always contained in the latter
group, and the quotient is trivial unless the analytic $p$-adic
$L$-function has a trivial zero, in which case it has corank one
(see \cite[\pages 108-109]{Greenberg1}).
It is for this reason that in our statement of the main conjecture we do not 
need to consider separately the case of trivial zeroes. 

The following deep theorem of Kato \cite{KKT}
establishes one divisibility of the main conjecture.

\begin{thm}[Kato]
\label{thm:Kato}
There is a $u \in \Lambda_{\O} \otimes \Qp$ such that
$$
L_p^{\alg}(f,\omega^i) \cdot u = L_p^{\an}(f,\omega^i).
$$
In particular, to verify the main conjecture for $f \otimes \omega^{i}$ 
it suffices to check that
$$
\mu^{\alg}(f,\omega^i) = \mu^{\an}(f,\omega^i)  \text{~~and~~}
\lambda^{\alg}(f,\omega^i) = \lambda^{\an}(f,\omega^i).  
$$
\end{thm}

Our results in Sections~\ref{sec:algebraic} and
\ref{sec:analytic L-functions} yield the following result.

\begin{thm}
\label{thm:mainthm}
Let $k$ be a finite field of characteristic $p$ and
let $\rhobar : \GQ \to \GL_{2}(k)$ be an irreducible, modular, $p$-ordinary
and $p$-distinguished representation; as always we fix a choice of
$p$-stabilization.  Suppose that
$$
\mu^{\alg}(f_0,\omega^i) = \mu^{\an}(f_0,\omega^i) = 0
\text{~~and~~}
\lambda^{\alg}(f_0,\omega^i) = \lambda^{\an}(f_0,\omega^i)
$$
for some $f_0$ in the Hida family attached to $\rhobar$ and some $i$.  Then 
$$
\mu^{\alg}(f,\omega^i) = \mu^{\an}(f,\omega^i) = 0
\text{~~and~~}
\lambda^{\alg}(f,\omega^i) = \lambda^{\an}(f,\omega^i)
$$
for every $f$ in the Hida family attached to $\rhobar$.
\end{thm}

Before giving a proof, we first
state an immediate corollary of Theorem~\ref{thm:mainthm} and
Kato's result.

\begin{cor}
\label{cor:maincor}
Let $\rhobar$ be as above and suppose that
$\mu^{\alg}(\rhobar,\omega^i) = \mu^{\an}(\rhobar,\omega^i) = 0$
for some $i$.
If the main conjecture holds for $f_0 \otimes \omega^i$ for
one form $f_0$ in the Hida family of $\rhobar$, then the main conjecture
holds for $f \otimes \omega^i$ for every form $f$ in the Hida family
of $\rhobar$.
\end{cor}

Since every Hida family contains a form of weight two, this corollary
in particular reduces the main conjecture to the case of weight two and
to the conjecture on the vanishing of the $\mu$-invariants.

The proof of Theorem~\ref{thm:mainthm} is based on the following lemma which
relates the invariants $\delta_{v}(\a,\omega^{i})$ and
$e_{\ell}(\a,\omega^{i})$ of section \ref{sec:alginv} and
Theorem~\ref{thm:lambda}.

\begin{lemma}
\label{lemma:difference}
Let $\a_1$ and $\a_2$ be minimal primes of $\T_\Sigma$.  For any
prime $\ell \neq p$
$$
\sum_{v | \ell} \delta_v(\a_1,\omega^i) - \delta_v(\a_2,\omega^i) =
e_\ell(\a_2,\omega^i) - e_\ell(\a_1,\omega^i)
$$
where the sum is taken over all primes $v$ of $\Q_\infty$ over $\ell$.
\end{lemma}

\begin{proof}
To prove the lemma it suffices to see that for any minimal prime $\a$ of $\T_\Sigma$, the sum
$$
e_\ell(\a,\omega^i) + \sum_{v | \ell} \delta_v(\a,\omega^i)
$$
is independent of $\a$.
For this, fix a classical newform $f$ on the branch $\T(\a)$.
Consider first the group $H^{1}(\Qiv,\Afi)$.  
Since $H^{1}(\Qiv,\Afi)$ is divisible (as $G_{v}$ has
$p$-cohomological dimension one) we have
$$
\dim_k H^1(\Qiv,\Afi)[\pi] =
\lambda \left(
\chr_{\Lambda_\O} (  H^1(\Qiv,\Afi)^{\vee} )\right)
$$
where $M^\vee = \Hom(M,\Qp/\Zp)$.

By
\cite[Proposition 2.4]{GV} we have
$$
\chr_{\Lambda_\O} \left( \oplus_{v | \ell} H^1(\Qiv,\Afi)^{\vee} \right) =
E_\ell(f,\langle \ell \rangle_i^{-1}) \cdot \Lambda_\O.
$$
Since $E_{\ell}(f,\langle\ell \rangle_i^{-1})$ is simply
$E_\ell(\a,\langle\ell \rangle_i^{-1}) \text{~mod~} \heightone_f$, 
we conclude that
\begin{equation} \label{eq:churchill}
\sum_{v \mid \ell} \dim_k H^1(\Qiv,\Afi)[\pi] = e_{\ell}(\a,\omega^{i}).
\end{equation}

Consider now the exact sequence
$$
0 \to \Afi^{\Gv}/\pi \to H^{1}(\Qiv,\Ab \otimes \omega^{i}) \to
H^{1}(\Qiv,\Afi)[\pi] \to 0.
$$
Since the first term has $k$-dimension $\delta_{v}(\a,\omega^{i})$ and
the second term is certainly independent of $\a$, the lemma now follows
from (\ref{eq:churchill}).
\end{proof}

\begin{proof}[Proof of Theorem \ref{thm:mainthm}]
Let $f$ be a form in the Hida family of $\rhob$.
The vanishing of the $\mu$-invariants of $f$ is immediate from
that for $f_{0}$ and
Theorems \ref{thm:prealgebraic} and \ref{thm:mu vanishing}.

By Lemma~\ref{lemma:difference} and
Theorems \ref{thm:algebraic} and \ref{thm:lambda} we see also that
$$
\lambda^{\alg}(f,\omega^i) - \lambda^{\alg}(f_0,\omega^i)
= \lambda^{\an}(f,\omega^i) - \lambda^{\an}(f_0,\omega^i).
$$
Since we are assuming the main conjecture for $f_0 \otimes \omega^i$, 
we have that
$$
\lambda^{\alg}(f_0,\omega^i) = \lambda^{\an}(f_0,\omega^i),
$$
so that it follows that
$$
\lambda^{\alg}(f,\omega^i) = \lambda^{\an}(f,\omega^i)
$$
as desired.
\end{proof}

\subsection{Raising the level}

We conclude our general discussion with some results
on branches of Hida families.
It is well-known (see \cite{diamondtaylor1,diamondtaylor2}) which
levels can occur among forms in the Hida family of $\rhob$.
In the next proposition we list the cases in which the invariant
$\delta_{v}(\cdot,\omega^{i})$ increases.

\begin{prop}
\label{prop:lambdachange}
Let $\rhob : \GQ \to \GL_{2}(k)$ be an
irreducible, $p$-ordinary and $p$-distinguished modular Galois
representation; assume further that $\rhob$ is ramified at $p$.
Let $\Sigma$ be some finite set of primes not containing $p$ and
let $\ell \neq p$ be a prime at which $\rhobar$
is unramified.  Set $\Sigma' = \Sigma \cup \{\ell\}$ and let
$\T_\Sigma$, $\T_{\Sigma'}$ be the Hida algebras associated to $\rhob$.
Let $a_{\ell} = \Trace \rhob(\Frob_{\ell})$ and
$c_{\ell} = \det \rhob(\Frob_{\ell})$, both viewed as elements of $k$.
Then for every branch $\T(\a)$ of $\T_\Sigma$:
\begin{enumerate}
\item If there is an $i$ such that
$a_{\ell} = \ell^{-i} + \ell^{1-i}$ and $c_{\ell} = \ell^{1-2i}$,
then there is a branch $\T(\b)$ of $\T_{\Sigma'}$ such that
$$N(\b) = N(\a)\ell \text{~~and~~} \delta_{v}(\b,\omega^{i}) = 1$$
for all places $v$ dividing $\ell$.

\item If $\ell \equiv 1 \pmod{p}$ and $a_{\ell} = c_{\ell} + 1$,
then there is a branch $\T(\b)$ of $\T_{\Sigma'}$ 
(distinct from the branch provided by (1)) such that
$$N(\b) = N(\a)\ell \text{~~and~~} \delta_{v}(\b,\omega^{i}) = 1$$
for all places $v$ dividing $\ell$ and all $i$.

\item If $\ell \equiv 1 \pmod{p}$, $a_{\ell} =2$, and $c_{\ell} =1$,
then there are distinct
branches $\T(\b_{1})$ and $\T(\b_{2})$ of $\T_{\Sigma'}$ such that
$$N(\b_1) = N(\b_2) = N(\a)\ell^2 \text{~~and~~} 
\delta_{v}(\b_{1},\omega^{i}) = \delta_{v}(\b_{2},\omega^{i}) = 2$$
for all places $v$ dividing $\ell$ and all $i$.

\item If $\ell \equiv -1 \pmod{p}$, $a_{\ell} = 0$, and $c_{\ell} = -1$,
then there is a branch $\T(\b)$ of $\T_{\Sigma'}$ such that
$$N(\b) = N(\a)\ell^2 \text{~~and~~} \delta_{v}(\b,\omega^{i}) = 1$$
for all places $v$ dividing $\ell$ and all $i$.
\end{enumerate}
\end{prop}

Note in particular that if (3) holds, then (1) and (2) also hold,
so that the proposition
exhibits four distinct non-minimal branches of $\T_\Sigma'$ for each
branch of $\T_\Sigma$.

\begin{proof}
If (1) holds, then one sees easily that the eigenvalues of $\Frob_{\ell}$
on $\rhob$ equal $\ell^{1-i}$ and $\ell^{-i}$.  It follows that there
exists a ramified representation
$$\tau_{\ell} : G_{\ell} \to \GL_{2}(\Zp)$$
such that $\bar{\tau}_{\ell} \otimes k \cong \rhob|_{G_{\ell}}$ and
$$\tau_{\ell} = \left( \begin{array}{cc}
\omega^{1-i} & * \\ 0 & \omega^{-i} \end{array} \right).$$
By \cite[Theorem 1]{diamondtaylor2} 
(which applies since $\rhob$ is modular of weight $2$ by
\cite[Theorem 1.1]{diamond1}) there exists a newform $g$
of tame level $N\ell$ in the Hida family of $\rhob$ such that
$\rho_{g}|_{I_{\ell}} \cong \tau_{\ell}|_{I_{\ell}}$.
Since $\rhob \otimes \omega^{i}$ is unramified at $\ell$ with
Frobenius eigenvalues $\ell$ and $1$, one computes easily that
$$\delta_{v}(A_{g},\omega^{i}) = 1$$
for any place $v$ dividing $\ell$.
In particular, the branch $\a_{1}$ on which $g$ lies satisfies the
requirements of the proposition.

The proofs for (2) (principal series), (3) (either special or
principal series) and (4) (supercuspidal) are entirely similar;
the conditions on $\ell$, $a_{\ell}$, and $c_{\ell}$ are simply those
obtained by combining the conditions of
\cite[p.\ 435]{diamondtaylor1} for ramified 
lifts of the appropriate type to exist
with the requirement that $\rhob(\Frob_{\ell})$ has a trivial eigenvalue.
Note also that $\omega|_{G_{\ell}}$ is trivial (resp.\
quadratic) in (2) and (3) (resp.\ (4)), so that the choice of $i$
is irrelevant.  Finally, the branches in (1) and (2) (resp.\ in (3)) are
distinct since one is special at $\ell$ and the other is principal series
at $\ell$.
\end{proof}

It would be interesting to determine when the branches of the
preceding proposition intersect.
The following proposition uses $\Lambda$-adic level raising to give
some insight into this question.

\begin{prop}
\label{prop:crossing}
We maintain the hypotheses of the previous proposition. 
Let $\a$ be a minimal prime of $\T_\Sigma$ such that 
$$
(T_\ell \text{~mod~} \a')^2 - \ell^{-2} \langle \ell \rangle (\ell+1)^2
$$
is not a unit in $\T(\a)^\circ$.  (Recall that $\a'$ is a minimal prime of
$\T^{\new}_{N(\a)}$ sitting over $\a$ and that $\T(\a)^\circ \cong \T^{\new}_{N(\a)}/\a'$.)
Then there exists some minimal prime $\b$ of $\T_{\Sigma'}$ with
$N(\b) = N(\a) \ell$ and a height one prime $\heightone$ of $\T_{\Sigma'}$
such that $\T(\a)$ and $\T(\b)$ cross at $\heightone$.  
(That is, both $\b$ and the pre-image of $\a$ under the
natural map $\T_{\Sigma'} \to \T_\Sigma$ contain $\heightone$.)
\end{prop}

\begin{proof}
Attached to the minimal prime $\a$ we have the $\Lambda$-adic modular form
$$
f(\frak a,q) = \sum_{n \geq 1} (T_n  \bmod \frak a') q^n.
$$
By \cite[Theorem 6C]{diamond3}, there exists  a finite extension
$R$ of $\T(\a)^\circ$,
a $\Lambda$-adic modular form $g(q) = \sum_{n \geq 1} b_n q^n \in R[[q]]$ new
of level $N(\a)\ell$,
and a height one prime $\heightone'$ of $R$ such that
\begin{equation}
\label{eqn:cong}
(T_n \bmod \a') \equiv b_n \pmod{\heightone'}.
\end{equation}
for each $n$ relatively prime to $N(\a)\ell$.

The $\Lambda$-adic form $g(q)$ corresponds to some minimal prime ideal $\b$ of 
$\T_\Sigma$.  Indeed, the map $\Hida \rightarrow R$ that sends $T_n$ to
$b_n$ factors through $\T^{\new}_{N(\a)\ell}$.  Let 
$\b' \subseteq \T^{\new}_{N(\a)\ell}$
denote the kernel of this map; it is necessarily a minimal prime since
both the source and the target are finite extensions of $\Lambda$. 
The preimage of $\b'$ in $\prod_{M | N(\Sigma')} \T^{\new}_{M}$ is
a minimal prime 
and thus, by Proposition \ref{prop:embedding}, corresponds 
to some minimal prime of $\T'_{N(\Sigma')}$. 
Since the residual representation attached to $g(q)$ equals $\rhobar$, this minimal
prime yields a minimal prime of $\T_{\Sigma'}$ which we denote by $\b$.
Note that by construction $N(\b) = N(\a)\ell$.

We thus have a map $\T_{\Sigma'} \rightarrow R$ 
(sending $T_{n}$ to $b_{n}$) with kernel $\b$.
If let $\tilde{\a}$ denotes the pre-image of $\a$ in $\T_{\Sigma'}$, then
we also have
a map $\T_{\Sigma'} \rightarrow R$ with kernel $\tilde{\a}$ (given by reducing
mod $\tilde{\a}$ and then embedding the image $\T(\a)$ into $R$).  By (\ref{eqn:cong}), 
the reduction of these two maps modulo $\heightone'$ are the same.  
If we let $\heightone$
denote the kernel of either of these maps modulo $\heightone'$,
then $\heightone$ is
a height one prime contained in both of 
$\b$ and $\tilde{\a}$, as desired.
\end{proof}

The ideal $\heightone$ of the previous proposition could potentially be
a prime ideal lying over the principal ideal $(p)$.
It would be interesting to determine whether or not these branches actually meet at
a prime ideal of residue characteristic zero.

\subsection{Examples}

\begin{example}
Set $p=11$ and let $f$ denote the weight $2$ newform associated to
the elliptic curve $X_{0}(11)$.  It follows from the fact that $f$ is the
only newform of weight $2$ and level dividing $11$ that 
$\T := \T_{\Sigma}(\rhobar_{f}) \cong \Lambda$
for $\Sigma = \emptyset$.  Thus for $k \geq 2$, there is a unique
newform $f_k$ of weight $k$ and 
level dividing $11$ that is congruent to $f$ modulo $11$.
For example, $f_{12}$ is the $11$-ordinary, $11$-stabilized oldform
of level $11$ attached to the Ramanujan $\Delta$-function.

The $p$-adic $L$-function of this family was studied in detail in
\cite{GS}.  Also, in \cite{Greenberg3}, the congruence between $X_0(11)$
and $\Delta$ was exploited to gain information about the Selmer groups 
of both forms.  We review below what is known in these examples and then
go on to study the other branches of this Hida family with tame conductor
greater than one.

We first verify the main conjecture for $f$.  Since $X_{0}(11)$
has split multiplicative reduction at $p=11$, the $p$-adic $L$-function
$L_{p}^{\an}(f,T)$ has a trivial zero at $T=0$.
By the Mazur, Tate, and Teitelbaum conjecture (proved by Greenberg
and Stevens \cite{GS}) we have an explicit formula for the derivative
of $L_{p}^{\an}(f,T)$ at $0$; that is,
$$
\frac{d}{dT} L_{p}^{\an}(f,T) \Big|_{T=0} = 
\frac{\L_{p}(f)}{\log_p(\gamma)} \cdot \frac{L(f,1)}{\Omega_f}.$$
(Here, $\gamma$ is a topological generator of $1+p\Zp$ that is
implicitly chosen by writing the $p$-adic $L$-function in the
$T$-variable; the appearance of the extra factor of $\log_p(\gamma)$
above is accounted for by the fact that we have written this formula
in the $T$-variable rather than the more standard $s$-variable.)
One checks that $\L_{p}(f)$ and $\log_p(\gamma)$ are exactly divisible by
$11$ while $L(f,1)/\Omega_f$ is an $11$-adic unit.  
Thus $L_{p}^{\an}(f,T)$ is a unit multiple of $T$ so that
$\lambda^{\an}(f)=1$ and $\mu^{\an}(f) = 0$.

Since $L(f,1)/\Omega_f$ is an $11$-adic unit, by
an Euler system argument of Kolyvagin,
the classical $11$-adic Selmer group $\Sel_{p}(X_{0}(11)/\Q)$
of the elliptic curve
$X_{0}(11)$ vanishes.  It thus follows from a refined control
theorem of Greenberg (\cite[Proposition 3.7]{Greenberg2})
that $\Sel_{p}(X_{0}(11)/\Qi)$ vanishes as well.
The latter group agrees with the Bloch--Kato Selmer group of $f$; as
$f$ has a trivial zero, it follows that the Greenberg Selmer group
$\Sel_{p}(\Qi,A_{f})$ is simply $\Q_{p}/\Z_{p}$ with trivial
Galois action.  
(This example is worked out in detail in \cite[Example 3]{Greenberg3}.)
Thus $L_{p}^{\alg}(f,T)$ is a unit multiple of $T$, so that
$\lambda^{\alg}(f) = 1$ and $\mu^{\alg}(f) = 0$.    
In particular, this verifies the main conjecture for $X_0(11)$ at $p=11$.

Corollary \ref{cor:maincor} now shows that
the main conjecture holds at $p=11$ for each form $f_{k}$ with $k \geq 2$,
and Theorem \ref{thm:mainthm} shows that
$$
\lambda^{\an}(f_k) = \lambda^{\alg}(f_k) = 1 \text{~for~} k \geq 2.
$$

Note that the $p$-adic $L$-function corresponding to $X_0(11)$ is the only
$p$-adic $L$-function in the Hida family with a trivial zero.  Nonetheless, 
the existence of this trivial zero forces every $p$-adic $L$-function in the
family to have at least one zero.

For $k \equiv 2 \pmod{10}$, the unique zero of $L_{p}^{\an}(f_k)$
can be explained by the functional equation of the $p$-adic $L$-function.  
Namely, the sign of this functional equation is $-1$ for such $k$ and thus
$L_{p}^{\an}(f_k)$ vanishes at the character that sends $x$ to 
$x^{\frac{k-2}{2}}$.  
Using this observation, we can determine
the two-variable $p$-adic $L$-function attached to $\rhobar_{f}$ (which
we denote by $L(\rhobar_{f})$).
By Theorem \ref{thm:lambda}, we have that $\lambda(L(\rhobar_{f})) = 1$  
and hence, under the identification of (\ref{eqn:second choice}),
$$L(\rhobar_{f}) = (T - a_0) \cdot u(T)$$ 
with $a_0 \in \T$ and $u(T) \in \T[[T]]$ a unit.
Since $L(\rhobar_{f})$ specializes at weight $k$ to $L_{p}(f_k)$ (Proposition
\ref{prop:two vs one}), we know that $a_0$ specializes at weight $k$
to the unique zero of $L_{p}(f_k)$.  In particular, for $k \equiv 2 \pmod{10}$,
we have that $a_0$ specializes to $\gamma^{\frac{k-2}{2}}$.
Thus, $a_0 = \gamma^{-1} \langle \gamma \rangle^{\frac{1}{2}}$ since this
element has the correct specialization for infinitely many $k$.  

So far we have only applied our results to minimal lifts of $\rhobar_{f}$; 
however,
Corollary \ref{cor:maincor} applies to every modular form lifting this representation.
For instance, 
we can use Proposition \ref{prop:lambdachange} to explicitly find primes $\ell$ to add
to the level that will produce modular forms that have higher $\lambda$-invariants and
for which we still know the main conjecture.

The first prime that satisfies condition (3) (and thus (1) and (2)) of Proposition 
\ref{prop:lambdachange} is $\ell = 1321$.  If we set $\Sigma = \{ 1321 \}$
then, $\T_\Sigma$ will have four distinct branches $\T(\a_1), \dots, \T(\a_4)$
with $\lambda(\a_1) = \lambda(\a_2) = 2$ and $\lambda(\a_3) = \lambda(\a_4) = 3$.  
(Here and in what follows, we do not distinguish between the analytic and algebraic $\lambda$-invariants
in cases where the main conjecture is known to be true.)
Note
that $1321$ is inert in the cyclotomic $\Z_{11}$-extension of $\Q$ and so there is
a unique prime $v$ of $\Qinf$ sitting over $1321$. 

The prime $\ell' = 2113$ is also inert in $\Qi/\Q$ and
satisfies condition (3) of Proposition \ref{prop:lambdachange}.
If $\Sigma' = \{ 1321, 2113 \}$, we then have that $\T_{\Sigma'}$ possesses
$20$ additional branches with $\lambda$-invariants ranging between $3$ and $5$. 
\end{example}

\begin{example}
We conclude by considering an example of
Greenberg and Vatsal in \cite{GV}.
Consider the elliptic curves
$$
E_1 : y^2 = x^3 + x - 10 \hspace{2.0cm} E_2 : y^{2} = x^3 - 584 x + 5444,
$$
of conductors $52$ and $364$ respectively, which are ordinary and congruent
mod $5$.  One computes that $\lambda^{\an}(E_{1})=\mu^{\an}(E_{1})=0$ so
that Kato's divisibility (Theorem \ref{thm:Kato})
yields the main conjecture for $E_{1}$ at $p=5$.  The main
conjecture for $E_{2}$ at $p=5$ follows and in this case
we have $\lambda(E_{2})=5$.  (This example is discussed in detail in \cite[\pages 22,44]{GV}.)

We now examine this congruence from the point of view of Hida theory.  
Let $f_{i}$ denote the newform of weight $2$ associated to $E_{i}$.
One checks that $f_{1}$ is not congruent modulo $5$
to any other modular forms of
weight two and level dividing $2^{2} \cdot 5 \cdot 13^{2}$
(using \cite{Stein}, for example).
It follows that we have $\T_{\Sigma} \cong \Lambda$ for
$\Sigma = \{2,13\}$.  
For consistency of notation, let us denote the unique irreducible
component of this space as $\T(\a_1)$.  
Moreover, since $\mu(f_1) = \lambda(f_1) = 0$,
the two-variable $p$-adic $L$-function attached to $\T_\Sigma$ is a unit
and $\lambda(\a_1) = 0$.

The prime $\ell = 7$ satisfies condition (1) of Proposition \ref{prop:lambdachange} (with $i=0$).  
Thus, if $\Sigma' = \{2, 7, 13 \}$, 
then $\T_{\Sigma'}$ contains an irreducible component $\T(\a_2)$
not contained in $\T_\Sigma$.  Moreover, since 
$f_2$ is the unique normalized newform congruent
modulo $5$ to $f_1$ with level 
dividing $2^2 \cdot 5 \cdot 7^2 \cdot 13^2$, it follows that $\T_{\Sigma'}$
has rank 2 over $\Lambda$.  
Hence, $\T(\a_2)$ is the only branch of $\T_{\Sigma'}$ not coming from
$\T_\Sigma$.    
By Theorems \ref{thm:algebraic} and \ref{thm:lambda}, we have that
$\lambda(\a_2) = 5$.
Since $f_2$ sits on the branch $\T(\a_2)$, it follows that $\lambda(f_2) = 5$.

We close with some questions.
Proposition \ref{prop:crossing} establishes that the branches
$\T(\a_1)$ and $\T(\a_2)$ must cross at some (non-classical)
height one prime, but
does not exclude the possibility that they cross at the prime $(p)$.  
Do they in fact cross
at a prime of residue characteristic zero, and if so 
could one compute $p$-adic approximations
of this prime?  How many such crossing points do these two branches share?
It appears at present 
that little is known about the shape of these Hida families 
when multiple branches appear (even in any particular case).

\end{example}



\providecommand{\bysame}{\leavevmode\hbox to3em{\hrulefill}\thinspace}

\end{document}